\documentclass[twoside]{article}

\usepackage{graphicx}
\usepackage{amsmath}
\usepackage{amssymb}
\usepackage{hhline}
\usepackage{epsfig}

\newtheorem{theorem}{Theorem}

\newtheorem{theoremb}{Theorem}
\newtheorem{theoremc}{Theorem}
\newtheorem{theoremd}{Theorem}
\newtheorem{cor}[theoremd]{Corollary}
\newtheorem{dfn}[theoremb]{Definition}

\newtheorem{lem}[theorem]{Lemma}
\newtheorem{prop}[theorem]{Proposition}
\newtheorem{rk}[theoremc]{Remark}

\newenvironment{proof}[1][Proof]{\textbf{#1. }}{\qed}
\newenvironment{Proof}[1]{\textbf{#1. }}

\newcommand\bib[1]{\bibitem[#1]{#1}}

\newcommand\abz{\hspace{12pt}}
\newcommand\qed{\phantom{\underline{y}}\hfill\hfill$\square$}

\newcommand\1{{\bf 1}}

\newcommand\h{h_{\text{\rm top}}}
\newcommand\hm{H_{\text{\rm mult}}}
\newcommand\hs{H_{\text{\rm sing}}}
\newcommand\e{\epsilon}
\renewcommand\l{\lambda}
\newcommand\z{\sigma}
\newcommand\op[1]{\mathop{\rm #1}\nolimits}
\newcommand\po{$\!\!\!{\text{\bf.}}$ }
\newcommand\ls{\mathop{\overline\lim}\limits_{n \rightarrow \infty}}
\newcommand\N{{\mathbb N}}
\newcommand\Z{{\mathbb Z}}
\newcommand\R{{\mathbb R}}
\newcommand\D{{\mathcal D}}

\newcommand\Y{{\mathcal Y}}
\newcommand\T{{\mathcal T}}
\newcommand\vp{{\varphi}}
\newcommand\hps{\hskip-16pt . \hskip2pt}
\newcommand\hpss{\hskip-13.5pt . \hskip2pt}
\newcommand{\weg}[1]{}

\makeatletter
\renewcommand{\@oddhead}{\hfil Dynamics and entropy of SOC\hfil}
\renewcommand{\@evenhead}{\hfil Boris Kruglikov, Martin Rypdal \hfil}
\makeatother

\begin{document}


\title{Dynamics and entropy in the Zhang model of Self-Organized Criticality}
\author{B. Kruglikov  \& M. Rypdal \\ ~ \\
{\small Institute of Mathematics and Statistics}\\
{\small University of Troms\o, N-9037 Troms\o, Norway}\\
{\small Boris.Kruglikov@matnat.uit.no;
Martin.Rypdal@matnat.uit.no} }
\date{}
\maketitle

\begin{abstract}
We give a detailed study of dynamical properties of the Zhang
model, including evaluation of topological entropy and estimates
for the Lyapunov exponents and the dimension of the attractor. In
the thermodynamic limit the entropy goes to zero and the Lyapunov
spectrum collapses.\!\!
\end{abstract}


\section*{Introduction}

In 1987 the concept of Self-Organized Criticality (SOC) was introduced by
Bak, Tang and Wiesenfeld \cite{BTW}. The attempt was to give an
explanation of the omnipresence of fractal structures and power-law
statistics in nature, and the claim was that certain physical
systems can self-organize into stationary states, reminiscent of equilibrium
system at the critical point, in the sense that one has scale invariance
and long range correlations in space and time.

SOC is proposed as an explanation for variety of phenomena in
nature, such as earthquakes, forest fires, stock markets and
biological evolution \cite{J}. However, most work has been devoted
to the study of idealized ''sandpile-like'' computer models, such
as the sandpile model \cite{BTW}, the abelian sandpile \cite{DR}
and the Zhang model \cite{Z} that, one believes, exhibit SOC in
the thermodynamic limit. Despite this effort, a satisfactory
understanding of the model is not yet achieved. Through numerical
investigation it was observed that in the thermodynamic limit,
observables have power-law distributions. More precisely, the
probability distribution of an observable $s$ has the form $P(s)
\sim 1/s^{\tau_s}$ in the thermodynamic limit. There is no widely
agreed upon method for computing the SOC-exponents $\tau_s$
numerically and, due to the incomplete understanding of the
dynamics of the models and lack of a formal treatment of the
thermodynamic limit, it is difficult to properly explain the
observed behavior. Hence it is not clear what the SOC-exponents
really tell us about the dynamics of the SOC models.

 \subsection{\hpss Discussion of the Zhang model}
In a series of papers by  Cessac, Blanchard and Kr\"uger
\cite{BCK} it was proposed that deeper understanding of SOC models
can be achieved by studying the models in the framework of
dynamical system theory. They showed how a particular model, the
Zhang model, could be formulated as a dynamical system of
skew-product type with singularities, where the randomness of the
external driving is described by a Bernoulli shift, and the
threshold relaxation dynamics is given by piecewise affine maps.

In this paper we present a detailed study of the dynamical system
defined in \cite{BCK}. We prove several basic properties, some of
which are already stated in \cite{BCK}, before discussing
fundamental dynamical properties. Depending on the parameters of
the model, we can observe fundamentally different types of
behavior.

For low values of the threshold energy (critical energy), the
dynamics can be relatively simple, since the singularities only
effect the dynamics in a finite number of time-steps. In such
situations we say that singularities are removable, and we show
that the system permits symbolic coding. We give examples of how
symbolic coding provides a complete description of the
dynamics as a topological Markov chain. Hence the dynamics is
chaotic, but the essential dynamical invariants are all inherited
from the Bernoulli shift factor. Moreover we can identify the
physical invariant measure, and hence understanding of the
statistical properties is reduced to the theory of Markov chains.

As we increase the critical energy the role of the singularities
becomes essential. Techniques based on codings are no longer
applicable and a very interesting dynamics emerges. The
dimensional characteristics of the attractor are also sensitive to
the parameters of the system, as we show by generalizing the Moran
formula for the iterated function system (IFS). In addition, we
observe the situation, when the dimension of the IFS-attractor
increases to the maximum, while the support of the SRB-measure
remains fractal.

To measure the complexity of the dynamics we study entropy and
Lyapunov exponents. We show that the system is hyperbolic, with
one positive exponent originating in the Bernoulli shift. However,
due to the presence of singularities the Ruelle inequality and the
Pesin formula are not directly applicable. We show that the metric
entropy of any SRB-measure equals the topological entropy almost
surely, and we evaluate the latter generalizing the technique
developed by Buzzi \cite{B1,B2}. The result is that the Pesin
formula and the variational principle hold a posteriori.

To give a satisfactory physical interpretation of the dynamics we
rescale time to prevent infinitely slow driving of the system. We
prove that for this physical system, the Lyapunov spectrum
collapses completely and that the entropy goes to zero in the
thermodynamic limit. This implies that the expanding (chaotic)
properties are lost, so that we may expect power-laws statistics
and long range correlation effects.

The statistical properties we obtain hold for any SRB-measure,
because the most input comes from the Bernoulli shifts. The
existence of SRB-measures is in fact still an open problem. From
the general theory of dynamical systems with singularities
\cite{KS,P1,ST}, we can give conditions that are sufficient for
the existence of SRB-measures, but it is not known if these
conditions hold for the majority of parameters. We expect this to
be true (it was also conjectured in \cite{BCK}) and derive some
statistical corollaries.

Apart from the physical importance of the Zhang model, it is
interesting from a mathematical point of view. It can be described
as a piecewise affine hyperbolic map of the form
 $$
F:\Sigma_N^+\times M\to\Sigma_N^+\times M,\quad
((t_0t_1t_2\dots),x)\mapsto((t_1t_2t_3\dots),f_{t_0}(x)),
 $$
where $\Sigma_N^+$ is the set of right infinite sequences from a
finite alphabet and $\{f_i\}$ a collection of piece-wise affine
non-expanding maps of $M$ to itself. Previously piecewise affine
expanding maps and piecewise isometries have been studied, but the
contracting property of the relaxation dynamics gives rise to some
difficulties. Therefore several methods are developed in this
paper, which hold far beyond the framework of the Zhang model.

 \subsection{\hpss Structure of the paper}
In Section \ref{sec_1} we describe the model and derive bounds on
the size and duration of avalanches. This enables us to
use the Poincar\'e return to reformulate the systems in a
skew-product form. Then we study the contraction property to
conclude hyperbolicity of the model (Theorem \ref{5}) and
describe, when degenerations occur (Theorem \ref{th_invert}; the
original Zhang setting $\e=0$ is not the only possibility). In
Section \ref{sec_2} we introduce the concept of removability of
singularities, which appears in the coding approach for the study
of the model.

Section \ref{sec_3} is devoted to the study of measure entropy and
Lyapunov spectrum. We prove in Theorem \ref{17} that the entropy
of an SRB-measure is always maximal. Section \ref{sec_4} concerns
the topological entropy. We evaluate it for the most parameter
values (Theorems \ref{140} and \ref{a.s.0}). We also discuss
nearly-Zhang models and show that the dynamical quantities do not
change. This is natural from the physical perspective, because SOC
should not be obtained through a fine tuning of parameters.

Section \ref{sec_5} briefly describes the dimension issues of the
model (Theorem \ref{theo-5} gives the asymptotic values), which
enters into all inter-relations involving entropy and
characteristic exponents. We demonstrate how the fractality occurs
in IFS-context, noting the difference due to singularities and
overlaps. In Section \ref{examples} we illustrate the most
important effects in the model by examples.

In Section \ref{sec_7} we discuss the thermodynamic limit and
attempt to explain appearance of the power-law statistics via a
reparametrization. Conclusion contains the physical implications
of the current investigation.

In Appendices \ref{app_A} and \ref{app_B} we provide bounds for
the entropy and the dimension, which are new in the presence of
singularities, overlaps and degenerations (this was designed for
an application to the Zhang model). The results are of interest in
its own and can be read independently.


\section{\hps Basic properties of the Zhang model}\label{sec_1}

In the Zhang model each site on the lattice is associated with a
non-negative real number, which we call the energy of the site.
The collection of energies is called an energy configuration, and
can be represented as a point in $N$-dimensional space, where $N$
is the number of sites in the lattice. If a configuration is
unstable, the overcritical sites will lose some of their energy to
their nearest neighbors, resulting in a new energy configuration.
This transformation on $\R^N$ is denoted by $f$. If a
configuration is stable, a site is chosen at random and an energy
quantum $\delta=1$ is added to this site. In \cite{BCK} it was
shown how the relaxation and random excitation can be formulated
as a map of skew-product type on an extended phase-space. This
extended phase-space has the configuration space as one factor,
and the set of all possible sequences of excitations as the other
factor. In \cite{BCK} it was also shown how one can reformulate
the dynamical system by considering the return maps to the set of
stable configurations. This gives a simplification, in the sense
that each avalanche is associated with an affine transformation.
The set of stable configurations is partitioned into domains,
where each domain corresponds to an avalanche.


\subsection{\hpss Relaxation}

Take $d,L \in {\mathbb N}$ and let $\Lambda \subset {\mathbb Z}^d$
be the cube $[1,L]^d$ of cardinality $N:=L^d=|\Lambda|$. Let
$\phi:\Lambda \rightarrow \Lambda':= \{1,\dots,N\}$ be a
bijection. We define a metric $d_{\Lambda}$ on $\Lambda$ by
$$
d_\Lambda ({\bf k},{\bf l})=\sum_{1 \leq n \leq d} |k_n-l_n|\,,
$$
and let $d_{\Lambda'}:=\phi_{*} d_\Lambda$. In the following we
omit primes when it is clear from the context that we are
considering the metric space $(\Lambda', d_{\Lambda'})$. Elements
of $\Lambda$ will be called sites. We say that sites $i$ and $j$
are nearest neighbors if $d_{\Lambda}(i,j)=1$. The boundary
$\partial \Lambda$ is defined as those sites $i \in \Lambda$ that
have less than $2d$ nearest neighbors.

Fix parameters $E_c>0$ and $\epsilon \in  [0,1)$ and define
$f:{\mathbb R}^N_{\geq 0} \rightarrow {\mathbb R}^N_{\geq 0}$ by
$$
f(x)_i=x_i-\theta(x_i-E_c)(1-\epsilon)x_i+\frac{1-\epsilon}{2d}
\sum_{d_\Lambda (i,j)=1} \theta(x_{j}-E_c)x_{j}\,,
$$
where
$$
\theta(a)= \begin{cases} 1 \text{ if } a>0 \\ 0 \text{ if } a \leq
0\end{cases}\,.
$$
Let $\|x\|_1=\sum_{i=1}^N |x_i|$ be the 1-norm on $\R^N$.
\begin{prop}\po \label{1}
For all $x \in \R_{\geq 0}^N$ we have
$$
\frac{1+\epsilon}{2} \|x\|_1\leq \|f(x)\|_1 \leq \|x\|_1\,,
$$
and $\|f(x)\|_1 = \|x\|_1$ if and only if $x_i \leq E_c$ for all
$i \in \partial \Lambda$. If there is $i \in \partial \Lambda$
such that $x_i>E_c$ then
$$
\|f(x)\|_1 \leq \|x\|_1 - \frac{1-\epsilon}{2d}E_c\,.
$$
\end{prop}
\begin{proof}
Let $\{x_{i_k} \}_{k=1}^m$ be the entries of the vector $x$ that
are greater than $E_c$. Let $n_{i_k}$ be the number of nearest
neighbors of $x_{i_k}$. Then
\begin{eqnarray*}
\|f(x)\|_1 &=& \sum_{i \in \Lambda} f(x)_i =\sum_{i \in \Lambda}
x_i-(1-\epsilon)\sum_{k=1}^m x_{i_k}+\frac{1-\epsilon}{2d}
\sum_{k=1}^m n_{i_k} x_{i_k} \\
&=&\sum_{i \in \Lambda} x_i -(1-\epsilon) \sum_{k=1}^m
(1-\frac{n_{i_k}}{2d})x_{i_k}\,.
\end{eqnarray*}
The statement follows from the fact that we always have $d \leq
n_{i_k} \leq 2d$, and $n_{i_k} = 2d$ if and only if $x_{i_k} \not
\in \partial \Lambda$.
\end{proof}
\\

We say that a site $i \in \Lambda$ of the configuration $x$ is
relaxed if $x_i \leq E_c$, and excited if $x_i>E_c$. A
configuration $x$ is called stable if all sites are relaxed. The
set of stable configurations is $M:=[0,E_c]^N$. For each
configuration $x$ we define $m(x)=\min \{n \geq 0\,|\,f^n(x) \in
M\}$.
\begin{prop}\po \label{2}
For all $x \in \R_{\geq 0}^N$ we have:
$$m(x) \leq
\frac{2dN}{1-\epsilon}\frac{\|x\|_1}{E_c}
\Big{(}\frac{2d}{1-\epsilon}+1\Big{)}^{\text{\em
diam}(\Lambda)/2}\,.
$$
\end{prop}

We need the following lemma ($[ \cdot ]$ denotes integer part):

\begin{lem}\po \label{3}
For $x \in \mathbb{R}^N_{\geq 0}$ and $n \in \mathbb{N}$ let
$\alpha_i(n,x)$ be the cardinality of the set $\{ l \leq
n\,|\,(f^l x)_i>E_c\}$.  Let $\gamma=[2d/(1-\epsilon)]+1$. If
$d_\Lambda(i,j)=1$, then $\alpha_j(n,x) \geq
[\alpha_i(n,x)/\gamma]$.
\end{lem}

\noindent
\begin{proof}
There is a finite increasing sequence $\{m_k\}$, such that
$(f^{m_k}(x))_i>E_c$. We claim that on each interval $(m_k,
m_{k+\gamma}]$ there is a number $m$ such that $(f^m(x))_j>
E_c$. In fact, in the opposite case
$$
\big{(}f^{1+m_{k+\gamma-1}}(x) \big{)}_j \geq \gamma
\frac{1-\epsilon}{2d} E_c>E_c\,.
$$
Since $[0,\alpha_i(n,x)]$ contains $\beta=[\alpha_i(n,x)/\gamma]$
disjoined such intervals, we get $\alpha_j(n,x) \geq \beta$. Thus
$\alpha_j(n,x) \geq [ \alpha_i(n,x)/\gamma ]$.
\end{proof}
\\

\noindent {\bf Proof of Proposition \ref{2}.} By applying
inductively Lemma \ref{3} we get:
$$
\alpha_j(n,x) \geq
\Big{[}\frac{\alpha_i(n,x)}{\gamma^{d_\Lambda(i,j)}}\Big{]}
$$
In fact, if $j=j_0,j_1,\dots j_k=i$ is a path with
$d_\Lambda(j_s,j_{s+1})=1$ and
$$
\Big{[} \frac{\alpha_i(n,x)}{\gamma^k} \Big{]}=t\,,
$$
then $\alpha_{j_k}(n,x) \geq t\gamma^k,\alpha_{j_{k-1}}(n,x) \geq
t\gamma^{k-1},\dots, \alpha_{j_0}(n,x) \geq t$. By Proposition
\ref{1}
$$
\alpha_j(n,x) \leq \frac{2d}{1-\epsilon} \frac{\|x\|_1}{E_c}
$$
for $j \in \partial \Lambda$, so
$$
\alpha_i(n,x) \leq \alpha(x):= \frac{2d}{1-\epsilon}
\frac{\|x\|_1}{E_c} \gamma^{\text{diam}(\Lambda)/2}\,.
$$
for all $i \in \Lambda$ and all $n \in \mathbb{N}$. Suppose
$f^m(x) \not \in M$ for all $m \leq T$. If $T>N\alpha(x)$, then
there must be a site $i \in \Lambda$ that is greater than $E_c$
for more than $\alpha(x)$ different times. This is impossible so
$m(x) \leq N\alpha(x)$. \qed


\subsection{\hpss Random excitations}\label{randex}

Define $\Sigma_N^+ = \Lambda^{\N}$ to be the set of
right-infinite $\Lambda$-sequences and let $\sigma_N^+:\Sigma_N^+
\rightarrow \Sigma_N^+$ be the left shift. We define a map
$\hat{f}:\Sigma_N^+ \times {\R}_{\geq 0}^N \rightarrow \Sigma_N^+
\times {\R}_{\geq 0}^N$  by
$$
\hat{f}({\bf t},x)=\begin{cases}  (\sigma_N^+ {\bf t},x+e_{t_0}) &
\,\mbox{if} \,\, x \in M \\ ({\bf t},f(x)) & \,\mbox{if} \,\, x
\not \in M
\end{cases}\,,
$$
where $e_{1},\dots,e_N$ is the standard basis in $ \mathbb{R}^N$.
We denote points in $\Sigma_N^+ \times \mathbb{R}_{\geq 0}^N$ by
$\hat{x}=({\bf t},x)$, and we define $\pi_u$ and $\pi_s$ to be the
projections to $\Sigma_N^+$ and $\R_{\geq 0}^N$ respectively.

\begin{prop}\po \label{4}
For all $\hat{x} \in \Sigma_N^+ \times \R_{\geq 0}^N$ it holds:
$$
\min \{ m \geq 0 \,|\,\forall i \in \Lambda \, \exists  m' \leq
m:\, (\pi_s \circ \hat{f}^{m'}(\hat{x}))_i>E_c\} \leq n(E_c,
\epsilon, \Lambda)\,,
$$
where
$$
n(E_c,\e, \Lambda)= N(NE_c+2) \Big{(}\Big{[}\frac{2d}{1-\e}
\Big{]}+1 \Big{)}^{\text{{\em diam}}(\Lambda)} \,,
$$
\end{prop}
\begin{proof}
In $N[E_c]+1$ time-steps, there must be an overcritical site.
Since in the relaxation process there is always an overcritical
site, then during arbitrary subsequent $N[E_c]+2$ time-steps an
exited site can be found. Hence after $N\xi(N[E_c]+2)$ time-steps
either all sites have been overcritical or there is a site that
has been overcritical at least $\xi$ times. However it follows
from the proof of Proposition \ref{2} that if one site is
overcritical
$$\xi=\Big{(}\Big{[}\frac{2d}{1-\e}
\Big{]}+1 \Big{)}^{\text{diam}(\Lambda)}$$
times, then all sites have been overcritical at least once.
\end{proof}
\\

For $x \in \mathbb{R}_{\geq 0}^N$ and $i \in \Lambda$ we define
$\tau(i,x):= \min \{n \in \N \,|\,f^n(x+e_i) \in M\}$. Proposition
\ref{2} assures us that this number is finite and
 $$
\max_{x\in M}\max_{i\in \Lambda}\tau(i,x)\le
\tau_m(E_c,\e,\Lambda)=N^2\Bigl(1+\frac1{NE_c}\Bigr) \Big{(}
\frac{2d}{1-\e}+1 \Big{)}^{\text{diam}(\Lambda)/2+1}.
$$
Thus we observe that neither $n(E_c,\e, \Lambda)$ nor
$\tau_m(E_c,\e,\Lambda)$ are uniformly bounded in $E_c$, but there
is the following alternative:

{\it There exists a constant $C_0$, not depending on the energy
$E_c$, such that either $n(E_c,\e, \Lambda)\le C_0$ or
$\tau_m(E_c,\e,\Lambda)\le C_0$.}

In fact, we can set $C_0=3N^2\bigl(\frac{2d}{1-\e}+1
\bigr)^{\text{diam}(\Lambda)+1}$. Thus we get that either
relaxation happen sufficiently fast or all the sites keep being
excited sufficiently often (uniformly in $E_c$).

But there does not exist such a bound uniform in $\e$ or $N$.


\subsection{\hpss Return maps} \label{return maps}

Let $\hat{x}=({\bf t},x) \in \Sigma_N^+ \times M$. For
$n=1,\dots,\tau(t_0,x)$ define $C_n(\hat{x})=\{ i \in \Lambda\,|\,
(\pi^s \circ \hat{f}^{n} \hat{x})_i>E_c \}$, and
$A(\hat{x})=\big{(}C_1(\hat{x}),\dots,C_{\tau(t_0,x)}
(\hat{x})\big{)}$.  We call $A(\hat{x})$ the {\em avalanche} of
the point $\hat{x}$. Let $\hat{M}:=\Sigma_N^+ \times M$ and define
an equivalence relation $\sim$ on $\hat{M}$ by
 $$
\hat{x} \sim \hat{y} \Leftrightarrow A(\hat{x})=A(\hat{y})\,.
 $$
This gives a partition of $\hat{M}$. From the definition it is clear
that $A(\hat{x})$ depends on $t_0$ and $x$ only. Hence partition
elements are of the form $[i] \times M_{ij}$, where
 $$
\forall i \in \Lambda:\, \bigcup_j M_{ij}=M
 $$
and $[i]=\{{\bf t} \in \Sigma_N^+\,|\,t_0=i\}$ is the cylinder of
the symbol $i$. We see that for each $i \in \Lambda$, the domains
$M_{i1},M_{i2},\dots$ are separated by segments of at most
$N!^{\tau_m}=\exp{(\tau_m(E_c,\e,\Lambda)\log N!)}$ hyperplanes.
Hence we have a finite number of domains $M_{i1},\dots,M_{iq_i}$
for each $i \in \Lambda$. By definition there is a unique
avalanche for each partition element $[i] \times M_{ij}$. We
denote this avalanche by $A_{ij}$. Its duration is
$\tau_{ij}:=\tau(i,x)$, for $x \in M_{ij}$, and define its size to
be $s_{ij}=\sum_{n=1}^{\tau_{ij}} |C_n|$.

We define the piecewise continuous map $F:\hat{M} \rightarrow
\hat{M}$  by
$$({\bf t},x) \mapsto (\sigma_N^+ {\bf t}, F_{t_0}x)$$ where
$F_i(x):=f^{\tau(i,x)}(x+e_i)$. We define $F_{ij}:=F_i
|_{M_{ij}}$.

 \begin{rk}\po\label{ps.vs.mt}
From a mathematical point of view the formulation $(\hat{M},F)$ is
a simplification compared to $(\Sigma_N^+ \times \R_{\geq
0}^N,\hat{f})$. However, the duration of avalanches are suppressed
so that all avalanches have the same duration. This is not
satisfactory from a physical point of view, and hence we call
$(\hat{M},F)$ the mathematical model and $(\Sigma_N^+ \times
\R_{\geq 0}^N,\hat{f})$  the physical model. We will later make a
rescaling of time in the physical model, so that the driving does
not become infinitely slow in the thermodynamic limit.
 \end{rk}

For each $x \in \mathbb{R}^N_{\geq 0}$ we define a matrix $Q(x)$
by
 $$
Q_{kl}(x)=\begin{cases}
\frac{1}{2d} \theta(x_l-E_c)   & \,\mbox{if} \,\, d_\Lambda(k,l)=1\,, \\
 0                             & \,\mbox{otherwise}\,,
\end{cases}
 $$
and a diagonal matrix $J(x)$ by
$J_{kl}(x)=(1-(1-\epsilon)\theta(x_l-E_c))\delta_{kl}$. Set
 \begin{equation}\label{defS}
S(x)=J(x)+(1-\epsilon)Q(x)
 \end{equation}
and observe that $f(x)=S(x)x$.  Let $x(1)=x+e_i$ and
$x(n)=f(x(n-1))$ for $n \in \{2,\dots,\tau(t,x)\}$. Then
$F_i(x)=L_i(x+e_i)$, where
$$
L_i(x)=S(x(\tau(x,i))) \dots S(x(1))\,.
$$
If $x,y \in M_{ij}$, then $\tau(i,x)=\tau(i,y)$ and the same
components of $x(n)$ and $y(n)$ are grater than $E_c$ for each
$n=1,\dots,\tau(t,x)$, so $L_i(x)=L_i(y)$. We define the linear
map $L_{ij}:=L_i(x)$ for $x \in M_{ij}$. We get
$F_i|_{M_{ij}}(x)=L_{ij}(x+e_i)$.

 \begin{dfn} \po \label{d1}
A sequence $\{(i_n,j_n)\,|\,1\leq n \leq \theta\}$ is said to be
admissible if
$$
\bigcap_{n=1}^\theta   (F_{i_{n-1} j_{n-1}}\circ\dots\circ F_{i_1
j_1})^{-1}(M_{i_{n}j_{n}})  \neq \emptyset\,,
$$
 \end{dfn}

 \begin{theorem}\po \label{5}
For all $i,j$ $\|L_{ij}\|_1 \leq 1$. Moreover for every constant
$c\in (0,1)$ there is a number $T \in \mathbb{N}$ such that for
every $\theta> T$ and admissible sequence $\{(i_n,j_n)\,|\,1\leq n
\leq \theta\}$ it holds:
$$
\|L_{i_\theta j_\theta} \dots L_{i_{1}j_{1}}\|_1<c\,.
$$
 \end{theorem}

\noindent
\begin{proof}
If $A$ is an $N \times N$ matrix we let $C_k(A)$ be its $k$-th
column. Observe that for any matrices $A$ and $B$ we have the
following formula:
 \begin{equation}\label{column}
\|C_k(AB)\|_1=\sum_l \|C_l(A)\|_1 B_{lk}\,.
 \end{equation}
By the construction: $\|C_k(S(x))\|_1 \leq 1$ for all $k \in
\Lambda$. Hence $\|L_{ij}\|_1 \leq 1$.

To prove the second statement we note that for $\e>0$ the diagonal
elements of the matrices $S(x)$ are non-zero and $\ge\e$.
Therefore
 $$
(S(x(m))S(x(m-1)))_{kl}\geq
\e\cdot\max\{S_{kl}(x(m)),S_{kl}(x(m-1))\}.
 $$
Moreover, $S_{kl}(x(m))>0$ if $x(m)_l>0$ and $d_\Lambda(k,l)=1$.
It follows that any admissible product $L_{i_\theta j_\theta}
\dots L_{i_1 j_1}$ of length $\theta \geq n(E_c, \epsilon,
\Lambda)$ is positive. By Proposition \ref{1} there must be at
least one column such that the sum over this column is less than
$1$, for some factor $L_{i_tj_t}$ and hence for the whole product.
Therefore the sum over each column of any admissible product of
length $2n(E_c,\epsilon, \Lambda)$ must be less than $1$. Let
$c_0<1$ be the maximal norm of all admissible products of length
$2 n(E_c,\epsilon,\Lambda)$.  For $k>k_0:=[\log c/\log c_0]+1$ we
have $c_0^k<c$ and hence $T=2k_0n(E_c,\epsilon, \Lambda)$ is the
required number.

The above argument does not apply to the case $\e=0$, and a different
proof must be given for this case (which actually works in general as well).
Take $\hat{x}\in \hat{M}$
and let $x(t) \in M$ be the projection of its orbit to $\hat{M}$.
Denote $S(x(t))$ by $S_t(\hat{x})$, and let
$$\tilde{S}_t(\hat{x})=S_t(\hat{x}) \cdots S_0(\hat{x})\,.$$
We make the following claims:
\begin{enumerate}
\item There exists $\bar{n} \in \N$ such that for all $l,m \in \Lambda$
and all $\hat{x} \in \hat{M}$ there is $t\le\bar{n}$ such that
$(\tilde{S}_t(\hat{x}))_{lm} \neq 0$.
\item For all $i \geq 0$ there exists $n_i \in \N$ such that
$\|C_m(\tilde{S}_t(\hat x))\|_1<1$ for all $t \geq n_i$, $\hat
x\in\hat M$ and all sites $m \in \Lambda$ with
$d_\Lambda(m,\partial \Lambda)\leq i$.
\end{enumerate}

The second claim for $i=\frac12\op{diam}(\Lambda)$
implies the statement of the theorem.

To see the first claim we fix $\hat x$ and let $U\subset\Lambda^2$
be the subset of the pairs $(l,m)$ with $(\tilde S_t)_{lm}=0$ for
all sufficiently large $t$. By inductively applying (\ref{column})
we see that the columns for $\tilde{S}_t(\hat{x})$ are non-zero
for all $t\geq 0$. So for all $\beta \in \Lambda$ there is $\alpha
\in \Lambda$ such that $(\alpha, \beta) \in \Lambda^2 \setminus
U$. Given sites $\alpha$ and $\beta$ we choose $t$ such that
$\tilde{S}_t(\hat{x})_{\alpha \beta} \neq 0$. Consider now column
$\alpha$ of the matrix
$\tilde{S}_{t+1}(\hat{x})=S_{t+1}(\hat{x})\tilde{S}_{t}(\hat{x})$.
If $\alpha$ is stable, i.e. $x(t+1)_\alpha\leq E_c$, then
$(S_{t+1}(\hat{x}))_{\alpha \alpha}>0$ and
$(\tilde{S}_{t+1}(\hat{x}))_{\alpha \beta} \neq 0$, so we just
repeat the argument. But the site $\alpha$ can not be stable for
more than $n(E_c,0,\Lambda)$ iterations. Hence we can with no loss
of generality choose $t$ such that $x(t+1)_\alpha > E_c$. Then the
column $\alpha$ of $S_{t+1}(\hat{x})$ has non-zero elements in all
position that correspond to neighbors of $\alpha$. Hence we obtain
that $(\alpha',\beta) \in \Lambda^2 \setminus U$ for all $\alpha'$
with $d_\Lambda(\alpha',\alpha)=1$. Any two points can be
connected by a path of neighbors, so $U=\emptyset$, and the first
claim follows. In fact, one can see that the bound $\bar n$ does
not depend on a choice of $\hat x$ and satisfies: $\bar n\le
\op{diam}(\Lambda)\cdot n(E_c,0,\Lambda)$.

To prove the second claim let us note that if
$\|C_k(\tilde{S}_{t}(\hat{x}))\|_1<1$, then
$\|C_k(\tilde{S}_{t+1}(\hat{x}))\|_1<1$ because by (\ref{column}):
$\|C_k(AB)\|\le\op{max}_l\|C_l(A)\|\cdot\|C_k(B)\|$.

We will use induction on $i$ starting from $i=0$. Take $k \in
\partial \Lambda$ and $t \leq \bar{n}$ such that
$(\tilde{S}_t(\hat{x}))_{kk}\neq 0$. If $x(t+1)_k > E_c$, then
$\|C_k(S_{t+1}(\hat{x}))\|_1<1$ and
$$
\|C_k(\tilde{S}_{t+1}(\hat{x}))\|_1=\sum_l
\|C_l(S_{t+1}(\hat{x}))\|_1 (\tilde{S}_{t}(\hat{x}))_{lk}<1\,,
$$
and so we have the desired inequality. If $x(t+1)_k \leq E_c$,
then $(S_{t+1}(\hat{x}))_{kk}=1$ and hence
$(\tilde{S}_{t+1}(\hat{x}))_{kk}\neq 0$. Then we repeat the
argument. Since no site can be stable for more than
$n(E_c,0,\Lambda)$ successive time-steps we obtain the claim for
$i=0$ with $n_0=n(E_c,0,\Lambda)+\bar{n}$.

Consider now the case $i>0$. For a site $m \in \Lambda$ with
$d(m,\partial \Lambda)=i$, we take $l\in \Lambda$ with
$d_\Lambda(l,m)=1$ and $d_\Lambda(l,\partial \Lambda)=i-1$. By the
first claim we find some $t\le\bar n$ such that
$(\tilde{S}_t(\hat{x}))_{lm} \neq \emptyset$, and by the induction
hypothesis for $t'\ge n_{i-1}$ we have:
 $$
\|C_l\big{(}S_{t+t'}(\hat{x})\cdots
S_{t+1}(\hat{x})\big{)}\|_1<1\,.
 $$
Using (\ref{column}) we obtain
$\|C_m\big{(}\tilde{S}_{t+t'}(\hat{x})\big{)}\|_1<1$. We can
choose $n_i=n_{i-1}+\bar{n}$.
\end{proof}

 \begin{lem}\po\label{6}
Let $E_c \geq \epsilon/(1-\epsilon)$. Then for any $\hat{x} \in
\hat{M}$, $n \in \N$ and $i,j \in C_n(\hat{x})$ we have
$d_\Lambda(i,j) \neq 1$.
 \end{lem}
 \begin{proof}
Take $\hat{x} \in M$ and let $E_n$ be the maximal energy of a site
in $C_n(\hat{x})$. Clearly $E_1 \leq E_c+1$ and
 $$
E_{n+1} \leq \max \big{\{} \max\{\epsilon E_n, E_c
\}+(1-\epsilon)E_n, E_c+1 \big{\}}\,.
 $$
From this we see by induction that
 $$
E_n \leq \max \big{\{} \frac{E_c}{\epsilon}, E_c+1
\big{\}}=\frac{E_c}{\epsilon}\,,
 $$
so $\epsilon E_n \leq E_c$ for all $n \in \N$, and this means that
a site cannot be overcritical in two successive time-steps (for
$\e=0$ the above argument does not work, but the statement holds
obviously).

All avalanches start with a single site. Let $C_1(\hat{x})=\{i\}$.
Then $d(i,j)=1$ for all $j \in C_2(\hat{x})$. This implies that
any two elements of $C_2(\hat{x})$ can be connected with a path of
length 2, so no two sites of $C_2(\hat{x})$ are nearest neighbors.
If there exists a path of even length between two points in
$\Lambda$, then all paths connecting these points are of even
length. Therefore we can repeat the argument proving by induction
that $d_\Lambda(i,j) \in 2 {\mathbb Z}$ for all $i,j \in
C_n(\hat{x})$.
\end{proof}

 \begin{prop}\po \label{7}
The linear maps $L_{ij}$ are all invertible whenever
$\epsilon\ge 1/2$ or $\e>0$ and $E_c \geq \epsilon/(1-\epsilon)$.
If we have $E_c \geq \epsilon/(1-\epsilon)$, then
$$
\det L_{ij}=\epsilon^{s_{ij}}\,.
$$
 \end{prop}
\begin{proof}
Take arbitrary $x \in \R^N_{\geq 0}$. First we observe that since
the sum over each column of $Q(x)$ is less than or equal to $1$,
we have $\|Q(x)v\|_1 \leq \|v\|_1$ for each $v \in \R^N$. This
implies that
\begin{eqnarray*}
\|S(x)v\|_1 &=&\|(J(x)+(1-\epsilon)Q(x))v\|_1 \\ &\geq&
\|J(x)v\|_1-(1-\epsilon)\|Q(x)v\|_1 \\
&\geq&(2\epsilon-1)\|v\|_1\,.
\end{eqnarray*} If $\epsilon>1/2$, then $S(x)v \neq 0$ for all $v \neq 0$, so we have
invertibility.

For $\e=1/2$ the claim follows since in the above chain of
inequalities at least one is strict if $v\ne0$. In fact, if
$J(x)v=\e v$, then $v_i=0$ for all relaxed sites $i$. We claim
that the equality $Q(x)v=v$ is impossible. To see this denote by
$\tilde Q$ the minor-matrix formed by the rows and columns of
$Q(x)$, corresponding to exited sites, and denote by $\tilde v$ be
the respective reduced vector. Then $\tilde Q\tilde v=\tilde v$.

Let $U$ be the set of overcritical sites $k$ with $v_k=\max v_l$
(we suppose it is positive, multiplying by $-1$ in the opposite
case). Choose a boundary site $k\in U$, i.e. the number of
neighbors $l$ to $k$ with $v_l=v_k$ is less than $2d$. Then:
 $$
v_k=\sum_l\tilde Q_{kl}v_l<v_k\sum_l\tilde Q_{kl}\le v_k.
 $$
This contradiction yields the result.

Finally consider the last statement about the case $E_c \geq
\epsilon/(1-\epsilon)$. It is proved by reducing the matrix
$S(x)$. If $x_i \leq E_c$, then column $C_i(S(x))$ equals
$(0,\dots,0,1,0,\dots0)^T$, where the 1 is in the $i^\text{th}$
position. We can start the decomposition of $\det S(x)$ with
column $i$, and hence we see that row $i$ and column $i$ can be
removed from $S(x)$ without changing the determinant. We remove
all rows and columns that correspond to relaxed sites. If
$\rho(x)$ is the number of overcritical sites of $x$, we get a
$\rho(x) \times \rho(x)$ matrix $S_{red}(x)$. If site $k$ is
overcritical then $J_{kk}(x)=\epsilon$. If $E_c \geq
\epsilon/(1-\epsilon)$, then it follows from Lemma \ref{6} that
all nearest neighbors of $k$ are relaxed. Hence column $k$ of
$Q(x)$ has only zero entries. This shows that
$S_{red}(x)=\mbox{diag}(\epsilon,\dots,\epsilon)$. Then
\begin{equation*}
\det S(x)=\det S_{red}(x)=\epsilon^{\rho(x)}\,,
\end{equation*} so $\det L_{ij}=\epsilon^{s_{ij}}$.
\end{proof}

\begin{rk}\po
In the original model of Zhang one has $\epsilon=0$ in which case
$\det L_{ij}=0$ if $L_{ij}\ne\1$. But it is not true that
non-trivial kernels can occur for $\epsilon=0$ only, contrary to
what was stated in \cite{BCK}. A simple counter-example is the
case $N=2$, $E_c=1/3$ and $\epsilon=1/3$. For $x_1>0$ and
$2x_1+3x_2 < 1$ we have:
$$
F_1 \Big{(} \begin{bmatrix} x_1 \\ x_2 \end{bmatrix} \Big{)}
=\frac{1}{9}
\begin{bmatrix} 2 & 3 \\ 2 & 3  \end{bmatrix} \begin{bmatrix} x_1+1 \\ x_2
\end{bmatrix}\,.
$$
and so $\det L_{12}=0$.
\end{rk}

Having non-degenerate maps in the model is more convenient from
the point of view of mathematical tools (though from a physical
viewpoint it can make no big difference between degenerate and
close-to-degenerate systems). Fortunately, degenerations occur
only for a negligible set of parameters.

 \begin{theorem}\po\label{th_invert}
The maps $L_{ij}$ are invertible for almost all $(\e,E_c)$. In
fact, they are invertible for the parameters complimentary to the
set $\Xi\subset[0,1)\times(0,\infty)$, which consists of a finite
set of vertical intervals for fixed $d$ and $N$.
 \end{theorem}

\begin{proof}
Fix an avalanche $A_{ij}$ and let $L_{ij}^\e$ be the corresponding
linear maps (we stress dependence on $\e$). These maps are the
compositions of elementary matrices $S^\e(x(\tau(x,i))) \dots
S^\e(x(1))$, with the factors from (\ref{defS})
 $$
S^\e(x(t))=\1+(\e-1)(\tfrac{dJ}{d\e}-Q)(x(t))
 $$
being polynomial in $\e$ and independent of the choice of $x=x(1)
\in M_{ij}$. The condition $\det L_{ij}^\e=0$ is equivalent to
$\det S^\e(x(t))=0$ for some $t$. Denoting by $\op{Sp}_-(T)$ the
negative part of the spectrum of $T$, we get:
$\e\in1+\op{Sp}_-(Q-\frac{dJ}{d\e})^{-1}$.

There are only finite number of possibilities for the matrix
$S^\e(x)$ (though a countable number for their compositions
$L_{ij}^\e$, the length of which grow as $E_c\to0$). Thus we
obtain $k=k(d,N)$ different values of $\e$ for which $\det
S^\e(x)=0$: $\{\e_a\}_{a=1}^k$. For each $\e_a$ there is the
maximal value $E_c^a$ of $E_c$ (finite if $\e_a\ne0$), where the
corresponding matrix $S^\e(x)$ can appear in the avalanche. Thus
the set of degenerate systems is
$\{(\e,E_c)\,|\,\e=\e_a,0<E_c<E_c^a\}$.
 \end{proof}
 \vskip4pt

By proposition \ref{7} $\Xi$ does not intersect the set
$\{\e\ge1/2\}\cup\{E_c\ge\e/(1-\e)\}$.


\section{\hps Removability of singularities and
coding}\label{sec_2}

The map $F$ may be considered as a piecewise affine map
$F:I\times M \rightarrow I \times M$, where $I=[0,1]$ and $F(t,
x)=(Nt \mod 1, F_{[Nt]}(x))$. The map $t \mapsto  Nt \mod 1$ is
not conjugated to $\sigma_N^+$ since the points $m/N^k \in I$ do
not have unique representations in $\Sigma_N^+$. However the sets
$\{m/N^k\}\times M \subset I \times M$ are singularities, and
following the standard approach for piecewise affine maps, should
be removed.

In some physical systems, like the Belykh family, the
singularities propagate, intersecting themselves transversally.
The Zhang model is not a general position system in this respect,
because singularities $\{m/N^k\}\times M \subset I \times
M=\hat{M}$ map into themselves, forming zero angle.


\subsection{\hpss Construction of attractors} \label{constofattr}

Define the (spatial) singularity set $S(F)=\cup_{ij} \partial
M_{ij}$. Then $U=M \setminus S(F)$  consists of a collection of
open connected sets $\mathcal Z=\{Z\}$.  Let $U_0:=U$ and
$$
U_n:=\bigcup_{i \in \Lambda} F_i ( U_{n-1})\cap U\,.
$$
We say that $x \in S(F)$ is a non-essential singularity of order
$m$ if there exists $\epsilon>0$ and $m>0$ such that
$$
\text{card}\{Z \in \mathcal Z\,|\, U_n\cap B_\epsilon(x)\cap Z
\neq \emptyset\} \leq 1
$$
for all $n>m$. Denote the set of non-essential singularities of
order $m$ by $NES(F;m)$ and let $NES(F):=\cup_{m \geq 0} NES(F;m)$
be the set of all non-essential singularities. Define $ES(F)=S(F)
\setminus NES(F)$ to be the collection of essential singularities.
Observe that there is a natural extension of $F$ to $V_0=U \cup
NES(F)$. In the following we let $F$ denote the extended map. As
above we define
$$
V_n:=\bigcup_{i \in \Lambda} F_i ( V_{n-1}) \cap V_0\,.
$$
Let ${\mathcal X}=\cap_{n \geq 0} V_n$ and $\D=\Sigma_N^+ \times
{\mathcal X}$. Clearly $F(\D)=\D$. The set $\Y=\overline{\mathcal
X}$ is called the {\em physical} (or spatial) attractor of $F$,
and $\mathcal A=\Sigma_N^+ \times \Y=\overline{\D}$ is the
extended attractor of $F$.

\begin{prop}\po \label{8}
$F|_\D$ is continuous.
\end{prop}

\begin{proof}
The set $\D$ intersects non-essential singularities only. Hence we
must show that if $x$ is a non-essential singularity in $\D$, then
the extension of each $F_i$ to $NES(F)$ is continuous at the point
$x$. Choose $m \in \N$ and $\varepsilon>0$ such that
$B_\varepsilon(x) \cap U_n$ intersects only one partition element
$Z\in\mathcal Z$ for $n> m$ and let $y\in
B_{\varepsilon/2}(x)\cap\D$. Then $B_{\varepsilon/2}(y)\subset
B_\varepsilon(x)$ and so $B_{\varepsilon/2}(y)\cap U_n$ intersects
the same partition element $Z$. So $x$ and $y$ are mapped by the
same affine map $F_i|_Z$ for each $i \in \Lambda$. The claim
follows.
\end{proof}
\vspace{5pt}

In general, the map $F$ does not have a continuous extension to
$\mathcal A$, but only to $\mathcal A\setminus(\Sigma_N^+\times
ES(F))$. Actually, if $x\in ES(F)\cap\Y$ lies on the boundary of
several continuity partitions for $F_i$, then there are several
extensions of $F$ to $(i,x)$. Thus we can continuously extend $F$
to $\mathcal A$ only when the essential singularities do not
intersect the attractor (are removable).


\subsection{\hpss Symbolic Coding}

If the singularities can affect the dynamics only for a finite
number of iterations, then the dynamics can be well
approximated by a topological Markov chain.

\begin{dfn}\po\label{d3}
We say that singularities are removable if there exists $m\in \N$
such that $S(F)=NES(F,m)$.
\end{dfn}

The physically most relevant observables $\phi:\hat{M} \rightarrow \R$ are those
that are determined by avalanches. We say that
$\phi$ is an {\em avalanche observable} if it is constant on
continuity domains $[i]\times M_{ij}$.

 \begin{theorem}\po \label{10}
If singularities are removable, then the map $F$ is well-defined
and continuous on $\mathcal A$ and there is a topological
Markov-chain $(\Sigma_A^+,\sigma_A^+)$ and a continuous
semi-conjugancy $g:\mathcal A \rightarrow \Sigma_A^+$ such that
for all $\hat{x}, \hat{y} \in {\mathcal A}$ and for all avalanche
observables $\phi$ we have:
 $$
g(\hat{x})=g(\hat{y}) \Rightarrow \phi(F^n(\hat x))=\phi(F^n(\hat
y))\,\, \forall n\geq 0\,.
 $$
The Markov-chain is determined by a
matrix $A$ which has a maximal eigenvalue equal to $N$.
 \end{theorem}

 \begin{rk} \po
It is clear that all properties related to distribution of
avalanche size, duration, area, etc. are invariant under a
semi-conjugancy such as this. Observe that for each avalanche
observable $\phi$ on $\mathcal A$, there is a unique observable
$\phi':\Sigma_A^+ \rightarrow \R$ such that $\phi=\phi' \circ g$.
Suppose we have a measure $\mu$ on $\mathcal A$, and let
$\nu=g_*\mu$. If $\phi$ is an avalanche observable on $\mathcal
A$, then the statistical properties of $\phi$ with respect to
$\mu$ are equivalent to the statistical properties of $\phi'$ with
respect to $\nu$. In this $(\sigma_A^+, \Sigma_A^+)$ is a good
approximation to $F|_{\mathcal A}$. The coding gives estimates on
entropy and growth of periodic points, but these estimates are
asymptotically no better than what we get from the trivial
semi-conjugancy $\hat{M} \rightarrow \Sigma_N^+$.
 \end{rk}

\begin{proof}
Singularities are removable so there exists an integer $m \in \N$
such that $\pi_s \circ F^m(\hat{M})$ only intersects trivial
singularities. Let $X_1,\dots, X_s$ be the closure of the
connected components of  $\pi_s \circ F^m(\hat{M})$. $F$ is well
defined and continuous on these components. Let $Y_1,\dots,Y_s$ be
the intersections of the components  $X_1,\dots,X_s$ with $\Y$. We
construct the partition ${\mathcal R}=\{[i] \times Y_k\}$ and
enumerate it so that ${\mathcal R}=\{R_1,\dots,R_r\}$, where
$r=Ns$.

Let $A=\|a_{ij}\|$ be the $r\times r$ matrix defined by the rule:
$a_{ij}=1$ if $F(R_i) \cap R_j \neq \emptyset$, and $a_{ij}=0$
otherwise. A sequence $R_{\omega_0}R_{\omega_1}\dots$ is legal if
$a_{\omega_{t-1} \omega_{t}}=1$ for all $t \in \N$. Define
$g:{\mathcal A} \rightarrow \Sigma_A^+$ by $g(\hat{x})= (\omega_0
\omega_1\dots \omega_t \dots)$, where $F^t(\hat{x}) \in
R_{\omega_t}$. To prove that $g$ is surjective it suffices to show
that for each legal sequence $R_{\omega_0}R_{\omega_1}\dots$,
there is a point $\hat{x} \in \pi_s \circ F^m(\hat{M})$ such that
$F^i(\hat{x}) \in R_{\omega_t}$ for all $i\in \N$. Note that each
$\omega$ can be written as a pair $(t,k)$, where $t\in
\{1,\dots,N\}$ and $k \in \{1,\dots,s\}$. Hence we can write
 $$
\bigcap_{n=0}^\infty F^{-n}( R_{\omega_n})=\bigcap_{n=0}^\infty
F^{-n}([t_n]\times Y_{k_n} )=\{{\bf t}\} \times
\bigcap_{n=0}^\infty F_{t_0}^{-1}\circ \dots \circ
F_{t_{n-1}}^{-1}(Y_{k_n})\,.
 $$
The continuous image of a connected set is connected, so for each
$i=1,\dots,N$ and each $k=1,\dots,s$ there is a unique $l \in
\{1,\dots,s\}$ such that $F_i(X_k) \subset X_l$.
This implies that we have a nested sequence
$$
Y_0 \subset F_{t_0}^{-1}(Y_{k_1}) \subset F_{t_0}^{-1}\circ
F_{t_1}^{-1}(Y_{k_1}) \subset \dots
$$
and hence the intersection is non-empty.

It is clear that $g_{\mathcal R}$ is continuous (see \cite{R} for
details). Since the partition ${\mathcal R}$ is a refinement of
the continuity partition the conjugancy will be injective up to
the classes of points that follow the same continuity domains.
Hence if $\phi(F^n \hat{x}) \neq \phi(F^n \hat{y})$ for some
avalanche observable $\phi$ and some $n \geq 0$, then $g(\hat{x})
\neq g(\hat{y})$.
\end{proof}

\noindent
\begin{rk}\po\label{rkk4}
Suppose we modify the Zhang model by using a full shift
$(\Sigma_N, \sigma_N)$ as the excitation factor. It is then
possible that the modified map $F$ is injective on $\Sigma_N
\times \Y$. Since we have strict attraction in the spatial factor
after a fixed number of iterations it is clear that we can then
obtain an injective coding, and hence a topological conjugancy.
However, if we make this modification it is not clear that
$\Sigma_N \times \Y$ equals the set
$$
\Omega=\overline{\bigcap_{n=-\infty}^\infty F^{n}(\Sigma_N \times
M)}\,.
$$
In fact if the maps $F_i|_\Y$ are all injective, then $F|_\Omega$
is invertible, but $F|_{\Sigma_N \times \Y}$ is typically
non-invertible. The reason for this is that, due to contraction, a
point $x \in \Y$ does not have preimages for all the maps $F_i$
and $F_i$ are invertible only on $F_i(\Y)\subset\Y$. So to obtain
invertibility we must turn to the attractor $\Omega$. From a
physical point of view the spatial attractor is of the great
interest, so it is desirable to have an attractor which is a
Cartesian product of the Bernoulli shift and the spatial attractor
$\Y$.
\end{rk}

We can always construct a coding of $F|_{\D}$ (even in
non-removable case) by choosing a partition ${\mathcal
R}=\{R_1,\dots R_r\}$, and taking $g_{\mathcal R}:\D \rightarrow
\{1,\dots,r\}^\N$ to be the map sending a point $\hat{x} \in \D$
to the unique sequence $\omega \in \{1,\dots r\}^\N$ such that
$F^t(\hat{x}) \in R_{\omega_t}$ for all $t \geq 0$. But there is
no reason, however, to expect $g_{\mathcal R}(\D)$ to be a
topological Markov chain, cf. \cite{BCK}.


\section{\hps Metric properties}\label{sec_3}

The natural volume on $\hat{M}$ is given by the product measure of
the uniform Bernoulli measure on $\Sigma_N^+$ and the Lebesgue
measure on $M$. By iterating this measure (and averaging) we can
construct SRB-measures. However it can happen that the measures
constructed are supported on essential singularities, where it is
not possible to define the dynamics in such a way that the measure
is invariant. Hence we must give some conditions to ensure the
existence of SRB-measures. If there is an SRB-measure it is
characterized by the fact that its projection to $\Sigma_N^+$
coincides with the uniform Bernoulli measure. From this it follows
that any SRB-measure is a measure of maximal entropy. In
situations where the system allows symbolic coding the SRB-measure
corresponds to the Perry measure on the topological Markov chain
$\Sigma_A^+$.


\subsection{\hpss Existence and characterization of SRB-measures}

Let $m=m^u \times m^s$, where $m^u=\mu_{\text{{\tiny Ber}}}$ is
the uniform Bernoulli measure on $\Sigma_N^+$, and
$m^s=\mu_{\text{{\tiny Leb}}}$ is the Lebesgue measure on $M$. We
say that an invariant Borel probability measure $\mu$ on $\hat{M}$
has the SRB-property if there exists a measurable invariant set $G
\subset \hat{M}$ such that
 \begin{enumerate}
\item $m(G)>0$
\item $m^u(\pi_u(G))=1$
\item All points $\hat{x} \in G$ are future generic with respect to $\mu$, i.e.
$$\frac{1}{n}
\sum_{t=0}^{n-1} \phi(F^t \hat x) \rightarrow \int \phi
\,d\mu\,,$$ for all $\hat{x} \in G$ and all continuous functions
$\phi:\hat{M}\to\R$.
 \end{enumerate}

For Axiom A attractors one can ensure the existence of measures
for which the set of generic points has full Lebesgue measure and
this is equivalent to saying that the canonical family of
conditional measures on the unstable manifolds are absolutely
continuous with respect to the Lebesgue measure. For
non-invertible maps one can in general only expect the set of
generic points to have positive measure and hence it is
unreasonable to require that $m(G)=1$. Condition 2 is (for
physical reasons) important in the Zhang model. It means that the
statistical properties do not depend on the choice of a generic
sequence ${\bf t}$ of excitations. (This is always implicitly
assumed in the numerical investigations of the Zhang model that
can be found in the physical literature.) Moreover, condition 2
will be satisfied for the SRB-measures that can be constructed by
iterating the measure $m$.

By a standard approach we can give conditions for existence of
SRB-measures that hold if singularities are removable, but it is
not known if these conditions hold in all non-removable
situations.

 \begin{prop}\po \label{11}
Let $(\e,E_c)$ does not belong to the negligible set $\Xi$ of
Theorem \ref{th_invert}. If there exists $n\geq 0$, $C>0$ and
$q>0$ such that
 $$
\forall \delta>0, \forall t \geq 0:\, m \Big{(} F^{-t}\big{(}
\Sigma_N^+ \times U_\delta( ES(F;n) \big{)} \Big{)} \leq C
\delta^q\,,
 $$
then there exists a set $\D \subset \hat{M}$ (constructed in \S
\ref{constofattr}), which may intersect singularities, and a
natural extension of $F$ to $\D$ such that $F(\D)=\D$. Moreover
the set $\D$ carries an $F$-invariant Borel probability measure
with the SRB-property.
 \end{prop}

 \begin{rk}\po
Proposition \ref{11} is a simple modification of the result of
Schmeling and Troubetzkoy \cite{ST}. In their paper the conditions
for existence are in general too restrictive for the Zhang model.
In fact, in Example A of \S \ref{examples} we show a situation
where the SRB-measure constructed in \cite{ST} does not exist, but
we clearly have existence of a physically relevant measure. The
reason for this paradox is that one in general remove all singular
points on the construction of the attractor, even if there is a
natural extension of $F$ to the points of singularity.
(Proposition \ref{11} obviously applies to this example since
$ES(F;5)=\emptyset$.)
 \end{rk}

\noindent
\begin{proof}
In \cite{ST} it is shown that a piecewise smooth map $f$ with
singularity set $S$ has a measure, not supported on singularities,
such that the set of generic points has positive Lebesgue measure.
They require that the following conditions are satisfied:
\begin{enumerate}
\item The restrictions of $f$ to each of its continuity domains
are diffeomorphisms onto their image.

\item The second differentials $D^2f_x$ does not grow too fast close
to singularities. (See \cite{ST} for a more precise formulation.)

\item $f$ is hyperbolic. In this context this means that there are
constants $C>0$ and $\lambda \in (0,1)$ such that for all  $x \not
\in S$ there is a splitting of the tangent space at $x$ into
subspaces $E^+(x)$ and $E^-(x)$. There are cones $C^+(x)$ and
$C^-(x)$ around $E^+(x)$ and $E^-(x)$ that are invariant under
$Df_x$ and $Df_x^{-1}$ respectively. The angles between the
$C^+(x)$ and $C^-(x)$ are bounded away from zero, and for all
points $x$ that do not intersect singularities in the first $n$
iterations it holds:
$$\|D_xf^n(v)\| \geq C^{-1}\lambda^{-n}\|v\|\,\,\text{for}\,\,\,
v \in C^+(x)\,,$$ and
$$\|D_xf^n(v)\| \leq C\lambda^{n}\|v\|\,\,\text{for}\,\,\,
v \in C^-(x)\,.$$
\item There exists $C>0$ and $q>0$ such that
$m(f^{-t}(U_\varepsilon(S))) \leq C \varepsilon^q$ for all
$\varepsilon>0$ and all $t\in \N$.
\end{enumerate}
We apply this result to the map $F|_{\Sigma_N^+ \times (U_n
\setminus ES(F;n))}$. The singularity set for this map is
contained in $\Sigma_N^+ \times ES(F)$, so by assumption condition
4 is satisfied. Condition 1 follows from Proposition \ref{7},
condition two is obviously satisfied since $F$ is piecewise affine
and condition 3 follows from Theorem \ref{5} with
$E^+=\R^1\oplus0$, $E^-=0\oplus\R^N$ and $C^\pm$ being the regular
cones around them (actually Theorem \ref{5} ensures hyperbolicity
for some iterate $F^T$, which implies the claim).

In \cite{ST} the measures are constructed by iterating $m$,
averaging and taking a weak limit. It is clear that, in the Zhang
model, any measure obtained in this way will satisfy condition 2
in our definition of an SRB-measure.
\end{proof}
\\

If an SRB-measure exists it can be characterized by a number of
different properties. From a physical perspective it is reasonable
to require that a relevant invariant measure should preserve the
uniform  Bernoulli structure on $\Sigma_N^+$. This corresponds to
the Lebesgue measure on $[0,1]$ in the alternative formulation of
the map $F$, and hence to absolutely continuous measure
conditional measures on the unstable space $[0,1]$.

 \begin{prop}\po \label{12}
If $\mu$ is an SRB-measure on $\D$, then $\mu^u:=(\pi_u)_* \mu$ is
the uniform Bernoulli measure on $\Sigma_N^+$.
 \end{prop}

\noindent
\begin{proof}
There is a set $A=\pi_u(G)$ of full $m^u$-measure, such that all
${\bf t} \in  A$ are generic with respect to $\mu^u$. Take a
continuous function $\phi:\Sigma_N^+ \rightarrow \mathbb{R}$. Then
 $$
\int \phi\,d\mu^u= \lim_{n \rightarrow \infty} \frac{1}{n}
\sum_{t=0}^{n-1} \phi((\sigma_\Lambda^+)^t {\bf t})=\int
\phi\,dm^u \,,
 $$
where the left equality holds for ${\bf t}\in A$ and the right one
for ${\bf t}\in B$ with $B\subset\Sigma_N^+$ a subset of full
$m^u$-measure (from Birkhoff ergodic theorem). Since $A\cap
B\ne\emptyset$, we get: $\int \phi \,d\mu^u=\int \phi\,dm^u$ for
all continuous functions $\phi$.
\end{proof}


\subsection{\hpss Measures of maximal entropy} \label{meaofmaxentr}

Suppose that there exists an invariant Borel probability measure
$\mu$ on $\hat{M}$. Let $\mu^u:=(\pi_u)_* \mu$ and let $\{
\nu_{\bf t}\}$ to be the canonical family of conditional measures
on the fibers $\pi_u^{-1}(\{\bf t \})$. By the Abramov-Rokhlin
formula
$$
h_\mu (F)=h_{\mu^u}(\sigma_N^+)+h_\mu(F|\sigma_N^+) \,,
$$
where
$$
h_\mu(F|\sigma_N^+;\mathcal Q) =  \lim_{n \rightarrow \infty}
\frac{1}{n} \int H_{\nu_{{\bf t}}} \Big{(} \bigvee_{k=0}^{n-1}
(F_{t_{k-1}} \circ\dots\circ F_{t_0})^{-1}({\mathcal Q}) \Big{)}
\, d\mu^u({\bf t})\,,
$$
for a partition $\mathcal Q$ and
$$
h_\mu(F|\sigma_N^+)=\sup_{\mathcal Q} h_\mu(F|\sigma_N^+;\mathcal
Q)\,,
$$
The supremum is taken over all finite measurable partitions
${\mathcal Q}$ of $\hat{M}$. This formula was originally proved
for product measures by Abramov and Rokhlin \cite{AR} and extended
to arbitrary skew products by Bogenschultz and Crauel \cite{BC}.

Below we use the notation $F_{\bf t}$ for the dynamics over a
pre-fixed sequence ${\bf t}=(t_0t_1\dots)\in\Sigma_N^+$. By
$n$-the iteration we mean the map $F^n_{\bf
t}=F_{t_{n-1}}\circ\dots\circ F_{t_0}$.

 \begin{theorem}\po \label{17}
If $\mu$ is an invariant Borel probability measure on $\hat{M}$,
then $h_\mu(F) = h_{\mu^u}(\sigma_N^+)$.
 \end{theorem}

\noindent
\begin{proof}
We prove the proposition by estimating $h_\mu(F|\sigma_N^+)$ from
above. Let ${\mathcal Q}$ be a partition of $M$ and
 $$
\mathcal Q_{t_0\dots t_{n-1}}:=\bigvee_{k=0}^{n-1} (F_{t_{k-1}}
\circ\dots\circ F_{t_0})^{-1}({\mathcal Q})=\bigvee_{k=0}^{n-1}
F_{\bf t}^{-k}({\mathcal Q})\,.
 $$
Fix $\varepsilon>0$ and choose the partition $\mathcal Q$ such
that $h_\mu(F|\sigma_N^+;\mathcal Q) +\varepsilon\geq
h_\mu(F|\sigma_N^+)$. The maps $F_i$ are non-expanding so it
follows from the Ruelle-Margulis inequality \cite{KH} that
 $$
\frac1n H_{\nu_{\bf t}} \Big{(}\bigvee_{k=0}^{n-1} F_{\bf
t}^{-k}({\mathcal Q})\Big{)}\longrightarrow h_{\nu_{\bf t}}(F_{\bf
t})=0.
 $$
Moreover the convergence is $\mu^u$-uniform and so the same holds
for the integrals. Another way to see it is via the multiplicity
notion of \S\ref{evalofentr} (then ${\mathcal Q}$ should be
subordinate to each continuity partition
$\{M_{ij}\,|\,j=1,\dots,q_i\}$):
 $$
h_\mu(F|\sigma_N^+;\mathcal Q)\le
\lim_{n\to\infty}\frac1n\log\max_{|{\bf t}|\leq n}
\op{mult}(\mathcal Q_{t_0\dots t_{n-1}}\cap\op{supp}(\nu_{\bf
t})).
 $$
Therefore $h_\mu(F|\sigma_N^+)\le \e$. Let $\e\to0$.
\end{proof} \\

 \begin{rk}\po
Theorem \ref{17} is a partial case of Theorem 4 from \cite{KR2}.
 \end{rk}

We say that an invariant measure $\mu$ is maximal if
$h_\mu(F)=\sup_\nu h_\nu(F)$, where the supremum is taken over all
invariant Borel probability measures on $M$. It follows from
Theorem \ref{17} that $h_\mu(F)\le\log N$. So if $\h(F)>\log N$,
the variational principle fails (this can happen for piece-wise
affine systems, see \cite{KR1,KR2}). But we show in
\S\ref{evalofentr} that the abnormal growth of $\h(F)$ does not
occur in the Zhang model, at least for generic values of
parameters $E_c,\e$.

 \begin{cor}\po\label{corr1}
Any SRB-measure on $\D$ has entropy $h_\mu(F)=\log N$ and is hence
a maximal measure.
 \end{cor}

 \begin{cor}\po
Suppose singularities are removable and that $\mu$ is an
SRB-measure on ${\mathcal A}$. Let $g:{\mathcal A} \rightarrow
\Sigma_A^+$ be the semi-conjugacy constructed in the proof of
Theorem \ref{10}. If $(\sigma_A^+,\Sigma_A^+)$ is topologically
transitive, then $g_*\mu$ is the Perry measure on $\Sigma_A^+$.
 \end{cor}
\begin{proof}
A transitive topological Markov chain has a unique measure of
maximal entropy. This measure is called the Perry measure
\cite{KH}.
\end{proof}


\subsection{\hpss Hyperbolic structure} \label{hypstr}

There are several ways to define Lyapunov exponents for the Zhang
model. The Zhang model can be represented as a piecewise affine
map, where Bernoulli shift is represented as the expanding map $t
\mapsto Nt \mod 1$ of the interval (see \S \ref{evalofentr}).
Hence it is clear that there is one positive Lyapunov exponent
$\chi^+_0=\log N$. We define the other exponents by introducing
the co-cycle $\T : \hat{M} \rightarrow GL(N,\R)$, defined by
$\T(\hat{x})=L_{ij}$, where $\hat{x} \in [i] \times M_{ij}$. For
$\hat {x} \in \hat{M}$ and $v \in \R^N \setminus \{0\}$ we define
$$
\chi(\hat x,v)= \mathop{\overline\lim}\limits_{n \rightarrow
\infty} \frac1n\log \frac{\| \T( F^{n-1}(\hat{x}) ) \dots \T(
F(\hat{x}) )  \T( \hat{x} ) v\|} { \|v\| }\,.
$$
It is a general fact that the function $\chi(\hat{x},\cdot)$ takes
at most $N$ different values $\chi^-_1(\hat{x})\ge\dots \geq
\chi^-_N(\hat{x})$.

 \begin{prop}\po \label{hyp}
For all $\hat{x} \in \hat{M}$ the Lyapunov spectrum is:
 $$
\chi^+_0=\log N>0>\chi^-_1(\hat{x})\ge \dots \geq \chi^-_N(\hat{x})
 $$
and for $(\e,E_c)$ outside the negligible set $\Xi$ from Theorem
\ref{th_invert}:  $\chi^-_N(\hat{x})>-\infty$.
 \end{prop}

 \begin{proof}
From Theorem \ref{5} we know that there exists $T \in \N$ and $c
\in (0,1)$ such that
$$
\|\T(F^{T-1}(\hat{x})) \dots \T(F(\hat{x}))\T(\hat{x})\|\leq c
$$
for all $\hat{x} \in \hat{M}$. It immediately follows
that $\chi(\hat{x},v) \leq T^{-1} \log c <0$.

For $\epsilon\ge1/2$ and arbitrary $E_c$ or for $\e>0$ and
$E_c\ge(1+\e)/(1-\e)$ all linear maps are invertible, and so for
all $\hat{x} \in \hat{M}$ and all $v \in \R^N \setminus \{0\}$ we
have $\chi(\hat{x},v)\geq \log k$, where $k=\min_{ij} \min
\text{Sp}(L_{ij})>0$.
 \end{proof} \\

If there exists a unique SRB-measure, then it follows from the
Osceledec theorem that there are numbers $\chi^-_1, \dots
,\chi^-_N$ such that $\chi_i(\hat{x})=\chi_i^-$ for Lebesgue
almost every $\hat{x} \in \hat{M}$. The numbers
$\chi_0^+,\chi_1^-,\dots, \chi^-_N$ are the Lyapunov exponents of
the Zhang model. From Proposition \ref{hyp} it follows that the
Zhang model is hyperbolic in the sense that the Lyapunov spectrum
consists of:
$$
\chi^+_0 =\log N > 0 > \chi^-_1 \ge \dots \geq
\chi^-_N>-\infty\,.
$$
We see from Corollary \ref{corr1} that the Pesin formula
$h_\mu(F)=\chi^+$ holds for any SRB-measure. If there is no
SRB-measure then the Lyapunov spectrum should be defined as
functions on $M$:
$$
\chi_i(x)=\int_{\Sigma_N^+}\chi_i({\bf t},x) \,d \mu_{\text{{\tiny Ber}}}\,,
$$
where $\mu_{\text{{\tiny Ber}}}$ is the uniform Bernoulli measure on $\Sigma_N^+$.


\subsection{\hpss Entropy of physical vs. mathematical
models}\label{305}

We can reformulate the Zhang system as the map $\hat f:\hat
B\to\hat B$, where $\hat B=\cup_{i\ge0}\hat f^i(\hat M)$ is a
compact $\hat f$-invariant subset of $\Sigma_N^+\times\R^N$. We
wish to compare this to the induced transformation $F:\hat
M\to\hat M$ (cf. Remark \ref{ps.vs.mt}).

If $\mu$ is an $F$-invariant Borel probability measure on $\hat
M$, then there is an associated $\hat f$-invariant Borel
probability measure $\hat\mu$ on $\hat B$ (and vice versa).
Abramov's theorem (\cite{Br}) relates the entropies of both
systems:
 \begin{equation}\label{umum}
h_{\hat\mu}(\hat f)=h_\mu(F)\cdot\hat\mu(\hat M).
 \end{equation}
One does not need to assume ergodicity and can allow degenerations
\cite{DGS}, as happens for the case of Zhang model. In ergodic
situation by the recurrence theorem of Kac \cite{Br} for a
$\mu$-generic point $\hat x\in\hat M$:
 \begin{equation}\label{mumu}
\frac1{\hat\mu(\hat M)}=\lim_{n\to\infty}\frac1n\sum_{k=0}^{n-1}
\tau(F^k\hat x),
 \end{equation}
where $\tau(\hat x)=\tau(i,x)$ is the avalanche time initiated by
addition of $e_i$ to $x\in M$ (see \S\ref{randex}). If
$\mu=\mu_\text{SRB}$ is a unique SRB-measure, the above point
$\hat x$ can be chosen Lebesgue generic. The resulting limit is
the average avalanche time $\bar\tau$ (we discuss it in more
details in \S\ref{Ssec}-\ref{S-sec}) and we obtain:
 $$
h_{\hat\mu}(\hat f)=h_{\mu}(F)/\bar\tau\quad\text{ resp. }\quad
h_{\hat\mu_\text{SRB}}(\hat f)=h_{\mu^u}(\sigma_N^+)/\bar\tau.
 $$

In general non-ergodic situation to get equality (\ref{mumu}) we
should integrate the terms in right-hand side and then we again
obtain the average avalanche size $\langle\tau\rangle$, but now it
is the space-average. Substituting this into (\ref{umum}) we get:
 \begin{equation}\label{mu-tau}
h_{\hat\mu}(\hat f)=h_\mu(F)/\langle\tau\rangle.
 \end{equation}

For SRB-measures this formula is indicated by the Ledrappier-Young
theorem \cite{LY}, because we have only one positive Lyapunov
exponent.


\section{\hps Topological entropy}\label{sec_4}

To calculate the topological entropy of $F$ we established in
\cite{KR2} a set of inequalities, using the technique developed by
J. Buzzi \cite{B1}, \cite{B2} for piecewise expanding maps and
piecewise isometries, see Appendix \ref{app_A}. The contraction in
the maps $F_{ij}$ provides difficulties, so several results were
generalized to fit the framework of the Zhang model. It is not
however true (as was widely believed, see  \cite{B2}) that the
contraction does not contribute to topological (contrary to
metric) entropy, the corresponding counter-example can be found in
\cite{KR1}. The Zhang model has a feature common to all such
examples \cite{KR2}, namely angular expansion, but still for most
values of the parameters this abnormal increase of the entropy
does not occur.


\subsection{\hpss Growth of the number of continuity domains}

Let $\mathcal P=\{[i] \times M_{ij}\}$ be the partition of
continuity for $F$, and enumerate the elements so that $\mathcal
P=\{P_1,\dots P_r\}$. Let
 $$
[P_{a_0}\dots P_{a_{n-1}}]:=\bigcap_{m=0}^{n-1} F^{-m}(P_{a_m})\,,
 $$
and
 $$
\mathcal P^n=\{[P_{a_0}\dots P_{a_{n-1}}] \neq \emptyset\,|\,
a_m=1,\dots,r\}\,.
 $$
We define the singularity entropy of $F$ by
 $$
\hs(F)=\lim_{n \rightarrow \infty} \frac{1}{n} \log
\text{card}(\mathcal P^n)\,.
 $$
 \begin{rk}\po\label{rk9}
Define a map $g:{\mathcal A} \rightarrow \Sigma_r^+$ by letting
$g(\hat{x})$ be the unique sequence $a_0a_1\dots$ such that
$F^n(\hat{x}) \in P_{a_n}$. Then
 $$
\hs(F)=\h(\sigma_r^+;g({\mathcal A})) \,.
 $$
For piecewise affine expanding maps and piecewise isometries it is
clear that this also equals the topological entropy, but due to
the contraction, this is not obvious in the Zhang model. In
addition, if singularities are not removable, the map $g$ has
discontinuities.
 \end{rk}

 \begin{dfn}\po \label{d4}
Call a point $x\in S(F)$ an unstable singularity if for all $i$
and $k\ne l$ we have: $\lim_{y\to x}F_{ik}(y)\ne \lim_{y\to
x}F_{il}(y)$.
 \end{dfn}

 \begin{theorem}\po \label{ssthm}
If all singularities $S(F)\cap {\mathcal Y}$ are unstable, then
$\h(F)=\hs(F)$.
 \end{theorem}

We need the following technical lemma:
 \begin{lem} \po \label{sslemma}
If the singularities in $\Y$ are unstable, then there exists a
constant $\gamma>0$ such that for all $\delta>0$, $x\in\Y\cap
M_{ik}$ and $y\in\Y\cap M_{il}$, $k \neq l$, we have:
 $$
d(x,y) < \delta \Rightarrow d(F_{i}(x),F_{i}(y))>\delta\,.
 $$
 \end{lem}

\noindent
\begin{proof}
Suppose that for all $\delta>0$ there exists $x\in M_{ik} \cap \Y$
and $y \in M_{il} \cap \Y$ such that $d(x,y)<\delta$ and
$d(F_{ik}(x),F_{il}(y))\leq \delta$. There exist sequences
$\{x_m\} \subset M_{ik}$ and $\{y_m\} \subset M_{il}$ such that
$d(x_m,y_m) \rightarrow 0$ and $d(F_{ik}(x_m),F_{il}(y_m))
\rightarrow 0$. The sequence $\{x_m\}$ has a convergent
subsequence $x_{m_n} \rightarrow z$. The point $z$ lies in $S(F)$
and $y_{m_n} \rightarrow z$. By the continuity of the maps
$F_{ik}$ and $F_{il}$ we have $F_{ik}(x_{m_n}) \rightarrow
F_{ik}(z)$ and $F_{il}(y_{m_n}) \rightarrow F_{il}(z)$. Since the
metric $d$ is continuous on $M \times M$ we have
 $$
d(F_{ik}(z), F_{il}(z))=
\lim_{n \rightarrow \infty}
d(F_{ik}(x_{m_n}),F_{il}(y_{m_n}))=0\,.
 $$
Hence $F_{ik}(z)=F_{il}(z)$. Contradiction.
\end{proof} \\

 \noindent {\bf Proof of Theorem \ref{ssthm}.}
Let
 $$
[P_{a_0}\dots P_{a_{n-1}}]=[i_0\dots i_{n-1}]\times K\,,
 $$
where $K \subset M$ is a convex polygon. Fix $\delta>0$ and set
$k=[\log1/\delta]$. Let $z_1,\dots,z_{m(\delta)}$ be a
$\delta$-spanning set for $K$. Chose ${\bf t} \in \Sigma_N^+$ and
define $N^k$ sequences
 $$
{\bf s}_{r_0 \dots r_{k-1}}= (i_0 \dots i_{n-1} r_0 \dots r_{k-1}
t_{n+k} t_{n+k-1} \dots ) \in \Sigma_N^+\,.
 $$
Since the maps $F_{ij}$ are contracting it is clear that the set
$\{\hat{x}=({\bf s}_{r_0 \dots r_{k-1}},z_l)\}$ is a
$(n,\delta)$-spanning set for $[P_{a_0}\dots P_{a_{n-1}}]$, and
since the minimum number of balls needed to cover a convex polygon
$K\subset M$ is bounded by $m(\delta) \leq C_N\delta^{-N}$ we see
that the number of $(n,\delta)$-balls to cover $[P_{a_0}\dots
P_{a_{n-1}}]$ is bounded by $m(\delta)N^k\cdot\#\{[P_{a_0}\dots
P_{a_{n-1}}]\}$. Therefore we get an estimate for the number of
$(n,\delta)$-balls to cover $\hat M$ and so
 $$
\h(F)\leq\lim_{n\to\infty}\frac{\log\bigl(
C_\delta\delta^{-N}\op{card}(\mathcal P^n)\bigr)}{n}=\hs(F)\,.
 $$

To see the opposite inequality let $A,B \in {\mathcal P}^n$ and
take $\hat{x}_1=({\bf t},x) \in A$ and $\hat{x}_2=({\bf s},y) \in
B$. Suppose that $A \neq B$ and that
$t_0=s_0,\dots,t_{n-1}=s_{n-1}$. Then there is $m<n$ such that
$\pi_s \circ F^m(\hat{x}) \in M_{t_mk}$ and $\pi_s \circ
F^m(\hat{y}) \in M_{t_ml}$, $k \neq l$. By Lemma \ref{sslemma}
there is $\gamma>0$ such that for all $\xi < \gamma$:
 $$
\max\{d(\pi_s \circ F^m(\hat{x}),\pi_s \circ F^m(\hat{y})),
d(\pi_s \circ F^{m+1}(\hat{x}),\pi_s \circ F^{m+1}(\hat{y})) \}
\geq \xi\,.
 $$
Therefore for $\delta$ sufficiently small, no $(n+1,\delta)$-ball
can contain points of both $A$ and $B$ and the minimal
$(n+1,\delta)$-spanning set has at least $\op{card}(\mathcal P^n)$
elements. Then $\h(F) \geq \hs(F)$. \qed

 \weg{
\begin{rk}\po
The claim $\h(F)=\hs(F)$ holds in a more general situation, for
instance, when the system is a skew-product $\Sigma_N^+\times
M$, $F({\bf t},x)=(\sigma_N^+{\bf t},f_{t_0}x)$ as in the Zhang
model. The essential feature of this extension is that the
singularities of the unstable factor do not intersect under
backward-iterations.
\end{rk}
 }

 \begin{theorem}\po \label{nostable}
For the Zhang models: $\h(\hat f)=\hs(\hat f)$, $\h(F)=\hs(F)$.
 \end{theorem}

\noindent
 \begin{proof}
By Remark \ref{rk9} all quantities are topological entropies. The
corresponding systems in the second equality are the Poincar\'e
return maps for the transformations of the first equality. The map
$F:\hat M\to\hat M$ is already the return map for $\hat f:\hat
B\to \hat B$ by the very construction.

To achieve the same claim for the symbolic system we extend the
partition $\mathcal{P}$ of $\hat M$ to a partition
$\mathcal{\tilde P}$ of $\hat B$, which on
$\Sigma_N^+\times(\R^N_{\ge0}\setminus M)\cap\hat B$ equals the
product of the standard partition $\Sigma_N^+=\cup_i[i]$ and the
partition of the spacial part by the hyperplanes $\{x_i=E_c\}$.
Denote by $s\le r+N(2^N-1)$ the number of elements of the new
partition $\mathcal{\tilde P}$.

Let $\mathcal{B}\subset\hat B$ be the $\hat f$-invariant closure
of $\mathcal{A}$ in $\Sigma_N^+\times\R^N_{\ge0}$. Define a map
$\hat g:\mathcal{B}\to\Sigma_s^+$ by letting $\hat g(\hat x)$ be
the unique sequence $b_0b_1\dots$ such that $\hat f^n(\hat
x)\in\mathcal{\tilde P}_{b_n}$. We wish to prove that
 \begin{equation}\label{tten}
\h(\hat f)=\hs(\hat f)=\h(\sigma_s^+|\hat g({\mathcal B})) \,.
 \end{equation}
For this it is sufficient to check that the singularities $S(\hat
f)\cap\pi_s(\mathcal{B})$ are unstable for $\hat f$ (the second
equality follows from the definition).

Consider a singular point $x\in\R^N_{\ge0}$. Let $y,z$ tend to $x$
by two different domains of the projected partition
$\mathcal{\tilde P}$ in $\R^N_{\ge0}$. Then we can subdivide
$\{1,\dots,N\}=A\cup B\cup C$, where $y_i=E_c+0$, $z_i=E_c-0$
(with the obvious notations instead of limits) for $i\in A$,
$y_i=E_c-0$, $z_i=E_c+0$ for $i\in B$ and $y_i,z_i$ belong to the
same side of $E_c$ for $i\in C$. Then for $\tilde S=S(y)-S(z)$ we
have:
 $$
\tilde S_{ij}=\left\{
 \begin{array}{ll}
\e-1, & \text{ if }i=j\in A,\\
1-\e, & \text{ if }i=j\in B,\\
\frac{1-\e}{2d}, & \text{ if }d_\Lambda(i,j)=1, j\in A,\\
\frac{\e-1}{2d}, & \text{ if }d_\Lambda(i,j)=1, j\in B,\\
0 & \text { otherwise.}
 \end{array} \right.
 $$
We should check that $\tilde S\cdot v\ne0$ for any vector
$v\in\R^N_{\ge0}$ with components $v_i=E_c$ for $i\in A\cup B$
(this implies that $x$ is unstable). Note that other components
$v_i$, $i\in C$ do not contribute to the product.

Let $i$ be a site from $A$ with not all neighbors from $A$ (if the
set $A$ is empty consider $B$). Denote the number of $A$-neighbors
of $i$ by $k_A<2d$ and the number of $B$-neighbors by $k_B\ge0$.
Then $(\tilde S\cdot v)_i=E_c(2d-k_A+k_B)\frac{\e-1}{2d}\ne0$.

Thus Theorem \ref{ssthm} implies (\ref{tten}). Moreover, the same
reasons yield a more general statement. Namely, since $\hat g$ is
a semi-conjugancy, we get:
 \begin{equation}\label{sten}
\h(\hat f|\mathcal{K})=\h(\sigma_s^+|\hat g({\mathcal{K}}))
 \end{equation}
for any subset $\mathcal{K}\subset \mathcal{B}$ (recall that $\hat
g$ may have discontinuities, but our arguments are not injured by
this fact). This subset needs not to be invariant, and in this
case we should use Bowen's definition of entropy \cite{Bo}. Thus
Katok's entropy formula \cite{K} implies that $h_{\hat\mu}(\hat
f)=h_{\hat g_*\hat\mu}(\sigma_s^+)$ for all $\hat f$-invariant
measures $\hat\mu$, which are not supported on singularities.

If the system $(\hat M,F)$ possesses a measure $\mu$ of maximal
entropy, then by the obtained result the second claim of the
theorem follows from (\ref{mu-tau}) and Kac's theorem \cite{PY}.
In general, we can apply the above arguments to the partition of
$\hat M$ by the subsets of equal return times and using the fact,
that both returns of $\hat f$ and $\hat g$ have the same
combinatorics, we get: $\h(F)=\h(\sigma_r^+|g({\mathcal A}))$
(see, for instance, the loop equation approach \cite{Pt}).
 \end{proof}


\subsection{\hpss Evaluation of topological entropy} \label{evalofentr}

It was predicted on the base of variational principle in
\cite{BCK} that topological entropy of the Zhang model is
$\h(F)=\log N$. However this principle does not apply because the
map is not well-defined (continuously) on the whole space (or
thanks to non-compactness if we remove the singularities). In
fact, there can be no invariant measures on the non-singular part
at all.

While we support the claim that $\h(F)=\log N$, it will not be
proved in full generality. We start with the asymptotic
statement.

\begin{figure}
\begin{center}
\includegraphics[width=5.5cm]{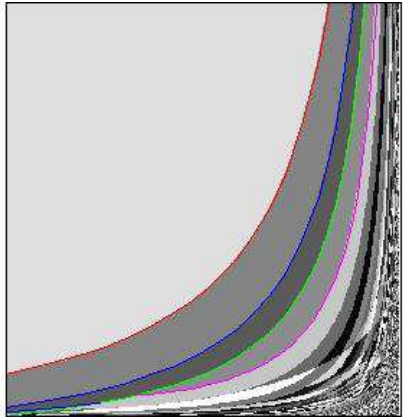}
\caption{{\small Shows the avalanche type domains on
$(\e,E_c)$-bifurcation diagram for $N=2$. The top domain is
$E_c>\frac{1+\e}{1-\e}$. The line from infinity to the origin is
$E_c=\frac{\e}{1-\e}$. We see infinitely many domains with
different avalanches, accumulating to two of the axes.}}
\label{bifdiag}
\end{center}
\end{figure}

Consider the bifurcation diagram on the $(E,\e)$ strip
$(0,+\infty)\times[0,1)$, where a point is critical if in its
neighborhood dynamics of the Zhang model can experience avalanches
of different types. Thus the strip is partitioned into different
avalanche type domains. The partition depends on $L,d$. For $N=2$
the diagram is shown on Figure \ref{bifdiag}.

\noindent

Note that for all $L,d$ there is the top avalanche type domain
representing the shortest avalanche time. For $N=2$ it is given by
the relation $E>\frac{1+\e}{1-\e}$. Also note that some domains
have $\e$-projection strictly smaller than the interval $[0,1)$.
The next statement concerns only the top domain and the avalanche
type domains that are adjacent to the line $\e=1$.

 \begin{theorem}\po \label{140}
For the Zhang model: $0\le\h(F)-\log N\le\theta(E,\epsilon)$,
where $\lim_{E_c\to\infty}\theta(E_c,\epsilon)=0$ ($\e$ fixed) and
$\lim_{\e\to1}\theta(E_c,\epsilon)=0$. In the latter case $E_c$
changes accordingly with $\e$ so that $(E_c,\e)$ belongs to the
same avalanche type domain (then $E_c\to\infty$, though this case
differs from the former).
 \end{theorem}

We can remove the $N$-rational points from $\Sigma_N^+$ and then
represent it as the subset of $I=[0,1]$ with points $n/N^k$ being
deleted (this process does not change the topological entropy). In
fact, we need to remove only points $n/N$ for the others will be
deleted by inverse iterations of the map of $I\times X$ with the
formula $F(t,x)=(Nt,f_{[Nt]}(x))$. Thus we represent our system as
a piece-wise affine partially hyperbolic map of the subset of
$\R^{1+N}$. This is important for an application of the results
from Appendix \ref{app_A}.

The proof uses the notions of multiplicity and multiplicity
entropy, due to G. Keller and J. Buzzi. Given a finite partition
${\mathcal Z}=\{Z_k\}$ of $\hat{M}$ we define
 $$
{\mathcal Z}^n=\{[Z_0\dots Z_{n-1}] \neq \emptyset\}\,,
 $$
where
 $$
[Z_0 \dots Z_{n-1}]=Z_0 \cap F^{-1}(Z_1) \cap \dots \cap
F^{1-n}(Z_{n-1})\,.
 $$
Define multiplicity of a partition by
 $$
\text{mult}({\mathcal Z})=\max_{\hat{x} \in \hat{M}}
\text{card}\{Z \in {\mathcal Z}\,|\,\overline{Z} \ni x\}\,,
 $$
If ${\mathcal Z}$ is the continuity partition of the map $F$ we
often denote the multiplicity of the ${\mathcal Z}$ by
$\text{mult}(F)$. Then it is clear that
$\text{mult}(F^n)=\text{mult}({\mathcal Z}^n)$. The multiplicity
entropy of $F$ is (the limit exists by subadditivity, cf.
\cite{KH})
$$
H_{\text{mult}}(F) = \sup_{{\mathcal Z}} \lim_{n \rightarrow
\infty} \frac{1}{n} \log \text{mult}({\mathcal Z}^n)\,,
$$
and we see that if ${\mathcal Z}$ is the continuity partition,
then
$$
H_{\text{mult}}(F) = \lim_{n \rightarrow \infty} \frac{1}{n} \log
\text{mult}(F^n)\,.
$$

It is clear that $\hm(F)=0$ if the singularities are removable, so
$\h(F)=\log N$ by the result in Appendix \ref{app_A}. For big
$E_c\gg1$ the singularities are generally non-removable, but still
we have the same effect asymptotically:\\

\noindent
 \begin{Proof}{Proof of Theorem \ref{140}}
We take $\theta(E_c,\e)=\hm(F)$. Since the singularities $t=n/N\in
I$ of the map $\sigma_N^+:t\mapsto Nt$ do not intersect in inverse
iterations, the multiplicity growth is only due to the spacial
maps $F_i:M\to M$. Thus using the notation $F_{\bf t}$ from
Section \ref{meaofmaxentr} for the dynamics over a prefixed
sequence of excitations we obtain $\hm(F)=\sup_{{\bf
t}\in\Sigma_N^+}\hm(F_{\bf t})$.

To show the first claim let us notice that when $\e=\op{const}$,
but $E_c\to\infty$, then the avalanches map the critical part of
the boundary $\{\exists i:x_i=E_c\}\subset\partial M$ far from
$\partial M$, namely to the distance $\sim\gamma(\e,L,d)E_c$, see
Figures \ref{pic1} for $N=2$ and \ref{pic2} for $N=3$. To reach
again the boundary and experience avalanche we need many shifts
$F_i$ by the basic vectors $e_i$. Thus the singularity can meet
only after big number of iterations. Since the initial picture of
singularities has bounded multiplicity (for $(\e,E_c)$ from the
top avalanche type domain), the multiplicity decreases at least as
$k/E_c$ so that it vanishes in the limit.

To prove the second claim we use the inequality
$\hs(F)\le\sum\rho_i(F)$ from Appendix \ref{app_A.3}. If the
avalanche type domain is fixed, the number of compositions of
matrices $S(x)$ in one avalanche (see \S\ref{return maps}) is
bounded. Every such a matrix tends to identity when $\e\to1$. Thus
all the linear parts $L_{ij}$ of avalanche maps $F_{ij}$ tend to
identity and so angular expansions $\rho_i(F)$ tend to zero.
 \qed
 \end{Proof}
 \vskip6pt

Notice that if $E_c$ is fixed, but $\e\to1$, then the number of
avalanches has unlimited grow and the previous argument do not
work. However, due to estimates of \S\ref{randex} on the maximal
avalanche length $\tau_m$ we conclude that $\theta(E_c,\e)\to0$ as
either $E_c\to\infty$ or $\e\to1$, both quantities being related
by the constraint $E_c\ge
N(1-\e)^{-\frac12\op{diam}(\Lambda)-\sigma}$ for some $\sigma>0$
(we don't require, but allow $N\to\infty$ as well). This statement
is stronger than in the theorem.

Now we are going to prove vanishing of $\theta(E_c,\e)$ a.e. in
the finite part.

 \begin{theorem}\po \label{a.s.0}
For generic $(E_c,\e)$ we have: $\h(F)=\log N$.
 \end{theorem}

We will prove the statement not only for the Zhang model, but also
for nearly Zhang models. By this we mean the following. The map
$F$ is a bundle over the Bernoulli shifts $\sigma_N^+$ with
factors $F_i$ being piece-wise affine partially contracting maps.
That is $M=\cup_j M_{ij}$ and for $F_{ij}=F_i|_{M_{ij}}$ we have
$F_{ij}(x)=L_{ij}(x)+b_{ij}$, $b_{ij}=L_{ij}e_i$. We are going to
make arbitrary small generic perturbations of the matrices
$L_{ij}$ (still with spectrum within unit ball) and vectors
$b_{ij}$ (one should care that $M_{ij}$ are mapped into $M$) and
prove the statement for this modified system.

Due to round-off errors there is no much difference between the
original and the perturbed systems in computer simulations. And
from the point of view of the experiment such perturbations
(instrument instability) are indispensable -- look at \cite{Ru}
for the discussion of physical relevance of variation of
parameters as the noise.

 \vskip4pt
\noindent
 \begin{Proof}{Proof of Theorem \ref{a.s.0}}
We claim that for generic $(E_c,\e)$ the singularities do not
multiply. Actually, some intersections of singularities are
deformable as we vary the parameters, just by the transversality
reasons, but the other disappear with small perturbations.

Namely, for a multi-index $\sigma=(\alpha_1,\dots,\alpha_k)$,
$\alpha_s=(i_s,j_s)$ coding an orbit, denote
$F_\sigma=F_{\alpha_k}\circ\dots\circ F_{\alpha_1}$ the
corresponding map along the orbit. The singularities of this map
are: $\op{Sing}(F_\z)=\cup_{r=1}^k
F_{\z_{[r]}}^{-1}\op{Sing}(F_{\alpha_r})$, where
$\z_{[r]}=(\alpha_1,\dots,\alpha_r)$.

If $z_0\in\op{Sing}(F_\z)$, then its orbit (with
multi-possibilities due to singularities: mapping a singular point
we extend the components of the map $F_\z$ in various ways) meets
several singularity planes, i.e. for some cuts $\z_s=\z_{[r_s]}$
of $\z$, $s=1,\dots,m$, we have the following system:
 \begin{equation}\label{masha}
F_{\sigma_1}(z_0)=z_1,\dots,\,F_{\sigma_m}(z_0)=z_1,\
l_{q_1}(z_1)=0,\dots,\,l_{q_m}(z_m)=0,
 \end{equation}
where $l_{q_s}(z)$ are equations for the singularity hyperplanes
of the corresponding map $F_i$ (there are also inequalities, which
we don't mention). When $m\le N$ there are occasions, when
(\ref{masha}) has a solution continuously depending on $(E_c,\e)$.
However for $m>N$ this is no longer the case. In fact, considering
nearly-Zhang models we see that for generic data the above system
(\ref{masha}) is characterized by a collection of non-trivial
polynomial in $\e$ equations (for each $E_c$).

More precisely, the set of $(L_{ij},b_{ij})$ giving trivial
polynomials has positive codimension and hence zero Lebesgue
measure. Uniting these sets over all choices of multi-indices
$(\z_1,\dots,\z_m)$ we see that the complement has full measure
and dimension and so a generic perturbation yields the data
$(L_{ij},b_{ij})$ from it. Since non-trivial polynomials have only
finite number of zeros, then for a generic nearly Zhang model and
every $E_c$ there is a countable subset of $\{0\le\e<1\}$, so that
the corresponding systems (\ref{masha}) have no solutions. This
means that multiplicity of $F_{\bf t}^n$ does not grow with $n$
and so $\hm(F_{\bf t})=0$.

Now let us look to the Zhang model. We restrict for simplicity of
exposition to the case of two sites $N=2$. In this case the linear
parts of the affine maps are compositions of $L_1=\1+(\e-1)A_1$
and $L_2=\1+(\e-1)A_2$ with
 $$
A_1=\begin{bmatrix}1 & 0 \\ -1/2 & 0\end{bmatrix},\quad
A_2=J^{-1}A_1J=\begin{bmatrix}0 & -1/2 \\ 0 & 1\end{bmatrix},
\text{ where }J=\begin{bmatrix}0 & 1 \\ -1 & 0\end{bmatrix}
 $$
(we exclude the obvious matrix $\1$, see Example A for details).
Notice that $\det L_1=\det L_2=\e$.

Suppose that (\ref{masha}) has a continuous solution in some
domain of $(E_c,\e)$. Then it is algebraic in $\e$ and linear in
$E_c$ (the latter is because the singularity lines within one
avalanche type domain shift with velocity 1 in the direction of
either $e_1=\begin{pmatrix} 1\\0 \end{pmatrix}$ or
$e_2=\begin{pmatrix} 0\\1 \end{pmatrix}$. Differentiating this by
$E_c$ we obtain:
 \begin{equation}\label{medved}
L_{\rho_1}(z_0')=z_1',\dots,L_{\rho_m}(z_0')=z_m',
 \end{equation}
where $L_{\rho_i}=L_{\rho_{i,t_i}}\circ\dots\circ L_{\rho_{i,1}}$
for $\rho_i=(\rho_{i,1},\dots,\rho_{i,t_i})$ is the linear part of
$F_{\z_i}$ ($\rho_i$ is different from $\z_i$ because $dF_{\z_i}$
is a composition of several maps $L_s$). The points $z_k$, $1\le
k\le m$, are constrained to the singularity lines and we can
suppose these are the lines $\{x_1=E_c\}$ or $\{x_2=E_c\}$ (all
other singularities are mapped to them within one avalanche). Thus
$z_k'=v_k+\psi_k w_k$, where $v_k=e_1\text{ or }e_2$ and
$w_i=Jv_i$ is the other basic vector, $\psi_k$ being an unknown
scalar.

Let $\xi_k=L_{\rho_k}^{-1}(v_k)$, $\zeta_k=L_{\rho_k}^{-1}(w_k)$;
these vectors depend meromorphically on $\e$. System
(\ref{medved}) is solvable iff the affine lines
$\xi_k+\psi_k\zeta_k$, $1\le k\le m$, in $\R^2$ have a common
point. We can suppose $m=3$.

Denote by $\Omega(\xi,\eta)=\langle J\xi,\eta\rangle$ the standard
symplectic form on $\R^2$. The above 3 lines intersect jointly iff
 $$
\Omega(\xi_1,\zeta_1)\Omega(\zeta_2,\zeta_3)+
\Omega(\xi_2,\zeta_2)\Omega(\zeta_3,\zeta_1)+
\Omega(\xi_3,\zeta_3)\Omega(\zeta_1,\zeta_2)=0.
 $$
Dividing by $\prod_{k=1}^3\Omega(\xi_k,\zeta_k)$ and using the
fact that
$\Omega(\xi_k,\zeta_k)=\det(L_{\rho_k}^{-1})=\e^{-|\rho_k|}$ we
get the equivalent equation:
 \begin{equation}\label{yest}
\Omega(\eta_1,\eta_2)+\Omega(\eta_2,\eta_3)+\Omega(\eta_3,\eta_1)=0.
 \end{equation}
Here $\eta_k=\e^{-|\rho_k|}L_{\rho_k}^{-1}w_k=\tilde
L_{\rho_k}w_k$, where $\tilde L_{\tau}=\tilde
L_{\tau_1}\circ\dots\circ\tilde L_{\tau_t}$ for
$\tau=(\tau_1,\dots,\tau_t)$ and $\tilde L_k=\e^{-1}L_k^{-1}=
\1+(\e-1)\tilde A_k$ is the adjunct matrix for $L_k$, which gives
 $$
\tilde A_1=\begin{bmatrix}0 & 0 \\ 1/2 & 1\end{bmatrix},\quad
\tilde A_2=\begin{bmatrix}1 & 1/2 \\ 0 & 0\end{bmatrix}.
 $$

Now equation (\ref{yest}) holds iff there exist 3 not
simultaneously zero numbers $\beta_1,\beta_2,\beta_3$ such that
$\beta_1+\beta_2+\beta_3=0$ and
$\beta_1\eta_1+\beta_2\eta_2+\beta_3\eta_3=0$. Since $\eta_k$ is a
polynomial matrix of degree $|\rho_k|$ and the products of $A_t$
are always proportional to $e_1$ or $e_2$ (depending on the
left-most factor) this last equation is never satisfied if the
multi-indices $\rho_k$ are different.
 \qed
 \end{Proof}

 \begin{rk}\po
To support the usage of nearly-Zhang models note that the whole
paradigm of SOC should allow generic perturbations of the data,
for if there is a fine tuning of parameters, the model is
unappropriate for physical explanation (of course, we should pass
to the thermodynamic limit, but in practice this only means some
large finite parameters).
 \end{rk}

Our computer experiments did not expose any exponential growth of
multiplicity in the Zhang model (though we see growth in
complications of singularities), so we suggest that $\hm(F)=0$ and
hence $\h(F)=\log N$ always. In addition, by the above discussion
we can disregard these exceptional values of $(E_c,\e)$ even if
there are any. This finishes discussion of topological entropy.


\section{\hps Geometry of the attractor}\label{sec_5}

The construction of the spacial attractor $\Y$ can be interpreted
as an iterated function system (IFS), where the maps $F_i$ are not
affine as usually considered, but piecewise affine. Hence one
might expect that attractors $\Y$ are fractal, but with various
size characteristics, like dimension and measure, depending on
parameters $E_c$ and $\epsilon$.


\subsection{\hpss Fractal structure}

\begin{figure}
\begin{center}
\includegraphics[width=5.5cm]{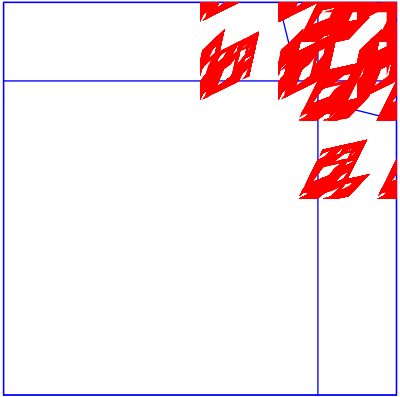}
\caption{{\small  Shows the set $U_{10}$ for $N=2$, $E_c=5$ and
$\epsilon=1/2$.}} \label{frac1}
\end{center}
\end{figure}

\begin{figure}
\begin{center}
\includegraphics[width=5.5cm]{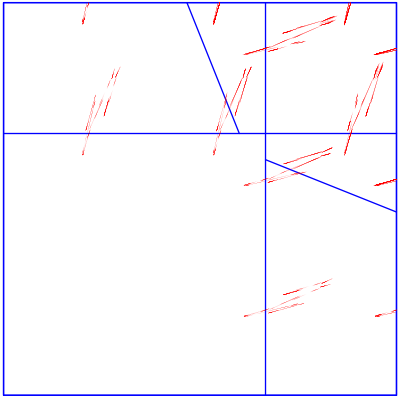}
\caption{{\small  Shows the set $U_{10}$ for $N=2$, $E_c=3$ and
$\epsilon=1/5$.}} \label{frac2}
\end{center}
\end{figure}

Computer experiments show that in certain cases the spacial
attractor $\mathcal{Y}$ has fractal structure, see e.g. Figures
\ref{frac1} and \ref{frac2}. We clearly see that $\Y$ consists of
self-similar pieces. However the pieces overlap, making evaluation
of the fractal dimension difficult. So we can provide only
estimates of the attractor's size.

Nevertheless we observe from our experiments that Hausdorff
dimension $\dim_\text{H}(\Y)$ and the Lebesgue measure
$\mu_\text{Leb}(\Y)$ of the attractor grow piece-wise
monotonically with $E_c$ and $\e$. Thus the following effects
occur in steps:

 \begin{itemize}
\item The dimension and the measure of $\Y$ vanish.
\item $\dim_\text{H}(\Y)$ is positive, while the measure is zero.
\item Both $\dim_\text{H}(\Y)$ and $\mu_\text{Leb}(\Y)$ are
positive.
\item The attractor $\Y$ contains an interior point.
 \end{itemize}

We will demonstrate the dimensional part in the next section,
while we disregard the observation about the measure. The reason
for this is that $\mu_\text{Leb}$ is not physically motivated and
we should look for an SRB-measure.

The experiments show that such a measure exists and has support
lying strictly inside $\Y$. (See Figure \ref{nonfrac} and Figure
\ref{orbit}. In Figure \ref{nonfrac} the attractor $\Y$ is shown
for the case $N=2$, $E_c=20$ and $\e=2/3$, and Figure \ref{orbit}
shows the orbit of a random initial condition. The latter
corresponds to the support of the SRB-measure.) This is possible
because the contraction rate of $f_{\bf t}$ is smaller for the
exceptional sequences ${\bf t}\in\Sigma_N^+$, than for a generic
one. Thus study of the IFS-attractor does not lead to conclusions
about ergodicity or uniqueness of the SRB-measure. Still it
provides an information about spacial distribution of the orbits
in the Zhang dynamics.

\noindent
\noindent


\subsection{\hpss Dimensional study of the attractor}

\begin{figure}
\begin{center}
\includegraphics[width=5.5cm]{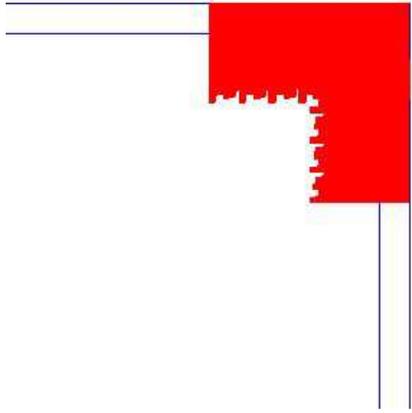}
\caption{{\small  Shows the set $\Y$ in the top right corner of
$M$ for $N=2$, $E_c=20$ and $\epsilon=2/3$. The dimension of $\Y$
is $2$ in this example.}} \label{nonfrac}
\end{center}
\end{figure}

\begin{figure}
\begin{center}
\includegraphics[width=5.5cm]{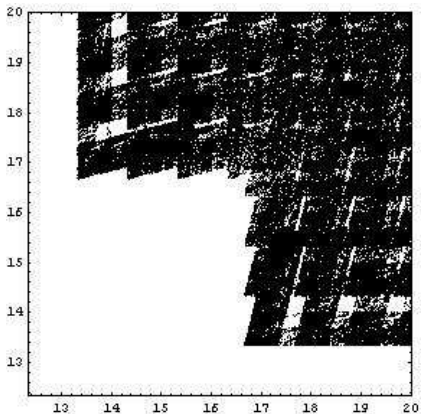}
\caption{{\small  Shows an orbit of a randomly chosen point (i.e.
the support of the SRB-measure) for $N=2$, $E_c=13$ and
$\epsilon=2/3$.}} \label{orbit}
\end{center}
\end{figure}

The fractal properties of $\Y$ do not hold for all values of
parameters $(E_c,\e)$. An example where $\Y$ has integer dimension
is shown in Figure \ref{nonfrac}.

It was noted in \cite{BCK} that Hausdorff (fractal) dimension of
the attractor is about to increase as $E_c$ grows. The arguments
were the following: For bigger $E_c$ the contraction rate
decreases, so the theory of iterated function system (IFS) implies
increasing of the Hausdorff dimension
 $$
{\frak D}_\mathcal{Y}(\e,E_c)=\dim_\text{H}(\mathcal{Y})
 $$
as a function of $E_c$. While this seems to be true, the statement
does not hold in precise sense. For instance, for $N=2$, $E_c\in
[\frac{1+\e}{1-\e},\frac{2}{1-\e}]$ the attractor is the set of 3
points, while it seems to have non-zero dimension for other
parameters (computer simulations clearly show this).

\noindent
\noindent

The problem is with the framework of IFS, where usually only
conformal maps are considered and certain regularity of their
mapping graph and overlaps is assumed. However we will show
validity of the claim in the asymptotic sense:

 \begin{theorem}\po\label{theo-5}
With fixed $d,L$ and generic $\e$ we have:
$\lim_{E_c\to\infty}{\frak D}_\mathcal{Y}(\epsilon,E_c)=N$.
Moreover, ${\frak D}_\mathcal{Y}(\epsilon,E_c)=N$ for big values
$E_c\gg1/\e$.

On the other hand for all $\e\in[0,1)$ it holds:
$\lim_{E_c\to0}{\frak D}_\mathcal{Y}(\epsilon,E_c)=0$.
 \end{theorem}

In the above statement "generic" means both full Lebesgue measure
and second Baire category. In fact, the equality holds for all
$\e$ outside a countable set. It seems though that the limit
statement is valid for all $\e$. Thus we see that ${\frak
D}_\mathcal{Y}(\e,E_c)$ is not strictly monotone in $E_c$ for its
large values as one might expect from the arguments cited before
the theorem. \\

 \noindent
\begin{proof}
Consider first the statement about big energies $E_c\gg1/\e$. Let
us start by demonstrating the idea of the proof on the example
$N=2$, see Figure \ref{pic1}.

The image of the vertical continuity domain $M_{13}$ of height 1
adjusted to the right-top corner is a trapezium with the slope
depending on $\e$ (see (\ref{3lines}) of Section \ref{examples}
for the numeration of domains). It is thin -- of constant length
of horizontal section equal $\e$ near its bottom side, but long --
with diameter approximately equal $E_c\frac{1-\e^2}4$. Thus if we
shift $\sim C/\e$ times this domain up and then all of the shifted
images horizontally to the right, so that its first coordinate
satisfies the inequality $E_c-1<x_1\le E_c$, then these shifts
cover an open domain, including a unit square, in the vertical
strip $K_1=\{x_1\in(E_c-1,E_c]\}$ (note that we can leave a copy
of the domain since this corresponds to the shift -- dropping of
energy $(x_1,x_2)\mapsto(x_1,x_2+1)$ on $M_{13}\subset K_1$ before
the avalanche). An easy calculation shows that such shifts cover
the whole upper part of $K_1$, strictly including the continuity
domain $M_{13}\subset K_1$ adjusted to $(E_c,E_c)$ and covering a
vertical part of $M_{12}$. A similar scenario happen to the second
coordinate. Thus in iterating the dynamics we will always have two
continuity domains $M_{13}$ and $M_{23}$ adjusted to $(E_c,E_c)$
and the adjacent parts of $M_{12},M_{22}$ lying in the attractor.

For general $L$ and $d$ we observe the same picture: With generic
$\e$ a continuity domain $M_{ij}$ adjusted to the upper-most
corner is mapped under avalanche-map $F_{ij}$ to the
trapezoid-like polyhedron with irrational slopes. Its shifts cover
then an open domain in each of the strips
$K_i=\{x_i\in(E_c-1,E_c]\}$ and so after more shifts -- the upper
part of this strip, whence the statement.

We illustrate this process on Figure \ref{pic2}. The 3 domains
adjusted to the corner $(E_c,E_c,E_c)$ are mapped into interior of
the cube and they have different irrational slopes (we picture
them of zero thickness that corresponds to large values of
$E_c\gg1$), so that their shifts cover a big open domain near the
faces adjusted to the above corner.

Consider now the second statement, $E_c\ll1$. To estimate the
fractal dimension from above we use the generalization of the
Moran's formula from Appendix \ref{app_B}. It implies that if the
IFS $f_1,\dots,f_N$ satisfies $\|f_i\|\le\delta$, then the
Hausdorff dimension of the attractor admits the following
estimate: ${\frak D}_\mathcal{Y}\le\log
{N\vartheta}/\log\frac1\delta$, where $\vartheta$ is the maximal
multiplicity of the continuity partitions for $f_i$.

Now we claim that as $E_c\to0$ we have: $\delta=\max\|f_i\|\le
(E_c)^\sigma$ for some $\sigma>0$. To see this let us estimate the
maximal duration of the avalanche
$\bar\tau_m=\max_{(i,x)\in\Lambda\times M}\tau(i,x)$. This
quantity tends to $\infty$ as $E_c\to0$, but not as fast as
$\tau_m(E_c,\e,\Lambda)\sim C_1/E_c$ (see \S\ref{randex}). Namely
we state that $\bar\tau_m\sim C_2\log1/E_c$. Actually, if we drop
energy 1 to an arbitrary site from a configuration in $M$, then in
a finite $E_c$-independent time all the sites become overcritical.
They remain overcritical, while the system does not loose a
substantial amount of energy. During this process the total energy
is dissipating in geometric progression with an average
contraction rate $1-\frac{1-\e}{N}<1$. So the duration of this
stage has asymptotic $C_3\log(1/E_c)$. The remaining time to
finish the avalanche has a smaller asymptotic.

By \S\ref{randex} $n(E_c,\e,\Lambda)<C_0$ for small $E_c$. Thus
the proof of Theorem \ref{5} implies that for certain
$E_c$-independent constant $C_n=kC_0$ for each sequence of $C_n$
steps in the avalanche process the product of the corresponding
$S$-matrices will have norm $\le c(\e)<1$, which is a uniform
estimate in $E_c\ll1$. The number of steps in one avalanche grows
as $\log1/E_c$. Therefore $\delta\le C_4
c(\epsilon)^{C_n^{-1}\log1/E_c}\le\exp(-\sigma\log1/E_c)$, where
$\sigma=\frac12C_n^{-1}\log\frac1{c(\e)}$ as was claimed in the
estimate.

Next we claim that $\vartheta\le\vp(E_c)$ with
$\vp=o((1/E_c)^\upsilon)$ $\forall\upsilon>0$. In fact, we
described above the avalanche process for small energies. The
first stage is finite and contributes only a bounded number of
singularity hyperplanes. In its seconds stage all of the sites are
excited, so there the corresponding number of singularity
hyperplanes equals the duration. The last stage is shorter of time
$\psi=o(\log1/E_c)$, but the number of singularity hyperplanes
grows faster, but still is bounded by $e^{\psi\log N!}
=o((1/E_c)^\upsilon)$ (see \S\ref{return maps}) for any
$\upsilon>0$.

Finally ${\frak D}_\mathcal{Y}(\e,E_c)\le\log
{N\vartheta}/\log\frac1\delta\le\frac{C'+\log\vp}{\sigma\log1/E_c}\to0$
as $E_c\to0$.
\end{proof}

\begin{cor}\po
For big $E_c\gg1/\e$ the system $(\mathcal{A},F)$ is not
topologically transitive.
\end{cor}
\noindent
\begin{proof}
Suppose $({\bf t},x)$ is a point with the dense orbit in
$\mathcal{A}$. We know that for big $E_c$ the Lebesgue measure of
the spatial part $\Y$ of the attractor $\mathcal{A}$ is a positive
number $\omega>0$. It follows from Theorem \ref{5} that under
iterations with fixed excitation sequence $\bf t$ the volume of
the spatial part $M$ decreases in geometric progression with the
number of avalanches. Thus after a finite number of steps it
becomes less than $\omega$. This iteration will be still a finite
number of polyhedra, so that its closure does not coincide with
$\Y$. Since it contains all the points $\pi_s(F^n({\bf t},x))$, we
obtain contradiction.
\end{proof}

\vspace{10pt}

In the case $N=2$ and $\e>1/2$, the value of $E_c$ starting from
which ${\frak D}_\mathcal{Y}(\e,E_c)=2$ can be calculated
precisely because even one shift of the sloped strip mentioned in
the above proof overlaps with itself and is sufficient for
obtaining an open domain in the attractor. This condition $\e>1/2$
together with $E_c\gg1$ from the theorem ideologically coincide
with the sufficient conditions for invertibility of the
differentials of avalanche maps (Proposition \ref{7}). This makes
an indication of a relation between this invertibility and
fractality of the attractor in the spirit of Ledrappier-Young
formula \cite{LY}. This
latter is however unappropriate in our situation. \\

\noindent {\bf Note on the usage of the Ledrappier-Young formula.}
This formula, essentially used in \cite{BCK} in the study of the
Zhang model, cannot be used for the map $F:\D \to \D$ since this
map is never invertible (in loc. cit. it was applied to $F^{-1}$).
In addition to invertibility the Ledrappier-Young theorem is based
on the SRB-property. For the map $F^{-1}$ this property is
equivalent to absolute continuity of the stable foliation for $F$
w.r.t.\ the measure $\mu$. If the measure has fractal support this
cannot happen. Therefore all formulas based on this property may
turn to be wrong. We demonstrate this in Example A of \S
\ref{examples}.

On the other hand, as we have just shown in the theorem, in
thermodynamic limit $E_c\to\infty$ the fractality is lost and so
the absolute continuity property is restored (but only for the
geometric attractor, the support of SRB-measure is smaller!). This
however does not help with non-invertibility of the factor
$(\Sigma_N^+,\sigma^+_N)$. Even if we change this factor to
invertible two-sided sequences $(\Sigma_N,\sigma_N)$, the system
remains non-invertible since not all points of the attractor
(which can be quite fat) admit negative iterations (Remark
\ref{rkk4}). In addition, a new negative Lyapunov exponent $-\log
N$ in the first factor appears and the formulas exploited in
\cite{BCK} become completely inadequate.


\section{\hps Examples} \label{examples}

In the examples below we consider the one-dimensional Zhang
model with two sites, $N=2$. \\

\noindent
{\bf Example A:} A computation shows that for $E_c \geq
(1+\epsilon)/(1-\epsilon)$ we have six domains of continuity $[i]
\times M_{ij}$, $i=1,2$. The domains $M_{1j}$ are given by

 \begin{equation}\label{3lines}
 \begin{array}{l}
M_{11} =\{x\in [0,E_c]^2\,|\,x_1+1\leq E_c \}       \\
M_{12}=\{x\in [0,E_c]^2\,|\,x_1+1> E_c \text{ and }
(1-\epsilon)(x_1+1)/2+x_2 \leq E_c \}\\
M_{13}=\{x\in [0,E_c]^2\,|\,x_1+1>E_c \text{ and }
(1-\epsilon)(x_1+1)/2+x_2 > E_c \}
 \end{array}
 \end{equation}
and the domains $M_{2j}$ are symmetric to these. The maps $F_{ij}$
are of the form $F_{ij}(x)=L_{ij}(x+e_i)$, where
 \begin{center}
\begin{tabular}{lll}
$L_{11}= \begin{bmatrix} 1 &  0  \\ 0 & 1 \end{bmatrix}$ &
$L_{12}= \begin{bmatrix} \epsilon & 0 \\ \frac{1-\epsilon}{2} & 1
\end{bmatrix}$ &
$L_{13}= \begin{bmatrix}  (\frac{1+\epsilon}{2})^2 &
\frac{1-\epsilon}{2}
\\ \epsilon \frac{\vphantom{2^{2^2}} 1-\epsilon}{2} & \epsilon \end{bmatrix}$
\end{tabular}
\end{center}

and
\begin{center}
\begin{tabular}{lll}
$L_{11}= \begin{bmatrix} 1 &  0  \\ 0 & 1 \end{bmatrix}$ &
$L_{12}= \begin{bmatrix} 1 & \frac{1-\epsilon}{2} \\ 0 & \epsilon
\end{bmatrix}$ &
$L_{13}= \begin{bmatrix} \epsilon & \epsilon \frac{1-\epsilon}{2}
\\  \frac{\vphantom{2^{2^2}} 1-\epsilon}{2} & (\frac{1+\epsilon}{2})^2 \end{bmatrix}$
\end{tabular}
\end{center}

The maps $F_{11}$ and $F_{21}$ correspond to avalanches of size
$0$, the maps $F_{12}$ and $F_{22}$ correspond to avalanches of
size $1$, and the maps $F_{13}$ and $F_{23}$ correspond to
avalanches of size $2$.

It was discovered in \cite{BCK} that the physical attractor $\Y$
has the following simple structure:
$$
\Y=\Big{\{} \big{(}\frac{1+\epsilon}{1-\epsilon},
\frac{\epsilon}{2-\epsilon} \big{)},
\big{(}\frac{\epsilon}{2-\epsilon},\frac{1+\epsilon}{1-\epsilon}
\big{)}, \big{(}\frac{1+\epsilon}{1-\epsilon},
\frac{1+\epsilon}{1-\epsilon} \big{)} \Big{\}}\,.
$$
for $E_c \in [\frac{1+\epsilon}{1-\epsilon},\frac{2}{1-\epsilon}]$
We denote these points by $a,b,c$ so that $\Y=\{a,b,c\}$. The maps
$F_1|_\Y$ and $F_2|_\Y$ are permutations of $\Y$:
$$
F_1|_\Y=\Big{(}\begin{matrix} a & b & c \\ b & c & a \end{matrix}
\Big{)}\,,\, F_2|_\Y=\Big{(}\begin{matrix} a & b & c \\ c & a & b
\end{matrix} \Big{)}\,.
$$
\noindent
\begin{figure}
\begin{center}
\includegraphics[width=10.5cm]{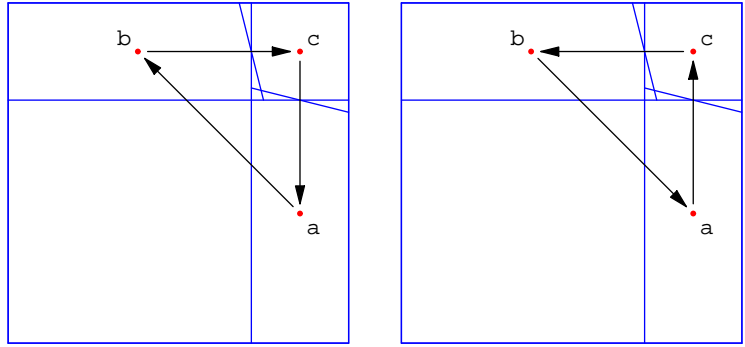}
\caption{{\small The figure shows the physical attractor
$\Y=\{a,b,c\}$ and the maps $F_1|_\Y, F_2|_\Y$ for $E_c=7/2$ and
$\epsilon=1/2$. The arrows on the picture to the left shows how
the points of $\Y$ are mapped under $F_1$, and the picture on the
right shows how the points are mapped under $F_2$.}} \label{figure 1}
\end{center}
\end{figure}

Figure \ref{figure 1} shows the physical attractor $\Y$ and the
maps $F_1|_Y$ and $F_2|_Y$ for $E_c=7/2$ and $\epsilon=1/2$. We
construct a partition ${\mathcal R}=\{R_1,\dots, R_6\}$ of
$\Sigma_N^+\times Y$ by
\begin{center}
\begin{tabular}{lll}
$R_1=[1]\times \{a\}$ & $R_2=[1]\times \{b\}$& $R_3=[1]\times \{c\}$ \\
$R_4=[2]\times \{a\}$ & $R_5=[2]\times \{b\}$& $R_6=[2]\times
\{c\}$
\end{tabular}\,.
\end{center}
We let $A=\|a_{ij}\|$ be the $6 \times 6$ matrix where $a_{ij}=1$
if $F(R_i) \cap R_j \neq \emptyset$ and $a_{ij}=0$ otherwise. It
is easy to verify that
$$
A=\begin{bmatrix}
0 & 1 & 0 & 0 & 1 & 0 \\
0 & 0 & 1 & 0 & 0 & 1 \\
1 & 0 & 0 & 1 & 0 & 0 \\
0 & 0 & 1 & 0 & 0 & 1 \\
1 & 0 & 0 & 1 & 0 & 0 \\
0 & 1 & 0 & 0 & 1 & 0 \\
\end{bmatrix}
$$
It is clear that $g_{\mathcal R}: \mathcal A \rightarrow
\Sigma_A^+$ is a topological conjugancy of the maps $F|_{\mathcal
A}$ and $\sigma_A^+$. The matrix $A$ is transitive and
$\text{Sp}(A)=\{-1,-1,0,0,0,2\}$. Hence $\h(F|_{\mathcal A})=\log
2$.

If $\mu$ is the SRB-measure on ${\mathcal A}$, with $(\pi_u)_*\mu$
being the uniform Bernoulli measure on $\Sigma_N^+$, then
$(g_{\mathcal R})_* \mu$ is the Perry measure on $\Sigma_A^+$.
With respect to this measure it is easy to see that the average
avalanche size is $\overline{s}_0=1$. It then follows from that
the sum of the negative Lyapunov-exponents is $\chi_1+\chi_2=\log
\epsilon$. For instance we see that for $E_c=11/2$ and
$\epsilon=2/3$ we have $h_{\mu}(F|_{\mathcal A}) >
|\chi_1|+|\chi_2|$. So even though $F_1|_\Y$ and $F_2|_\Y$ are
invertible, the Ruelle inequality (and therefore the Pesin
formula) cannot be reversed (though this was argued in
\cite{BCK}).

We also remark that for $E_c=(1+\epsilon)/(1-\epsilon)$ we have
$\Y \subset S(F)$, so the standard construction of SRB-measure
will fail in this case. Still there is clearly a natural invariant
measure.

The example with a trivial physical attractor can be generalized
to all $N$ for $d=1$. Then $\Y$ consists of $N+1$ points
$z_0,z_1,\dots,z_N$ given by
$$
z_n=\frac{1+\epsilon}{1-\epsilon}(1,1,\dots,1)-e_n
$$
where $e_0=0$ and $e_1,\dots,e_N$ is the standard basis in  $\R^N$.  \\

\noindent {\bf Example B:}
\noindent
\begin{figure}
\begin{center}
\includegraphics[width=10.5cm]{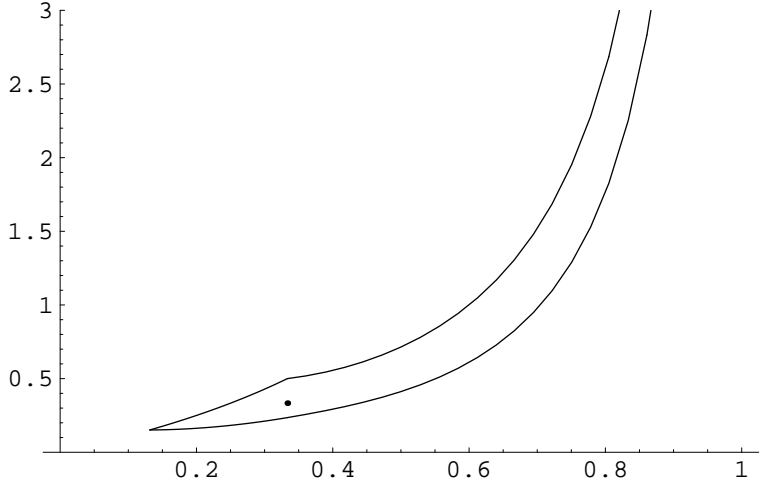}
\caption{{\small The region in the $(\e,E_c)$-plane where the
avalanches are the same as for $E_c=\e=1/3$. The point $(1/3,1/3)$
is shown in the interior of the region.}} \label{avalachedomain}
\end{center}
\end{figure}
For $E_c=\epsilon=1/3$ the dynamics is very simple. The map $F_1$
has two domains of continuity: $ M_{11}=\{(x,y)\in M\,|\,y \leq
-2x/3+1/3\}$ and $M_{12}=M \setminus M_{11}$. The domains of
continuity for $F_2$ are given by symmetry. The maps are given by
$F_{11}(x)=L_{11}(x+e_1)$ and $F_{12}(x)=L_{12}(x+e_1)$, where
$$
L_{11}=\frac19 \begin{bmatrix} 2 & 3 \\ 2 & 3
\end{bmatrix}\,\,\text{ and }\,\,
L_{12}=\frac{2}{27} \begin{bmatrix} 2 & 3 \\ 2 & 3
\end{bmatrix}\,.
$$
We see that all four $F_{ij}$ are mappings to the diagonal line in
$M$. In fact, $F_{11}(M_{11})=F_{21}(M_{21})=[p_1,p_2]$, where
$p_1=(2/9,2/9)$ and $p_2=(1/3,1/3)$. The interval is contained in
$M_{12} \cap M_{22}$ since the two lines of singularity intersect
the diagonal in the point $(1/5,1/5)$. The images of $M_{12}$ and
$M_{22}$ also coincide and is an interval $[p_1,p_3]\subset
[p_1,p_2]$, where $p_3=(22/81,22/81)$. This shows that
singularities are removable after one iteration. It is easy to see
that for all ${\bf t} \in \Sigma_N^+$ the dynamics will contract
to the fixed point $P=(4/17,4/17)$. Hence the attractor of the
system is ${\mathcal A}=\Sigma_N^+ \times \{P\}$. The dynamics on
the attractor is of course conjugated to $\sigma_N^+$.

The example can easily be extended to a neighborhood of
$(1/3,1/3)$ in the $(\e,E_c)$-plane. In the region
$$
\max \Big{\{}\frac{1-2\e+13\e^2}{5+12\e-13 \e^2},
\frac{\e(7\e^2-6\e+3)}{7\e^3-10 \e^2+7 \e-4} \Big{\}} \leq E_c
\leq \min \Big{\{}\frac{\e}{1-\e}, \frac{1-2\e+5\e^2}{1+4\e-5\e^2}
\Big{\}}
$$
we have the same avalanches as for $E_c=\e=1/3$. This region is
shown in Figure \ref{avalachedomain}. For all points in the region
the maps depend continuously on $\e$ and the lines of singularity
depend continuously on $\e$ and $E_c$. The condition for an atomic
spacial attractor is that the images of the domains do not
intersect the singularities. For $E_c=\e=1/3$ the images are
bounded away from the lines of singularity, and hence there exists
an open neighborhood of $(1/3,1/3)$ in the $(\e,E_c)$-plane where
the same holds, i.e. the attractor is of the form $\Sigma_N^+
\times \{P\}$ (where $P$ is a point in $M$) and the dynamics is
conjugated to $\sigma_N^+$. \\

\noindent {\bf Example C:} Let us consider $E_c=1/3$ and
$\epsilon=1/2$. In this case there are 28 domains of continuity
and 28 corresponding maps. A computer program is written to
compute the sets $U_n$. The program uses exact calculations of the
edges of the polygons that make up $U_n$, and hence it can be used
to give rigorous ''proof by computer'' of removability of
singularities. By using the program we obtain that singularities
are removable. In fact $U_5 \cap {\mathcal S}=\emptyset$. The set
$U_5$ consists of 13 connected components. Figure \ref{figure 2}
shows the set $U_5$ and the lines of singularity, and Figure
\ref{figure 3} is a schematic illustration of how these connected
components are situated with respect to the lines of singularity.
\noindent
\begin{figure}
\begin{center}
\includegraphics[width=5.5cm]{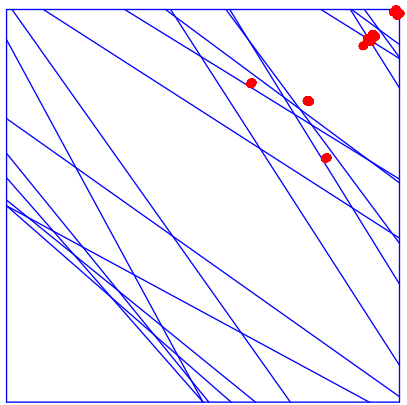}
\caption{{\small Shows the set $U_5$ for $N=2$, $E_c=1/3$ and
$\epsilon=1/2$. The points on the attractor are magnified in order
to make them visible in the figure, and hence it looks as if they
intersect singularities, but in fact they do not.}} \label{figure
2}
\end{center}
\end{figure}
\noindent
\begin{figure}
\begin{center}
\includegraphics[width=5.5cm]{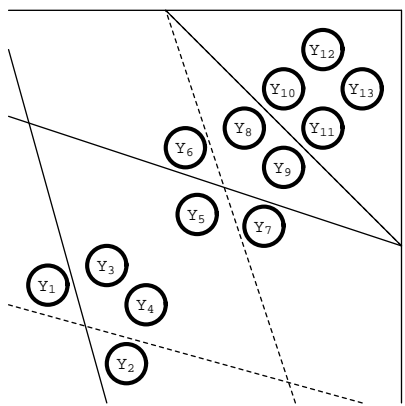}
\caption{{\small Shows how the 13 spatial partition elements
are situated with respect to the lines of singularity.}} \label{figure 3}
\end{center}
\end{figure}
\begin{figure}
\begin{center}
\includegraphics[width=5.5cm]{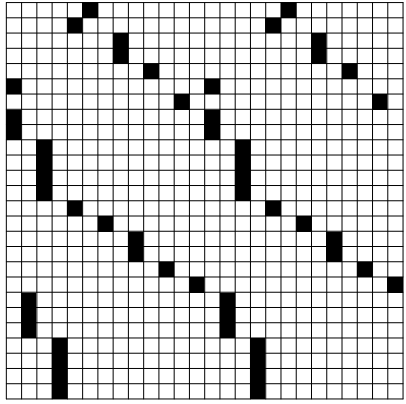}
\caption{{\small Shows the matrix $A$ for the coding of Example C.
Black squares are 1-s and white squares are 0-s.}} \label{figure
4}
\end{center}
\end{figure}
\noindent The intersection of the connected components of $U_5$
with $\Y$ are denoted by $Y_1,\dots ,Y_{13}$. Then we construct the
partition $\mathcal R=\{[i] \times Y_j\}$, and enumerate the
elements so that $\mathcal R=\{R_1,\dots,R_{26}\}$, where
\begin{center}
\begin{tabular}{llll}
$R_1=[1]\times Y_1$ & $R_3=[1]\times Y_2$  & \dots &
$R_{25}=[1]\times Y_{13}$
\\
$R_{2}=[2]\times Y_1$ & $R_4=[2]\times Y_2$  & \dots &
$R_{26}=[2]\times Y_{13}$
\end{tabular}\,.
\end{center}
We construct the $26 \times 26$ matrix $A=\|a_{ij}\|$ by letting
$a_{ij}=1$ if $F(R_i) \cap R_j \neq \emptyset$, and $a_{ij}=0$
otherwise. After making the computations, the matrix $A$ becomes
as shown in Figure \ref{figure 4}. The black squares represent
ones and white squares represent zeros. Direct computation shows
that the matrix $A$ is transitive.

Since singularities are removable the map $g_{\mathcal
R}:{\mathcal A} \rightarrow \Sigma_A^+$ is an avalanche conjugancy
between $F_{\mathcal A}$ and $\sigma_A^+$. The SRB-measure
projects to the Perry measure on $\Sigma_A^+$, so it is possible
to calculate properties such as average avalanche size. In this
example a computation gives $\overline{s}_0=123/17$.  The spectral
radius of the matrix $A$ is $2$, and hence
$\h(\sigma_A^+)=\log 2$.  \\

\begin{figure}
\begin{center}
\includegraphics[width=5.5cm]{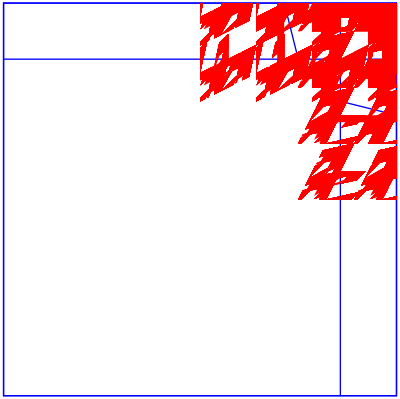}
\caption{{\small Shows the set $U_{20}$ for $E_c=7$ and $\e=1/2$. In this example
singularities are non-removable.}} \label{figure 5}
\end{center}
\end{figure}

\noindent {\bf Example D:} In the previous examples singularities
are removable. This is however not always the case. Figure
\ref{figure 5} shows the set $U_{20}$ for $N=2$, $E_c=7$ and
$\epsilon=1/2$. This is an example where singularities are
non-removable.
\begin{figure}
\begin{center}
\includegraphics[width=5.5cm]{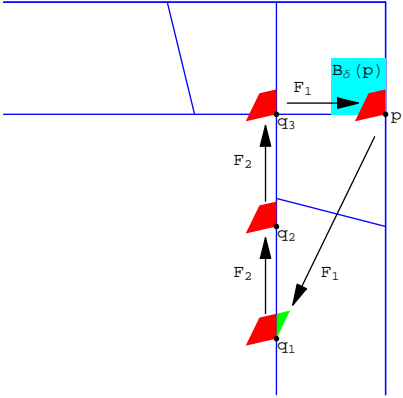}
\caption{{\small Illustration of the fact that singularities
are non-removable for $E_c=7$ and $\e=1/2$.}} \label{figure 6}
\end{center}
\end{figure}
\noindent We will in the following show that singularities are
non-removable for $E_c=(3+\epsilon)/(1-\epsilon)$.

Since  $(3+\epsilon)/(1-\epsilon)>(1+\epsilon)/(1-\epsilon)$, the
domains of continuity are given by the formulas presented in
Example A. Take the point $p=(E_c, E_c-1) \in M_{13}$. See Figure
\ref{figure 6}. Observe that $p$ lies on the horizontal line
$x_1=E_c-1$ and hence $p \in S(F)$. Clearly $p$ is in the interior
of $M_{13}$, so
$$
F_1(p)=F_{13}\Big{(}\frac{3+\epsilon}{1-\epsilon},
\frac{3+\epsilon}{1-\epsilon}-1\Big{)}=
\Big{(}\frac{3+\epsilon}{1-\epsilon}-1,
\frac{3+\epsilon}{1-\epsilon}-2 \Big{)}=p-(1,2)\,.
$$
Denote $q_1:=F_1(p)$ and observe the it lies on the singularity
line $x_2=E_c-1$. On the other hand $q_1$ is in the interior of
$M_{21}$, so $F_2(q_1)=p-(1,1)$. Denote $q_2:=F_2(q_1)$. This
point also lies in the interior of $M_{21}$. Let
$q_3:=F_{2}(q_2)=p-(1,0)$. It is clear that $F_{11}(q_3)=p$, and
in this sense $p$ is a periodic point. However, $F_{1}(q_3)$ is
not well defined since $q_3 \in \partial M_{12} \cap  \partial
M_{13}$.

In the following we let $\langle z_1,\dots,z_n \rangle$ denote the
open convex polygon with edges $z_1,\dots,z_n$. Define
$B_\delta(p)=\langle p, p+(0,\delta),p+(-\delta,
\delta),p+(-\delta,0)\rangle$. For small $\delta>0$ we have
$B_\delta(p) \in M_{13}$, and hence
$$
F_1(B_{\delta}(p))=F_{13}(B_{\delta}(p))=\langle q_1,q_1+\delta
a,q_1+\delta b,q_1+\delta c \rangle\,,
$$
where
$$
a=\Big{(} \frac{1-\epsilon}{2}, \epsilon \Big{)}\,,\,\, b=\Big{(}
\frac{1-4\epsilon-\epsilon^2}{4}, \frac{1+\epsilon}{2} \epsilon
\Big{)}\,,\,\, c=\Big{(}-\big{(} \frac{1+\epsilon}{2}\big{)}^2,
-\frac{1-\epsilon}{2} \epsilon \Big{)} \,.
$$
The polygon $F_1(B_{\delta}(p))$ intersects the singularity line
$x_1=E_c-1$. See Figure \ref{figure 6}. A simple computation shows
that
$$
F_1(B_{\delta}(p)) \cap M_{11}=\langle q_1,q_1+\delta
a',q_1+\delta b,q_1+\delta c \rangle\,,
$$
where
$$a'=\Big{(}0, \big{(} \frac{2\e}{1+\e} \big{)}^2 \Big{)}\,.$$
It then follows from the above discussion that
$$
F_1\circ F_2^2(F_1(B_{\delta}(p)) \cap M_{11})= \langle p,p+\delta
a',p+\delta b,p+\delta c \rangle\,.
$$
It is then easy to verify that for all $\delta>0$ there is
$\gamma>0$ such that $B_\gamma(p) \subset \pi_s\circ
F^4(\Sigma_N^+ \times B_\delta(p))$. So for each $n \in \N$ there
is $\delta_n>0$ with $B_{\delta_n}(p) \subset U_n$. The image of $
B_{\delta_n}(p)$ under $F_1$ intersects singularities, and its
closure contains the point $q_1$, which thus is an essential
singularity. Clearly the points $p$, $q_2$ and $q_3$ are also
essential singularities.


\section{\hps Statistical properties}\label{sec_7}

In order to evaluate the entropy and Lyapunov spectrum
of the physical model
in the thermodynamic limit we need to derive several
estimates for the asymptotic behavior of observables like
avalanche size, avalanche duration and "waiting-time" between
avalanches. The results are derived using only the uniform Bernoulli
measure on $\Sigma_N^+$, and hence hold for any SRB-measure and for
time-averages.

In \S \ref{sect13} we define the thermodynamic limit as the double
limit $E_c \rightarrow \infty$, $L=\sqrt[d]{N} \rightarrow
\infty$, contrary to \cite{BCK}, where the thermodynamic limit is
defined as the limit $L \rightarrow \infty$ only. This is
important as the quasi-classical limit since equivalently means a
fixed energy, but the energy quantum of Section \ref{sec_1}
$\delta\to0$. As $E_c \rightarrow \infty$ we must make a scaling
of time in the physical model. Otherwise the influx of energy to
the system will go to zero. With this new scaling we show that the
entropy goes to zero and the Lyapunov spectrum is collapsing.

We do not provide strict mathematical proofs, but still think
important to include the discussion of our results from the
physical point of view.


\subsection{\hpss Statistics of observables}\label{Ssec}

Let $\tau$ be the coordinate measuring the duration of avalanche
and let $\omega$ correspond to the interval between avalanches
(minimal value 1). We will also study the observable $s$ -- the
avalanche size (defined in \S \ref{return maps}).
While in the first case we
consider only actual avalanches, so that $\tau>0$, in the second
we make distinction between $s_0$ -- all avalanches including the
trivial case of under-critical state ($s=0$) and $s_+$ -- the
actual avalanches, so that $s_+>0$.

The reason for introducing two different avalanche size
observables is the following: $s_0$ plays a crucial role in
mathematical investigation of the model (see \S\ref{L.S.}), while
$s_+$ is important from physical perspective. In \S\ref{sect13} we
will see that physical observables should allow a thermodynamic
limit.

Denote by $\bar\tau$, $\bar\omega$, $\bar s_0$, $\bar s_+$ the
corresponding mean time-average quantities, each of which is a
function on the space-factor $M$ and is defined as follows:
 $$
\bar\sigma(x)=\int_{\Sigma_N^+}\lim_{k\to\infty}\frac1k
\sum_{i=0}^{k-1}\sigma\bigl(\hat f^i({\bf t},x)\bigr)\,d\mu_{\text{{\tiny Ber}}}
 $$
with $\mu_{\text{{\tiny Ber}}}$ being the Bernoulli measure on the
one-sided shifts (by Birkhoff ergodic theorem the time-average
limit exists almost everywhere and is measurable). This function
is invariant in the sense:
$\bar\sigma(x)=\frac1N\sum_{i=1}^N\bar\sigma(\hat f(i,x))$.

Whenever the system is ergodic with respect to an invariant measure $\mu$ the
function $\bar\sigma$ is constant $\mu$-a.e.
and equals the space (ensemble) average
$$
\langle\sigma\rangle_\mu=\int\limits_{\Sigma_N^+\times M}
\hspace{-8pt}\sigma\,d\mu.
$$
If the system has a unique invariant SRB-measure the function
$\bar\sigma$ is constant $\mu_{\text{\tiny Ber}} \times
\mu_{\text{\tiny Leb}}$-a.e. and equals the space (ensemble)
average
$$
\langle\sigma\rangle=\int\limits_{\Sigma_N^+\times M}
\hspace{-8pt}\sigma\,d\mu_{\text{\tiny SRB}}.
$$

We will need the maximal values of these observables in a sequel,
which we denote by $\tau_{\text{max}}$, $\omega_{\text{max}}$ and
$s_{\text{max}}=\max s_0=\max s_+$ respectively. We shall
calculate their asymptotics in $L$ and $E_c$.

Denote by $\varphi_1\sim\varphi_2$ the asymptotic equivalence
relation meaning that the ratio $\varphi_1/\varphi_2$ has
subexponential grow/decay. We denote the equivalence by $\approx$,
when the limit of the ratio is 1.

 \begin{lem}\po\label{lmm}
For $E_c\gg1$ the maximal avalanche time and size have the
asymptotics: $\tau_{\text{max}}\sim L^{\gamma_\tau}$,
$\omega_{\text{max}}=L^d E_c$, $s_{\text{max}}\sim
L^{d+\gamma_s}$, where $\gamma_\tau=\gamma_s=1$ for $d=1$ and
$1<\gamma_\tau,\gamma_s<d$ for $d>1$.
 \end{lem}
 \begin{proof}
Consider at first the simple case $d=1$.
The maximum avalanche duration and size are achieved when all
sites contribute to the avalanche, i.e.\ their energies are
sufficiently big and there is one site with energy greater than
$E_c-1$ to initiate the avalanche. Actually, if some site has
small energy, it will serve as a boundary and the avalanche wave
reflects from it (to be explained below).

Assuming $E_c>\frac\e{1-\e}$ we know from Lemma \ref{6} that in
the avalanche process each site, whenever overcritical, relaxes
until in the next step it receives a sufficient portion of energy
to become overcritical and relax etc. In other words, the sites
blink, being under- and over-critical in turn. But in this process
they make overcritical their neighbors and the process propagates
as a wave, with only difference that its front excites new sites,
while in the traversed region there remain blinking overcritical
sites.

This wave spreads along the interval $B_L^1=[1,L]\subset\Z$
towards its boundary and then it reflects from it, bearing now
relaxation. In fact, as the front wave reaches the boundary it
losses a substantial part of the energy on the boundary sites,
which thus cannot be recovered and remain undercritical for the
rest of this avalanche. They influence their neighbors to stop
being critical and so forth. Thus we obtain the reflected wave
that, in contrast with the first one, turns overcritical sites
into relaxed. When the wave hits itself, the avalanche process
stops.

It is clear that the duration of this avalanche is
$\tau_{\text{max}}\approx L$. The number of involved sites
corresponds to the area of the triangle with a side $B_L^1=[1,L]$
and height $L/2$ (recall that each site is overcritical only half
time of its blinking period), whence $s_{\text{max}}\approx
L^2/4$.

Consider now the case of dimension $d>1$. Here the scenario of
maximal avalanche is more complicated and consists of three stages
(the proof is similar to the case $d=1$, but quite lengthy and will be suppressed).
Again the maximum avalanche duration and size are achieved when
all sites contribute to the avalanche, though now if some isolated
site has smaller energy it serves as a boundary only once but then
on the next several waves it receives the required portion of energy
and follows the general scheme of motion.

At the first stage a site is excited and it initiates the
rhombus-shape wave (a cube in the Manhattan metric, see Figure
\ref{xxx}) that spreads to the boundary of the cube
$B_L^d=[1,L]^d\subset\Z^d$ (in time $\approx L/2$ if the center of
the rhombus is placed near the center of the cube).
\noindent
\begin{figure}
\begin{center}
\includegraphics[width=5.5cm]{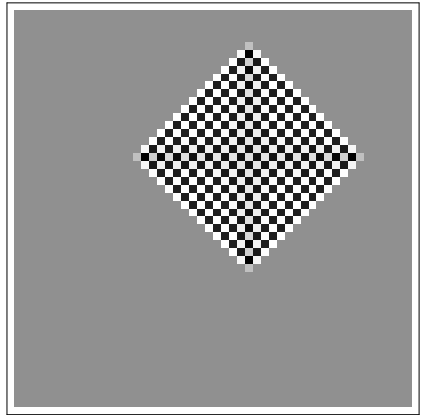}
\caption{{\small The picture is a "snapshot" of the lattice $\Lambda$ for
$d=2$, $L=10$, $E_c=7$ and $\epsilon=1/2$. A single site in the marginally
stable configuration has been exited, and a great avalanche is unfolding in a
rhombus shape. This is the first stage of this avalanche. The white squares are
overcritical sites, the gray squares are sites with energy just bellow $E_c$ and
the black squares have energy approximately equal $\epsilon E_c$.}} \label{xxx}
\end{center}
\end{figure}
The second stage begins as the wave reaches the boundary face and
reflects from it (See Figure \ref{zzz}). The reflected wave is
almost momentary overthrown by the coming overcritical wave, which
again reflects from the boundary, come now deeper into the
interior of the cube, but is overthrown too etc. If one looks
along the boundary face, the reflected wave travel along towards a
vertex with preserved form (like a soliton) and then disappears
into this vertex (there occur strong interactions with other waves
in this corner). But if one looks into the perpendicular
direction, the collection of reflected waves oscillates (each
reflected wave enters deeper and deeper into the cube)
contributing to the avalanche duration the sum $1+2+3+\dots$,
which stops with the end of the second stage (we do not specify
the sum precisely because after some oscillations the wave front
becomes more and more eroded by the interactions between
overcritical and relaxing waves; this impairs the sum and
decreases the exponent, but not too drastically).
\noindent
\begin{figure}
\begin{center}
\includegraphics[width=5.5cm]{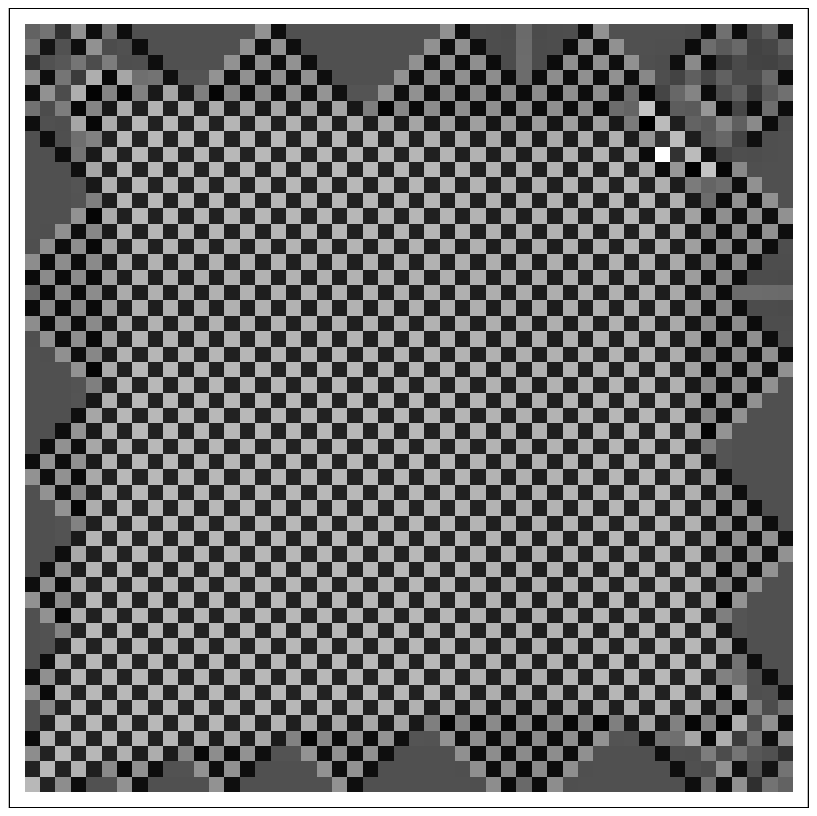}
\caption{{\small The picture shows a "snapshot" of the second
stage of the avalanche shown in Figure \ref{xxx}. We can see the
well-shaped (soliton-like) waves of energy near the boundary and
observe how their form begins being eroded near the vertices.}}
\label{zzz}
\end{center}
\end{figure}

The third stage begins as the main body of the overcritical sites
becomes disconnected and the avalanches behaves like worms
crawling along the high energy fractal-like collection of states.
We illustrate this in Figure \ref{xyz}.
\noindent
\begin{figure}
\begin{center}
\includegraphics[width=5.5cm]{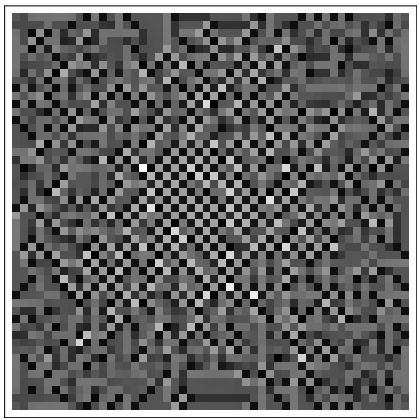}
\caption{{\small  The picture shows a "snapshot" of the third
stage of the avalanche shown in Figure \ref{xxx}. The energy
configuration has a fractal structure where the avalanche can
sustain in regions of high energy.}} \label{xyz}
\end{center}
\end{figure}
From the description of the maximal avalanche process it is clear
that the asymptotic exponents $\gamma_\tau,\gamma_s$ do not depend
on the energy $E_c\gg1$. Let us denote (the sequences are
increasing, so the limits exist):
$$
\gamma_\tau=\lim_{L\to\infty}\frac{\log\tau_{\text{max}}}{\log L}
\quad\text{ and }\quad d+\gamma_s=\lim_{L\to\infty}\frac{\log
s_{\text{max}}}{\log L}.
$$

Duration of the first stage of avalanche is $\sim L$. The above
arguments show that the exponent $\gamma_\tau'$ of the second
stage is $>1$. The last stage can only increase it:
$\gamma_\tau\ge\gamma_\tau'$. To see that $\gamma_\tau<d$ we note
that in average the number of critical sites on the boundary is
about $\kappa_L\sim L^{\gamma_\kappa}$ with $0<\gamma_\kappa<d$.
Thus we have a constant flow of energy out of the system with the
average speed $>\frac{1-\e}{2d}\kappa_L E_c$, and the inequality
$E_{\text{tot}}\le L^dE_c$ proves the claim.

To estimate the maximal avalanche size exponent $\gamma_s$
consider again the second stage of the above scenario. The number
of involved sites corresponds to the volume of the prism $\Pi_L$
over the cube $B_L^d$ with height
$\ell\approx\frac12\tau_{\text{max}}\sim L^{\gamma_\tau'}$, i.e.
 $$
L^{d+\gamma_s'}\sim\op{Vol}_{d+1}(\Pi_L)=\dfrac
\ell{d+1}\cdot\op{Vol}_d(B^d_L) \sim L^{d+\gamma_\tau'}.
 $$
Thus $\gamma_s\ge\gamma_s'=\gamma_\tau'>1$. Inequality
$\gamma_s<d$ follows from the inequality  from above for the
duration exponent $\gamma_\tau$.

The maximal value of the waiting time is obvious.
 \end{proof}


\subsection{\hpss Asymptotic of the statistical data}\label{S-sec}

Now we can study statistics of the avalanche data asymptotically
(as for the thermodynamic limit).

We will need the following technical statement (informal only for
$E_c>1$):
\begin{lem}\po\label{lml}
Almost every (w.r.t. a random excitation sequence) spatial
trajectory returns to the cube $B_0=(E_c-1,E_c]^N$.
\end{lem}
\begin{proof}
It follows from Proposition \ref{4} that $K=\{x\in M\,|\,\exists
i:x_i\in(E_c-1;E_c]\}$ is a return set, i.e. every trajectory
$F^n({\bf t},x)$ meets $\Sigma_N^+ \times K$. Partition the
spatial part $\Y$ of the attractor according to the hyperplanes
collection $\mathcal{H}=\cup_{i,m\in\N}\{x_i=E_c-m\}$.

Each such a part can be shifted by an excitation sequence to the
cube $B_0$. The probability of all such sequences (where the
avalanche starts from the set $B_0$ of maximal energy in all
sites) is positive. Let $\rho>0$ denote the minimum of these
probabilities over the finite set of all partition elements of $Y$
by $\mathcal{H}$. Then the probability of not entering $B_0$ in
$k$ successive avalanches is less than $(1-\rho)^k$.

Therefore since the number of avalanches tend to infinity as we
iterate the dynamics, the measure of trajectories staying away from
$B_0$ is zero.
\end{proof}

\begin{theorem}\po\label{tyh}
We have the following asymptotic estimates valid as $E_c\to\infty$
(and $N\gg1$ fixed) or $L=\sqrt[d]{N}\to\infty$ (and $E_c\gg1$
fixed):
\begin{enumerate}
\item[\rm(i)] $\bar\tau\sim L^{\gamma_\tau}$;
\item[\rm(ii)] $\bar\omega \sim E_c$;
\item[\rm(iii)] $\bar s_0\sim L^{d+\gamma_s}/E_c$;
\item[\rm(iv)] $\bar s_+\sim L^{d+\gamma_s}$,
\end{enumerate}
where $\gamma_\tau,\gamma_s$ are the same exponents as in Lemma
\ref{lmm} (thus $\gamma_\omega=\lim\limits_{L\to\infty}\frac
{\log\bar\omega}{\log L}=0$).
\end{theorem}
\begin{proof}
The maximal avalanche size is achieved for a certain configuration
of states $V\subset M$, which we can bound as follows: $U_1\subset
V\subset U_d$, where
 $$
U_j=\{x\in M\,|\,\forall i:
\bigl(1-j\tfrac{1-\e}{2d}\bigr)E_c<x_i\le E_c\text{ and }\exists
i_0: x_{i_0}>E_c-1\}.
 $$
Denote also $\tilde U_j=\{x\in M\,|\,\forall i:
\bigl(1-j\tfrac{1-\e}{2d}\bigr)E_c<x_i\le E_c\}\supset U_j$.

To estimate the measure of the sites leading to the maximal
avalanche we consider preimages of $\Sigma_N^+\times\tilde U_j$
under the map $F$. It is clear that one needs
$k_\e\in[\frac2{1-\e},\frac{2d}{1-\e}]$ different backwards
iterations $F^{-i_s}$, $s=1,\dots,k_\e$, to cover the spatial
attractor $\mathcal{Y}$. Since the measure $\mu$ is $F$-invariant,
we get for its $\pi_s$-push-forward: $\mu_s(\tilde
U_j)\approx\rho_j/k_\e$, which is $E_c$-independent. Thus we get
the same exponent $(d+\gamma_s)$ for $\bar s_+$ as for
$s_{\text{max}}$.

The same arguments yield the asymptotic of $\bar\tau$.

To obtain the asymptotic of $\bar\omega$ in $L$ we note that since
the amount of lost energy is $<C L^{\gamma_\kappa+\delta}E_c$
(where $\gamma_\kappa<d$ is the quantity from the proof of Lemma
\ref{lmm}), the average remained energy in a site of configuration
obtained from a maximal one after an avalanche is $(L^dE_c-C
L^{\gamma_\kappa+\delta}E_c)/L^d\approx E_c$. Thus the waiting
time does not grow with $L$.

The asymptotic of $\bar\omega$ in $E_c$ is quite different: If $N$
is fixed but $E_c$ grows, then any state from $\partial M$ becomes
at distance $\theta_N\cdot E_c$ after some relatively small number
of iterations.

Let us first demonstrate the idea in the simple case $N=2$. For
critical energy $E_c>\e/(1-\e)$ the picture of avalanches is shown
on Figure \ref{pic1}. We see that in a few steps of the dynamics
the configuration becomes far from $\partial M$, i.e. it strongly
contracts in all directions. Actually, it is possible to imagine
the situations when the point is mapped to the vertical strip and
then is shifted horizontally for a long time by excitations of the
first site, but probability of this event exponentially goes to
zero as $E_c\to\infty$.
\noindent
\begin{figure}
\begin{center}
\includegraphics[width=7.5cm]{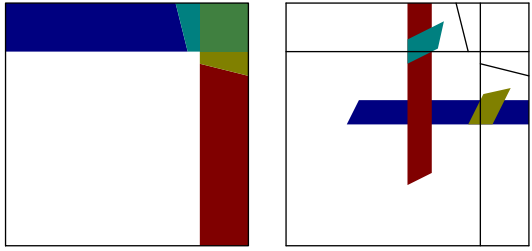}
\caption{{\small The figure illustrates how the continuity
domains are mapped under $F_1$ and $F_2$ for $N=2$ and
$E_c>\e/(1-\e)$.}}\label{pic1}
\end{center}
\end{figure}
Thus in a relatively short time the point from $\partial M$ is
mapped into the square $[0,\theta_2E_c]^2$, where the constant
$\theta_2$ is $E_c$-independent. To achieve the boundary $\partial
M$ again it needs $\sim(1-\theta_2)2E_c$ random excitations.

The similar picture happens for $N=3$, see Figure \ref{pic2}.
\noindent
\begin{figure}
\begin{center}
\includegraphics[width=7.5cm]{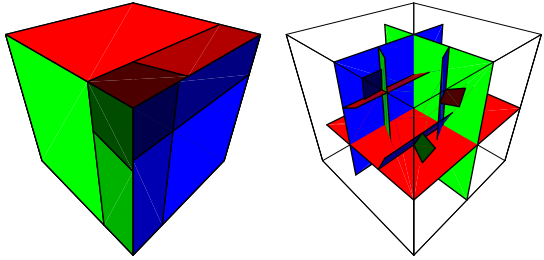}
\caption{{\small The figure illustrates how the faces
of the continuity domains are mapped under
$F_1$, $F_2$ and $F_3$ for $N=3$ and $E_c>\e/(1-\e)$. The
point of view is at $(E_c,E_c,E_c)$.}} \label{pic2}
\end{center}
\end{figure}
In the general case Theorem \ref{5} insures that after some (few)
number of steps we get strong contraction, so that the point
becomes in the cube $[0,\theta_NE_c]^N$, i.e. far from the
boundary $\partial M$. Thus we need $\sim(1-\theta_N)NE_c$
excitations to make it overcritical and this implies the claim.
Note that the asymptotic for $\bar\omega$ is not for the double
limit $E_c\to\infty$, $L\to\infty$, but for two partial limits
only.

To obtain the estimate for $\bar s_0$ we need to estimate the
conditional measure $\mu_s(U_j|\tilde U_j)$, which coincides with
the probability that a randomly chosen configuration $x\in \tilde
U_j$ and site $i$ satisfy: $x_i>E_c-1$. This probability is
$\approx b_1/E_c$. Thus $\mu(\Sigma_N^+\times V)\sim
\sigma_\e/E_c$ and the mean avalanche size is $\langle
s_0\rangle\sim L^{d+\gamma_s}/E_c$.

This is however the space-average (the arguments below work also
for $\bar s_+$, $\bar\tau$). We would obtain the same for the
time-average of Lebesgue a.e.-initial condition if we have an
SRB-measure, or for $\mu$-a.e. if we have an ergodic measure
$\mu$. But we cannot guarantee existence of an SRB-measure or an
ergodic measure. However, for a non-ergodic measure, we can
decompose $\mathcal{A}$ into ergodic components, where the
Birkhoff theorem works. By Lemma \ref{lml} each ergodic component
(with non-trivial contribution) intersects in its spatial part the
top-energy cube $B_0$ and so for each of it the same asymptotic of
space-averages holds with universal exponents but maybe different
coefficients. Therefore we obtain the required asymptotic for
$\bar s_0$.

Another way to get the last asymptotic is via the formula $\bar
s_0=\frac{\bar\omega\cdot 0+1\cdot\bar s_+}{\bar\omega+1}$.\!\!
\end{proof}
\vspace{0pt}

For $d=1$, $E_c\in[\frac{1+\e}{1-\e},\frac2{1-\e}]$ we already
have $\bar\tau\sim N=L^1$ as the theorem states, but $\bar s_0\sim
N=L^1$, $N\to\infty$, which shows in the last respect the critical
energy $E_c=2/(1-\e)$ is small.

In \cite{BCK} the estimate $\gamma_\tau>1$ was predicted for all
$d$, while this is a feature of the cases $d>1$. In the latter
cases the analytic calculation of exact values of exponents
$\gamma_\tau,\gamma_s$ is a difficult problem.

\begin{rk}\po
The difference between cases $d=1$ and $d>1$ demonstrated in the
theorem is known in the physical literature. The former case is
usually considered as the trivial SOC-model.
\end{rk}


\subsection{\hpss Thermodynamic limit}\label{sect13}

By thermodynamic limit of an observable $\phi$ we understand the
double limit
$$
[\phi]_\infty=\lim_{\begin{array}{c} \scriptstyle L\to\infty\vspace{-5pt}\\
\scriptstyle E_c\to\infty\end{array}}\phi
$$
if it exists. It is assumed that for physical observables this
limit exists. In \cite{BCK} only limit $L\to\infty$ was
considered, though then the value of energy $E_c$ could serve as
an essential parameter, which is not desirable in the
SOC-paradigm. However it was suggested there that consideration of
$E_c\to\infty$ can be helpful.

As an example of non-physical observable we expose $s_0$ (in
\S\ref{Ssec} we called it mathematically relevant): The double
limit does not exists because the repeated limits are different:
 $$
0=\lim_{L\to\infty}\lim_{E_c\to\infty}\bar s_0\ne
\lim_{E_c\to\infty}\lim_{L\to\infty}\bar s_0=+\infty.
 $$
But $s_+$ and $\tau$ are good physical observables, for the
thermodynamic limits exist:
 $$
[\bar\tau]_\infty=\infty,\quad [\bar s_+ ]_\infty=\infty.
 $$
From \S\ref{305} (use iteration arguments in the topological case)
and \S\ref{evalofentr} we obtain:
 $$
[\h(\hat f)]_\infty=[\h(F)/\bar\tau]_\infty=0 \quad\text{ and }
\quad[\h(F)]_\infty=\infty.
 $$

But with $\omega$ the situation is different because the proof
(rather than the vague statement of part (ii)) of Theorem
\ref{tyh} implies:
$\lim_{L\to\infty}\lim_{E_c\to\infty}\bar\omega=\infty$, while
$\lim_{E_c\to\infty}\lim_{L\to\infty}\bar\omega$ is finite. Since
$\omega$ is definitely physically relevant observable one needs
the following reparametrization: $\omega\mapsto\omega/E_c$, which
corresponds to contraction of the waiting time via the following
ansatz:

We let the energy quantum added at a unit time to the system equal
$\delta=\hbar$ (instead of 1 as before), but speed up time in the
waiting intervals respectively: $E_c=E_0/\hbar$,
$\omega_\text{new}=\omega\hbar$. Then the thermodynamic limit (the
space part $M$ can be quantized similarly via $L=[l/\hbar]$ with
$l$ a finite length) corresponds to the quasi-classical limit
$\hbar\to0$.

The duration of avalanche was suppressed in our definition of
dynamics to length one. So if we want to find the entropy of the
physical system, where each step of avalanche has time-duration
one, we should multiply it by the probability of dropping energy
into the system. For every trajectory this equals
$\frac{\bar\omega}{\bar\omega+\bar\tau}$. This ratio behaves
differently as $L\to\infty$ or $E_c\to\infty$, so we need the
reparametrization described above.

In \S \ref{evalofentr} we calculated the entropy of the "return"
Zhang model $\h(F)=\log N$ (this was proved almost surely, but
even with the possible entropy growth for some exceptional
parameters the thermodynamic limit below is unaltered). But after
reparametrization it changes. Denoting by $h^\text{Zhang}$ the
entropy of the reparametrized system we get:
 $$
h_\mu^\text{Zhang}=h_{\mu^u}(\sigma_N^+)\cdot
\langle\frac{\bar\omega_\text{new}}{\bar\omega_\text{new}
+\bar\tau}\rangle,\quad
 \h^\text{Zhang}=d\cdot\log L\cdot
\bigl(\frac{\bar\omega_\text{new}}{\bar\omega_\text{new}
+\bar\tau}\bigr)_{\text{max}}.
 $$
This implies:
 $$
[h_\mu^\text{Zhang}]_\infty=[\h^\text{Zhang}]_\infty=0.
 $$
Therefore the expanding property is lost in the thermodynamic
limit for the original physical system, as was already noticed in
\cite{BCK} for a bit different situation.

\begin{rk}\po
Notice that in the reparametrized system
$\bar\omega_\text{new}\ll\bar\tau$, which is counter-intuitive for
certain SOC-examples (sandpile, earthquakes etc, where one expects
$\bar\omega\gg\bar\tau$). This indicates that the Zhang system
should be modified by introducing the local contraction of time
depending on the avalanche size or speed. We will not consider
such gradient-type models here.
 \end{rk}


\subsection{\hpss Lyapunov spectrum}\label{L.S.}

In \S \ref{hypstr} we showed that the Zhang model is hyperbolic with one
positive exponent and the remainder of the spectrum negative. This hyperbolicity
is lost in the thermodynamic limit.
\begin{prop}\po
For $E_c \geq \e/(1-\e)$ we have
$\sum_{i=1}^N\chi_i^-=\bar s_0\log\e$.
\end{prop}
\begin{proof}
We cannot ensure the existence of an SRB-measure, so
 both the Lyapunov exponents
and $\overline{s}_0$ should be seen as functions on $M$. From the
general theory of Lyapunov exponents we know that
$$
\sum_{i=1}^N \chi_i^-(\hat{x}) =
\lim_{n \rightarrow \infty} \frac1n \sum_{t=1}^{n-1}
\log \det {\mathcal T}(F^t\hat{x})\,.
$$
From Proposition \ref{7} we know that the formula
$\det(L_{ij})=\e^{s_{ij}}$ holds for $E_c \geq \e/(1-\e)$. Hence
\begin{eqnarray*}
\sum_{i=1}^N \chi_i^-(x)&=&
\sum_{i=1}^N \int_{\Sigma_N^+}
\chi_i^-({\bf t},x)
\,d\mu_{\text{{\tiny Ber}}} \\
&=&\int_{\Sigma_N^+} \sum_{i=1}^N
\chi_i^-({\bf t},x)
\,d\mu_{\text{{\tiny Ber}}} \\
&=& \int_{\Sigma_N^+} \lim_{n \rightarrow \infty}
\frac{1}{n} \sum_{t=0}^{n-1} \log (\det {\mathcal T})(F^t({\bf t},x))\,
d\mu_{\text{{\tiny Ber}}} \\
&=& \int_{\Sigma_N^+} \lim_{n \rightarrow \infty}
\frac{1}{n} \sum_{t=0}^{n-1} s({\bf t},x)\,d\mu_{\text{{\tiny Ber}}} \log \e \\
&=& \overline{s}_0(x) \log \e\,.
\end{eqnarray*}
\end{proof}

 \begin{cor}\po
$|\overline{\chi\mathstrut}{^-}|=\frac{\bar s_0}N\log\frac1\e\to
0$ as $E_c\to\infty$, but $|\overline{\chi\mathstrut}{^-}|\to
\infty$ as $L\to\infty$.
\end{cor}
~ \\
\noindent
Thus we should study differently the following cases:
~\\

\noindent
{\bf 1. $E_c\to\infty$, but $L$ (and $d$) fixed.} Since
$\lim_{E_c\to\infty}\bar s_0=0$, the negative part of the Lyapunov
spectrum collapses: $\lim_{E_c\to\infty}\chi^-_i=0$ for every
$1\le i\le N$. The hyperbolicity is lost, but the positive
exponent $\chi^+_0=\log N$ survives. In particular, the entropy
does not collapses. \\

\noindent
{\bf 2. $L\to\infty$, but $E_c\gg1$ fixed.} Here only a bounded
piece $\chi_1^-,\dots,\chi_k^-$ of the Lyapunov spectrum
collapses, $k=\op{const}$. But the number of elements of this
spectrum grows and in average $|\overline{\chi\mathstrut}{^-}|\to
\infty$. In particular, $|\chi^-|_\text{max}\to \infty$ as
$L\to\infty$. Again, the positive exponent $\chi_0$ and entropy
are preserved and though we loose hyperbolicity there are many
non-degenerate Oscelledec modes. Moreover they prevail over
collapsing modes and so essentially the hyperbolicity is preserved
as well. \\

\noindent
{\bf 3. Reparametrized model.}
This was introduced in \S\ref{sect13} and require the
renormalization: multiplication of waiting time by the function
$\frac{\bar\omega}{\bar\omega+\bar\tau}$ along the trajectory. In
this case the Lyapunov spectrum collapses to zero in any limit
$E_c\to\infty$ and $L\to\infty$ and the hyperbolicity is
completely lost. \\

Thus the exponential grow of the statistics is suppressed and we
can observe power law statistic as is the basic idea of
SOC-phenomenon. The corresponding SOC-exponents are related to the
asymptotic of the Lyapunov spectrum (as discussed in \cite{BCK}),
but are difficult to calculate analytically.


\section{\hps Conclusion}

It is of importance to the paradigm of SOC to understand the
mechanisms behind the behavior one observes numerically in the
Zhang model and other sandpile models, and the aim of this paper
is to provide a first step to a rigorous mathematical
understanding of the dynamics of Zhang model.

Due to the singularities and non-invertibility of the model there
existed very few applicable results, and hence we had to modify
known results and develop some new methods in order to describe
the dynamical properties of the model.

The result of this work is that the singularities play a modest
role in the sense they do not change the main dynamical
characteristics. However the effects of the singularities can be
seen in the rich fractal structure of the spacial attractor.

Our analysis allows to take the thermodynamic limit of the main
dynamical quantities, showing that the entropy vanishes and the
Lyapunov-spectrum collapses after re-scaling the model. From a
physical point of view this is interesting. Typically chaotic
dynamics (positive entropy and positive Lyapunov exponents) is an
indication of exponential speeds of mixing (short decay of
correlations) \cite{Ba}, which is not compatible with the
power-law statistics of the SOC-hypothesis. The loss of
hyperbolicity in the thermodynamic limit hence supports the
critical behavior observed numerically \cite{J}, \cite{Z},
\cite{GD}.

To conclude: We have shown that the Zhang model is a chaotic
hyperbolic dynamical system, where all the entropy is produced by
the random driving of the system and, due to singularities, the
orbit structure is richer than for a topological Markov chain. The
hyperbolicity indicates that under weak conditions there exists an
SRB-measure (self-organization). In the thermodynamic limit the
hyperbolicity is lost and we may expect power-law statistics
(criticality). In practice the systems that are studied have
finite size and finite critical energy. Hence they are chaotic,
but with small entropy. The SOC-hypothesis is that these weakly
chaotic systems have SRB-measures with the rates of convergence to
these measures being exponential but slow compared to a unit step
in an avalanche, causing the prevailing of power-law statistics.

Finally note that modifications of the Zhang model are possible.
For instance, different amounts of energy $\delta_i$ can be added
to different sites $i\in\Lambda$ in the excitation process. This
corresponds to the rectangular form of $\Lambda\subset\Z^d$ and
uniform quantum of energy $\delta=\hbar$. We can also consider
other spacial configurations $\Lambda$. In this way nearly Zhang
models considered in Section \ref{evalofentr} are natural. Another
approach is to use a random amount of energy, i.e. stochastic
$\delta:\Lambda\to\R_+$ (this idea was used in the plasma physics
\cite{KK}), which can also be presented as a skew product with
piecewise affine fibers, but now over a solenoidal system. In all
these theories the ideas from the present paper work well (though
the thermodynamic limit will be sensitive to the form of
$\Lambda$).


\appendix


\section{Topological entropy of piecewise affine maps}\label{app_A}

\subsection{\hpss The Buzzi theorem and its generalization}\label{app_A.1}
The statement as it is done in \cite{B1} does not apply to our
situation. In \cite{B2} Buzzi noted that it extends to isometries
and contractions. In fact, the assertion holds always, but since
we cannot make a simple reference we write an adapted proof for
the convenience of the reader.

 \begin{theorem} \po \label{buzzi}
Let $X\subset \R^d$ be a bounded polytope and $f: X \rightarrow X$
a piecewise affine map. Then
 $$
\hs(f) \leq \lambda^+(f) + H_{\text{{\em mult}}}(f)\,,
 $$
where
 $$
\lambda^+(f)=\mathop{\overline\lim}\limits_{n \rightarrow \infty}
\sup_{x \in X} \frac{1}{n}\max_k \log \|\Lambda^k d_x f^n\|\,.
 $$
 \end{theorem}
\begin{proof}
Let ${\mathcal P}=\{P_i\}$ be the continuity partition and
$f_i:=f|_{P_i}$ be affine maps. Fix $\e>0$ and
$T=T(\e)\ge\frac{d}\e\log(\sqrt{d}+1)$ such that for $n\ge T$ we
have:
 $$
\text{mult}({\mathcal P}^n)\leq \exp((H_{\text{mult}}(f)+\e)n),\
\|\Lambda df^n\| \leq \exp((\lambda^+(f)+\e)n).
 $$
Take $r=r(\e)$ to be compatible with the partition ${\mathcal
P}^T$ (for $f^T$), i.e. any $r$-ball intersects maximally
$\op{mult}(\mathcal{P}^T)$ partition elements.

We will prove that each non-empty cylinder $C({\bf
a})=[P_{a_0}\dots P_{a_{lT-1}}]$ of length $|{\bf a|}=lT$ can be
partitioned into a collection $Q({\bf a})=\{W\}$ satisfying the
following properties:
\begin{enumerate}
\item $\sum_{|{\bf a}|=lT} \text{card}(Q({\bf a}))
\leq C_0 \exp((\lambda^+(f)+H_{\text{mult}}(f)+3\e)lT)$
\item $\text{diam}(f_{a_{lT-1}}\circ \dots \circ f_{a_0}(W)) \leq r\,.$
\end{enumerate}

Let us prove this claim by induction assuming it holds for some
$l\ge0$. The base of induction is obvious and $C_0$ is the minimal
cardinality of an $r/2$-ball cover.

Take a partition element $W \in Q({\bf a})$ that is used to cover
the cylinder $[P_{a_0}\dots P_{a_{lT-1}}]$. By the induction
hypothesis it has diameter less than $r$, so it can be continued
to cover a non-empty cylinder of length $(l+1)T$ in at most
$\text{mult}({\mathcal P}^T)$ ways. So to cover the cylinders
$[P_{a_0}\dots P_{a_{lT-1}}P_{b_0} \dots P_{b_{T-1}}]$ we make a
division of $W$:
 $$
W=\bigcup_{i=1}^\gamma W'_i\,\,,\gamma \leq \text{mult}({\mathcal
P}^T)\,.
 $$
Let
 $$
W''_i=f_{a_{Tl-1}}\circ\dots \circ f_{a_0}(W'_i)
 $$
and
 $$
W'''_i=f_{b_{T-1}}\circ\dots \circ f_{b_0}(W''_i)\,.
 $$
By the assumption $\text{diam}(W''_i)<r$ for all
$i=1,\dots,\gamma$, but the sets $W'''_i$ may have greater
diameter. We need to divide the sets $W'''_i$ so that they have
diameter less than $r$, and then pull this refinement back to the
partition of sets $W'_i$.

Let $L$ be the differential of $f_{b_{T-1}} \circ \dots \circ
f_{b_0}$ on $W_i''$. We can assume that $L$ is symmetric and take
$\{e_k\}$ to be a basis of eigenvectors corresponding to
eigenvalues $\lambda_1,\dots,\lambda_d$. Let $\{v_k\}$ be a basis
in the vector subspace corresponding to $W'''_i$. We can choose
this basis to be orthonormal and triangular with respect to
$\{e_k\}$. Divide $W'''_i$ by the hyperplanes
 $$
\psi_j(x)\stackrel{\text{def}}=\langle v_j,x\rangle=p
\frac{r}{\sqrt{d}},\quad p\in \Z,\quad j=1,\dots,d\,.
 $$
This defines cells $\tilde W$ of diameter less than $r$. Since
$\psi_j(W'''_i)=\psi_j(L(W''_i))$ has $\op{diam}\le|\lambda_i|r$,
the number of cells $\tilde W$ needed to cover $W'''_i$ is less
than or equal to
 $$
(\sqrt{d}+1)^d|\lambda_1|^+ \dots |\lambda_d|^+ \leq
(\sqrt{d}+1)^d \|\Lambda d f^n\| \leq \exp \big{(}(\lambda^+(f) +
2 \e)T \big{)}\,,
 $$
where $|\lambda_i|^+=\max\{|\l_i|,1\}$.

Therefore the total cardinality of the new partition is less than
or equal to
 \begin{eqnarray*}
\op{mult}({\mathcal P}^T)\exp \big{(}(\lambda^+(f)+2\e)T\big{)}
\exp \big{(}(\lambda^+(f)+H_{\text{mult}}(f)+3\e)lT \big{)}  \\
\leq \exp \big{(}(\lambda^+(f)+H_{\op{mult}}(f)+3\e)(l+1)T
\big{)}\,.
 \end{eqnarray*}
This proves the statement.
\end{proof} \\

The theorem holds as well for most degenerate piece-wise affine
systems, but there can be problems with $\hm(f)$. Namely the
latter is not defined if the image of continuity domain contains a
boundary face of a continuity domain. But if we assume the
image and the faces always meet transversally, no problems occur
and the above theorem applies literally.

In degenerate cases of the Zhang model, the above requirement
holds for most parameters. For instance, if $N=2$ and $\e=0$, then
all $E_c\notin\Z$ satisfy the request. For integer $E_c$ the
theory fails, but one can look just to the whole image set, which
has dimension $<N$: its Poincar\'e return map is piece-wise affine
and it satisfies the requirements. In the above example $N=2$,
$\e=0$ the dynamics is confined to two one dimensional lines,
whence the multiplicity entropy is zero (by dimensional reasons)
and $\h(F)=\log 2$ in this case.


\subsection{\hpss Entropy of conformal piecewise affine skew-products}\label{app_A.2}

We say that an affine map is conformal modulo degenerations if the
image on some subspace transversal to the kernel is mapped
conformally to its image. A piecewise affine map is said to be
conformal modulo degenerations if all its affine components are
conformal modulo degenerations. We will assume that degenerations
satisfy the transversality requirement of \ref{app_A.1}.

In \cite{KR1} we noticed that $\hm=0$ for piece-wise affine
conformal maps. This easily extends to allow degenerations. Now we
consider a more general situation of skew-product systems of
Zhang's type.

 \begin{theorem}\po \label{14}
Let $f_i:X\to X$ be piece-wise affine non-strictly contracting and
conformal modulo degenerations, $i=0,\dots,N-1$. Define
$F:\Sigma_N^+ \times X\to\Sigma_N^+ \times X$ by the formula
$F({\bf t},x)=(\sigma_N^+({\bf t}),f_{t_0}(x))$, ${\bf
t}=t_0t_1\dots\in\Sigma_N^+$, $x\in X\subset \R^d$. Then we have:
$\h(F)=\log N$ (so that a-posteriori the variational principle
holds).
 \end{theorem}

 \begin{rk}\po
For $N=2$ and $\e=0$ the affine components have rank 1, and hence
Theorem \ref{14} shows $\h(F)=\log N$. The same holds for
$E_c=\e=1/3$.
 \end{rk}

\noindent
 \begin{proof}
It follows from Theorem \ref{buzzi} that it suffices to prove that
$\hm(F)=0$.

To achieve the desired equality note that preimages of the
time-like singularity planes $\{t=n/N\}$ never intersect under
inverse iterations of $F(t,x)=(Nt \mod 1,f_{[Nt]}(x))$. Thus
$\hm(F)=\sup_{\bf t}\hm(f_{\bf t})$. We will prove that
$\hm(f_{\bf t})=0$ for all ${\bf t}\in\Sigma_N^+$.

A piecewise affine map can be considered as an ordered triple
$(X,{\mathcal P},f)$, where $X$ is a polytope in $\R^d$,
${\mathcal P}=\{P_i\}$ is a partition of $X$ made up of pairwise
disjoined polytopes (with certain faces of boundary included, so
that the whole boundary is distributed between polytopes) and
$f_i:=f|_{P_i}:P_i \rightarrow X$ are affine maps. We let
$X'=\cup\op{Int}(P_i)$ and $\op{Sing}(f)=X\setminus X'$.

For a piecewise affine map $(X,{\mathcal P},f)$ and a point $x \in
X$ construct a piecewise affine map $(X_x,{\mathcal P}_x, f_x)$,
called the differential of $f$ at $x$, by letting
\begin{enumerate}
\item $X_x=\{ y \in \R^d \,|\, \exists \e_0>0 \,\text{s.t.}
\forall \e \in (0,\e_0):\,x+\e y \in X\}\subset\R^d$ is the
tangent cone to $X$.

\item ${\mathcal P}_x$ is the partition of $X_x$
consisting of non-empty sets
 $$
P_x=\{ y \in \R^d \,|\, \exists \e_0>0 \,\text{s.t.} \forall \e
\in (0,\e_0):\,x+\e y \in P\}\,,
 $$
where $P \in {\mathcal P}$.

\item $f_x:P_x \rightarrow X_{f(x)}$ is the collection of maps
 $$
f_x(y)= \lim_{\e \rightarrow 0^+} \frac{f(x+\e y)- \lim_{\delta
\rightarrow 0^+} f(x+\delta y)}{\e}\,,\quad P_x\in{\mathcal
P}_x\,.
 $$
\end{enumerate}

Consider iterated differentials and denote
$f_{x_1,\dots,x_n}:=(\dots (f_{x_1})_{x_2}\dots)_{x_n}$. Note that
$(f^n_{\bf t})_x=(f_{t_{n-1}})_{f_{\bf t}^n(x)}\circ
\dots(f_{t_1})_{f^1_{\bf t}(x)}\circ(f_{t_0})_x$.


 \begin{prop}\po \label{15}
There exists a constant $C\in\R_+$ such that for any subspace
$W\subset\R^d$ and any $(x_1,\dots,x_r)\in X\times W^{r-1}$ we
have:
$$
\op{mult}((f^n_{\bf t}|_W)_{x_1,\dots,x_r})\le
C\,\Bigl(\sup_{V\subset\R^d}\sup_{y_1,\dots,y_{r+1}}\op{mult}
(f_{\bf t}|_V)_{y_1,\dots,y_{r+1}}\Bigr)^n,\qquad \forall n\ge0,
$$
where the collection of points $(y_1,\dots,y_{r+1})$ runs over
$X\times V^r$ with the condition
$\op{rank}(y_2,\dots,y_{r+1})+\op{codim}V
=\op{rank}(x_2,\dots,x_r)+\op{codim}W+1$.
 \end{prop}

The theorem follows from this, because for
 $$
\mu(r)=\mathop{\overline\lim}\limits_{n\to\infty}\frac1n
\log\hspace{-25pt} \sup_{\begin{array}{c}
\scriptstyle y_1\in X,y_2\dots,y_r\in V\subset\R^d\\
\scriptstyle \op{rank}(y_2,\dots,y_r)+\op{codim}V=r-1
 \end{array}}\hspace{-25pt}
 \op{mult}(f_{\bf t}^n|_V)_{y_1,\dots,y_r}
 $$
we have: $\hs(f_{\bf t})=\mu(1)\le\mu(2)\le\dots\le\mu(d+1)=0$.

To prove the proposition note that when $r=0$ we have (in this
case we do not need $W,V$):
$$
\op{mult}(f_{\bf t}^n)\le\Bigl(\sup_{x\in X} \op{mult}(f_{\bf
t})_x\Bigr)^n=\Bigl(\max_{0\le j<N}\sup_{x\in X}
\op{mult}(f_j)_x\Bigr)^n.
$$
Let $r\ge1$. Consider the continuity partition $\mathcal{P}^{\bf
t}_{x_1,\dots,x_r}$ for $(f_{\bf t})_{x_1,\dots,x_r}$, which is
just the continuity partition
$\mathcal{P}^{(t_0)}_{x_1,\dots,x_r}$ of
$(f_{t_0})_{x_1,\dots,x_r}$ (one iteration), and let
$\mathcal{P}^{n,\bf t}_{x_1,\dots,x_r}$ for $(f^n_{\bf
t})_{x_1,\dots,x_r}$ be the iterated partition. Note that the
latter is the collection of all non-empty intersections
$$
P^{t_0}\cap (f_{t_0})_{P^{t_0}}^{-1}(P^{t_1})\cap
(f_{t_0})_{P^{t_0}}^{-1}(f_{t_1})_{P^{t_1}}^{-1}(P^{t_2})\cap\dots
\cap(f_{t_0})_{P^{t_0}}^{-1}\dots(f_{t_{n-1}})_{P^{t_{n-1}}}^{-1}(P^{t_n}),
$$
where $P^{t_i}$ are elements of
$\mathcal{P}^{(t_i)}_{x_1,\dots,x_r}$ and $f_P$ denotes the
restriction of the (differential of the) map to the corresponding
continuity domain. Every element of these partitions is invariant
under the shift by vectors from $\op{span}(x_2,\dots,x_r)$.
Therefore it intersects the unit sphere $S_1(x_2,\dots,x_r)^\perp$
in the orthogonal complement. Consider the induced partition on
the sphere and refine it so that every element has diameter no
greater than $\varepsilon$. Denote by $n(\mathcal{P}^{n,\bf
t}_{x_1,\dots,x_r}, \varepsilon)$ the minimal cardinality of such
a refinement. Let also $m(\mathcal{P}^{(t_i)}_{x_1,\dots,x_r},
\varepsilon)$ be the maximal number of elements of
$\mathcal{P}^{(t_i)}_{x_1,\dots,x_r}$ that an $\varepsilon$-ball
$B(y,\varepsilon)\cap S_1^\perp$ of $S_1(x_2,\dots,x_r)^\perp$ can
meet.

Denote by $n(\mathcal{P}\cap W,\varepsilon)$, $m(\mathcal{P}\cap
W,\varepsilon)$ the corresponding quantities in the subspace $W$.
Then from the above formula for the iterated partition:
 $$
n(\mathcal{P}^{n+1,\bf t}_{x_1,\dots,x_r}\cap W, \varepsilon)\le
n(\mathcal{P}^{n,\bf t}_{x_1,\dots,x_r}\cap W, \varepsilon)\cdot
m(\mathcal{P}^{(t_{n+1})}_{y_1,\dots,y_r}\cap V, \varepsilon),
 $$
where $y_1=f_{\bf t}^n(x_1)$, $y_2=f_{\bf t}^n{}'(x_2),\dots$,
$y_r=f_{\bf t}^n{}'(x_r)$ and $V=f_{\bf t}^n{}'(W)$ with $f_{\bf
t}^n{}'=(f_{\bf t}^n)_{x_1,\dots, x_r}$. Therefore
 $$
n(\mathcal{P}^{n+1,\bf t}_{x_1,\dots,x_r}\cap W, \varepsilon)\le
n(\mathcal{P}^{n,\bf t}_{x_1,\dots,x_r}\cap W, \varepsilon)\cdot
\sup_V\sup_{y_1,\dots,y_r}m(\mathcal{P}^{(t_{n+1})}_{y_1,\dots,y_r}\cap
V,\varepsilon),
 $$
where the supremum is taken over all $V\subset\R^d$ and
$(y_1,\dots,y_r)\in X\times V^{r-1}$ such that codimension of
$\langle y_2\dots,y_r\rangle$ in $V$ equals codimension of
$\langle x_2\dots,x_r\rangle$ in $W$.

Since for a fixed $\varepsilon$ the number
$n(\mathcal{P}^{t_0}_{x_1,\dots,x_r}, \varepsilon)$ is finite and
$$
m(\mathcal{P}^{(t_i)}_{y_1,\dots,y_r}\cap V,\varepsilon)\le
\sup_{y\in S_1(y_1,\dots,y_r)^\perp\cap V}
|\mathcal{P}^{(t_i)}_{y_1,\dots,y_r}\cap B(y,\varepsilon)\cap V|,
$$
the claim follows from the following statement. Fix $i\in[0,N)$.

 \begin{lem}\po\label{16}
There exists $\varepsilon>0$ (depending only on $i$) such that for
all $V\subset\R^d$ and all $(y_1,\dots,y_r)\in X\times V^{r-1}$,
with $\op{rank}(y_2,\dots,y_r)<\dim V$, and $y\in
S_1(y_2,\dots,y_r)^\perp\cap V$ there exists
$y'\in(y_2,\dots,y_r)^\perp\cap V$ satisfying:
 $$
|\mathcal{P}^{(i)}_{y_1,\dots,y_r}\cap B(y,\varepsilon)\cap V|\le
\op{mult}((f_{i}|_V)_{y_1,\dots,y_r,y'}).
 $$
 \end{lem}
This statement, modulo our notations and restrictions to $V$, is
proved in \cite{B2}. The proposition and hence the theorem follow.
\end{proof}


\subsection{\hpss Estimates on entropy by angular expansion rates}\label{app_A.3}

It is possible to estimate the effect of angular expansion on the
topological entropy of a piecewise affine map $f:X \rightarrow X$
by its spherizations. Define the piecewise smooth map
$d_x^{(s)}\!f:S T_x X \to S T_{f(x)}X$ given at $x\in X'$ by the
formula
 $$
d_x^{(s)}\!f(v)=\frac{d_x f(v)}{\|d_x f(v)\|}\,.
 $$
For $x\in\op{Sing}(f)$ and $v\not\in T_x\op{Sing}(f)$ (the tangent
cone) we let $d_x^{(s)}\!f(v)=\lim\limits_{\e\to+0}d_{x+\e
v}^{(s)} f(v)$. For other $(x,v)\in STX$ the map is not defined.
The angular expansion of $f$ is exactly the expansion in the
fibers of its spherization.

If $d_x f$ is degenerate we restrict to the orthogonal component
of its kernel, and consider the map
 $$
d_x f|_{\op{Ker}(d_x f)^\bot}:\op{Ker}(d_x f)^\bot \rightarrow
\op{Im}(d_x f)\,.
 $$
Then the map $S_x(f)=d_x^{(s)} f|_{\op{Ker}(d_x f)^\bot}:S
\op{Ker}(d_x f)^\bot \rightarrow S \op{Im}(d_x f)$ between
$(\op{rank}(d_x f)-1)$-dimensional spheres is given by the formula
 $$
v \mapsto \frac{d_x f|_{\op{Ker}(d_x f)^\bot}(v)}{\|d_x
f|_{\op{Ker}(d_x f)^\bot}(v)\|}.
 $$

For $i<d$ we define
 $$
\rho_i(f)=\ls \frac1n \sup_{(x,v)} \max_{0 \leq k \leq i} \log \|
\Lambda^k d_vS_x(f^n) \|.
 $$
Let $m_*=\min_x \op{dim} \op{Ker} (d_x f)$ and $d_*=d-m_*=\max_x
\op{rank} (d_x f)$, where $d$ is the dimension of $X$. The numbers
$\rho_i(f)$ can be non-zero only for $i<d_*$.

We have: $\rho_0(f)=0$. The number $\rho_{1}(f)$ measures the
maximal exponential rate with which angles can increase under the
map $f$. The numbers $\rho_i(f)$ for $i<d$ measure the maximal
rate of expansion of the restrictions to $i$-dimensional spheres.
If $f$ is conformal, then $\rho_i(f)=0$ for all $i$.

\begin{theorem}\!\!{\bf (\cite{KR2}).} \label{th2}
For piece-wise affine maps $\hm(f)\le\sum_{i=1}^{d_*-1}
\rho_i(f)$.
\end{theorem}

We define the maximal expansion rate
$$
\lambda_{\op{max}}(f)=\ls \sup_{x} \frac1n \log \|d_x f^n\|\,,
$$
and the minimal finite expansion rate
$$
\lambda_{\op{min}}(f)=-\ls \sup_{x} \frac{1}{n} \log \| (d_x
f^n|_{\op{Ker}(d_x f^n)^\bot})^{-1}\|\,.
$$
In \cite{KR2} we show that
$$
\rho_i(f)\leq i \Big{(}\lambda_{\op{max}}(f)-\lambda_{\op{min}}(f)
\Big{)}\,.
$$
This gives the following result (the same bound holds for
$\h(f)$):
\begin{theorem}\po \label{thnew}
For a piecewise affine map $f$ it holds:
$$
\hs(f) \leq \lambda^+(f)+\frac{d_* (d_*-1)}{2}
\Big{(}\lambda_{\op{max}}(f)-\lambda_{\op{min}}(f) \Big{)}\,.
$$
\end{theorem}


\section{\hps Generalization of the Moran formula}\label{app_B}

Consider an IFS $(M^d,f_1,\dots,f_N)$ on a Riemannian manifold
$M$, where the maps $f_i$ can possess singularities, but we assume
that they are mild in a sense that the number of continuity
domains is finite, any of them has piece-wise smooth boundary and
the map, restricted to any of the domains, smoothly extends to the
adjacent singularities (this is the case of the Zhang model).

Remark that the IFS can be interpreted as the dynamical system
$(\Sigma_N^+\times M,F)$, $\hat f({\bf t},x)=(\sigma_N^+{\bf
t},f_{t_0}x)$. Attractor $\mathcal{Y}$ of the IFS can be defined
via the attractor of the extended system $F$, which has the form
$\mathcal{A}=\Sigma_N^+\times\mathcal{Y}$.

Let $\|\cdot\|$ be the norm on $TM$ generated by the metric on
$M$. Denote
$$
s_i^+=\max_{x\in M}\|d_xf_i\|,\qquad s_i^-=\bigl(\max_{x\in
M}\|d_xf_i^{-1}\|\bigr)^{-1}.
$$
We assume that the maps are non-degenerate (this is just for
simplicity of arguments) and strictly contracting, so that
$0<s_i^-\le s_i^+<1$.

Let $\eta=\max\limits_{x\in M}\#\{i\,|\,x=f_i(y_i)\text{ for some
}y_i\in \mathcal{Y}\}$ be the maximal multiplicity of overlaps and
$\varkappa_i=\max\limits_{x\in
\mathcal{Y}}\#\{y\in\mathcal{Y}\,|\,x=f_i(y)\}$ be the maximal
multiplicity of self-overlaps (we assume it is finite) on the
attractor. Denote also by $\vartheta_i$ the multiplicity of the
continuity partition for $f_i|_\mathcal{Y}$, i.e. the maximal
number of continuity domains intersecting the attractor and
meeting at one point of it.

\begin{theorem}\po
Let $\underline{D}=\alpha$, $\overline{D}=\beta$ be the solutions
of the equations
$$
\sum_{i=1}^N\tfrac1{\varkappa_i}|s_i^-|^{\alpha}=\eta,\qquad
\sum_{i=1}^N\vartheta_i|s_i^+|^{\beta}=1.
$$
Then the Hausdorff dimension of the attractor satisfies:
$$
\underline{D}\le \dim_\text{H}(\mathcal{Y})\le\overline{D}.
$$
\end{theorem}

In addition to Hausdorff dimension we will need some other
dimensional characteristics (see \cite{P2} for details). Denote by
$\mathfrak{N}(X,\delta)$ the minimal cardinality of covers of $X$
by balls of radius $\delta$. Then the lower and upper box
dimensions are defined by the formula:
 $$
\underline{\dim}_B(X)=\mathop{\underline\lim}\limits_{\delta\to+0}
\log\mathfrak{N}(X,\delta)/\log\tfrac1\delta,\quad
\overline{\dim}_B(X)=\mathop{\overline\lim}\limits _{\delta\to+0}
\log\mathfrak{N}(X,\delta)/\log\tfrac1\delta.
 $$
When these quantities are equal, their value is also called
fractal dimension.

Consider a Borel probability measure $\mu\in\mathcal{M}(X)$ (an
SRB-measure on the attractor can be taken in the SOC-context if it
exists). The upper and lower pointwise dimensions are defined then
as
 $$
\underline{d}_\mu(x)=\mathop{\underline\lim}\limits_{\delta\to+0}
\log\mu(B(x,\delta))/\log\delta,\qquad
\overline{d}_\mu(x)=\mathop{\overline\lim}\limits _{\delta\to+0}
\log\mu(B(x,\delta))/\log\delta.
 $$
When they are equal and constant a.e. the measure $\mu$ is called
exact-dimensional. This is precisely the case, when
$\op{supp}\mu=X$ and we have equality in the general chain of
inequalities (together with (\ref{dim-ineq1}) below):
 \begin{equation}\label{dim-ineq}
\op{ess.}\op{inf}\underline{d}_\mu(x)\le \op{dim}_{\text{H}}(X)\le
\underline{\dim}_B(X)\le \overline{\dim}_B(X),
 \end{equation}
where by essential infimum we mean its upper bound taken over all
subsets $U\subset X$ of measure 1 (and similar for
$\op{ess.}\op{sup}\overline{d}_\mu(x)$). The last two inequalities
are known and the first one follows from the inequality
$\op{ess.}\op{inf}\underline{d}_\mu(x)\le
\op{dim}_{\text{H}}(\mu)$ (\cite{P2}), where
$\op{dim}_{\text{H}}(\mu)=\lim_{\delta\to+0}\op{inf}\{
\op{dim}_{\text{H}}(Z)\,|\,\mu(Z)>1-\delta\}$.

Other dimensional characteristics of the measure are defined
similarly and satisfy:
 \begin{equation}\label{dim-ineq1}
\underline{\dim}_B(\mu)\le \overline{\dim}_B(\mu)\le
\op{ess.}\op{sup}\overline{d}_\mu(x).
 \end{equation}
Note that $\overline{\dim}_B(\mu)\le \overline{\dim}_B(X)$, while
the quantities $\op{ess.}\op{sup}\overline{d}_\mu(x)$ and
$\overline{\dim}_B(X)$ are in general incomparable.

The known formulas for the Hausdorff and other dimensions are
generalizations of Moran's result (\cite{M,H}) and are based on
the Bowen's equation (using the idea of coding); in this case one
usually obtains exact-dimensionality \cite{P2}. In the SOC-context
coding becomes problematic in the presence of singularities
(unless the properties of the SRB-measure are clarified) and thus
we cannot easily establish exact-dimensionality or formula for the
dimension.

We prove instead the inequality of the theorem for all the various
dimensions from (\ref{dim-ineq}) and (\ref{dim-ineq1}), which we
denote just by $\dim(\mathcal{Y})$:
 $$
\underline{D}\le \dim(\mathcal{Y})\le\overline{D}.
 $$

\noindent
\begin{proof}
Let us consider at first the upper box dimension
$\overline{\dim}_B(\mathcal{Y})$. The function $\mathfrak{N}$
satisfies the inequalities:
$$
\tfrac1{\varkappa_i}\mathfrak{N}(\mathcal{Y},\delta/{s_i^-})\le
\mathfrak{N}(f_i(\mathcal{Y}),\delta)\le
\vartheta_i\mathfrak{N}(\mathcal{Y},\delta/{s_i^+}).
$$
The inequality from above is obtained as follows. Let
$\mathcal{S}=\{x_j\}$ be a $\delta$-spanning set, i.e. a
collection of points from $X$ with $U_\delta(\mathcal{S})=X$. Then
$f_i(\mathcal{S})=\{f_i(x_j)\}$ may fail to be a $\delta\cdot
s_i^+$-spanning set thanks to singularities. Whenever
$\delta\ll1$, every $\delta$-ball intersects maximally
$\vartheta_i$ domains of continuity for $f_i$ meeting
$\mathcal{Y}$. Then we need to add maximally $\vartheta_i$ points
for each ball $U_\delta(x_j)$ intersecting singularities. The
inequality from below is proved similarly.

Now we have:
$$
\mathcal{Y}=f_1(\mathcal{Y})\cup\dots\cup f_N(\mathcal{Y})
$$
and the same for $U_\delta$-neighborhoods. This implies:
\begin{equation}\label{780}
\mathfrak{N}(\mathcal{Y},\delta)\le
\sum_{i=1}^N\mathfrak{N}(f_i(\mathcal{Y}),\delta)\le
\sum_{i=1}^N\vartheta_i\mathfrak{N}(\mathcal{Y},\delta/s_i^+).
\end{equation}
Denote $\sigma(\delta)=\mathfrak{N}(\mathcal{Y},\delta)
\delta^{\overline{\dim}_B(\mathcal{Y})}$. This functions grows
sub-polynomially:
\begin{equation}\label{848}
\mathop{\overline\lim}\limits_{\delta\to+0}\frac{\log
\sigma(\delta)} {\log 1/\delta}=0.
\end{equation}

 \begin{lem}\po\label{l34}
Let $\l_i>1$ be some numbers and $p_i>0$ be some probabilities,
$\sum_{i=1}^Np_i=1$. Then (\ref{848}) implies:
 $$
\mathop{\underline\lim}\limits_{\delta\to+0}\frac{\sum_{i=1}^N
p_i\sigma(\lambda_i\delta)}{\sigma(\delta)}\le1.
 $$
 \end{lem}

\noindent
\begin{proof}
Suppose the lower limit is $>\kappa>1$. Then for every
sufficiently small $\delta$ there exists $i\in[1,N]$ such that
$\sigma(\lambda_i\delta)\ge\kappa\sigma(\delta)$.

Denote $\bar\lambda=\max_{1\le i\le N}\lambda_i$. Let
$C=\max_{\delta\in[1/\bar\lambda,1]}\sigma(\delta)$. Then:
$$
\sigma(\delta)\le\frac1\kappa\sigma(\lambda_{i_1}\delta)\le
\frac1{\kappa^2}\sigma(\lambda_{i_1}\lambda_{i_2}\delta)\le
\dots\le\frac1{\kappa^{s(\delta)}}
\sigma(\lambda_{i_1}\dots\lambda_{i_{s(\delta)}}\delta),
$$
where $s(\delta)$ is the first number such that
$\lambda_{i_1}\dots\lambda_{i_{s(\delta)}}\delta\in[1/\bar\lambda,1]$.
This number can be estimated as follows:
$s(\delta)\ge-\log\delta/\log\bar\lambda-1$, whence:
$$
\frac{\log\sigma(\delta)}{\log1/\delta}\le \frac{\log
C-s(\delta)\log\kappa}{\log1/\delta}\le
\frac{\log(C\kappa)}{\log1/\delta}-\frac{\log\kappa}{\log\bar\lambda}.
$$
Therefore $\overline\lim_{\delta\to+0}\frac{\log
\sigma(\delta)}{\log1/\delta}\le
-\frac{\log\kappa}{\log\bar\lambda}<0$ and we get a contradiction.
This proves the lemma.
\end{proof} \\

Now to obtain the inequality from above for
$\overline{\dim}_B(\mathcal{Y})$ divide (\ref{780}) by
$\mathfrak{N}(\mathcal{Y},\delta)$. Denoting
$\varpi=\sum_{i=1}^N\vartheta_i|s_i^+|^{\overline{\dim}_B(\mathcal{Y})}$,
$\lambda_i=1/{s_i^+}$ and
$p_i=\vartheta_i|s_i^+|^{\overline{\dim}_B(\mathcal{Y})}/\varpi$
we get:
 $$
\frac1\varpi\le\sum_{i=1}^Np_i\frac{\sigma(\lambda_i\delta)}
{\sigma(\delta)}.
 $$
Thus Lemma \ref{l34} implies that $1/\varpi\le 1$ or
 $$
1\le \sum_{i=1}^N\vartheta_i|s_i^+|^{\dim\mathcal{Y}}
 $$
and the first claim
$\overline{\dim}_B(\mathcal{Y})\le\overline{D}$ follows from the
contraction $|s_i^+|<1$. The same arguments show another statement
that $\overline{d}_\mu(x)\le\overline{D}$.

The inequality from below follows from
 $$
\mathfrak{N}(\mathcal{Y},\delta)\ge
\frac1{\eta}\sum_{i=1}^N\mathfrak{N}(f_i(\mathcal{Y}),\delta)\ge
\frac1{\eta}\sum_{i=1}^N\frac1{\varkappa_i}
\mathfrak{N}(\mathcal{Y},\delta/s_i^-),
 $$
which implies $\sum_{i=1}^N\tfrac1{\varkappa_i}|s_i^-|^
{\underline{\dim}_B(\mathcal{Y})}\le\eta$ and the same for the
lower pointwise dimension: $\underline{d}_\mu(x)\ge\underline{D}$
a.e.

In this case we should define
$$
\sigma(\delta)=\mathfrak{N}(\mathcal{Y},\delta)
\delta^{\underline{\dim}_B(\mathcal{Y})}\text{ or }
\sigma(\delta)=\bigr(\op{ess.}\op{inf}-\log\mu(B(x,\delta))\bigl)
\delta^{\underline{\dim}_B(\mu)}
$$
respectively and use
\begin{lem}\po\label{l34+}
Let $\l_i>1$ be some numbers and $p_i>0$ be some probabilities,
$\sum_{i=1}^Np_i=1$. Then:
$$
\mathop{\underline\lim}\limits_{\delta\to+0}\frac{\log
\sigma(\delta)} {\log 1/\delta}=0 \ \ \Longrightarrow\ \
\mathop{\overline\lim}\limits_{\delta\to+0}\frac{\sum_{i=1}^N
p_i\sigma(\lambda_i\delta)}{\sigma(\delta)}\ge1.
$$
\end{lem}
~\\
\noindent This is proved similarly to Lemma \ref{l34}. The
inequalities for the Hausdorff dimension follows now from
(\ref{dim-ineq}).
 \end{proof} \\
\abz

\noindent {\bf Acknowledgements.} We thank J. Schmeling and Y.
Pesin for several stimulating discussions and references. M.
Rypdal thanks E. Mj\o lhus for some useful questions and comments.
We are grateful to the organizers of the Clay Mathematics
Institute/MSRI Workshop on Recent Progress in Dynamics (2004),
where we finished the final stage of the paper.


\end{document}